\RequirePackage{rotating}
\documentclass[trsc,nonblindrev,11pt]{informs3} 

\usepackage{multicol}
\usepackage{comment}
\usepackage{subcaption}
\usepackage{float}

\usepackage{amsmath,hyperref}
\usepackage{amssymb}
\usepackage{booktabs}
\usepackage{amsfonts}
\usepackage{rotating} 
\usepackage{multirow}
\usepackage{latexsym}
\usepackage{graphicx}
\usepackage{tabularx}
\usepackage[graphicx]{realboxes}
\usepackage{bbm}
\usepackage{mathtools}
\usepackage{bm}
\usepackage[shortlabels]{enumitem}
\usepackage{algorithm}
\usepackage{tabu}
\usepackage{framed}
\usepackage{ifthen}
\usepackage[dvipsnames]{xcolor}
\usepackage[margin=1in]{geometry}
\usepackage[noend]{algpseudocode}
\usepackage{fontenc}

\OneAndAHalfSpacedXI 

\usepackage{natbib}
 \bibpunct[, ]{(}{)}{,}{a}{}{,}%
 \def\bibfont{\small}%
\TheoremsNumberedThrough     

\EquationsNumberedThrough    

\MANUSCRIPTNO{2} 

\begin{document}


\RUNAUTHOR{Cummings {\it et al.}}

\RUNTITLE{{\color{black}Transportation Pricing Alliance Design}}

\TITLE{
{\color{black}Multimodal Transportation Pricing Alliance Design: Large-Scale Optimization for Rapid Gains}}

\ARTICLEAUTHORS{%
\AUTHOR{Kayla Cummings\thanks{Massachusetts Institute of Technology}, Vikrant Vaze\thanks{Dartmouth College}, \"Ozlem Ergun\thanks{Northeastern University}, Cynthia Barnhart\footnotemark[1]}
} 

\ABSTRACT{%
Transit agencies have the opportunity to outsource certain services to  established Mobility-on-Demand (MOD) providers. Such alliances can improve service quality, coverage, and ridership; reduce public sector costs and vehicular emissions; and integrate the passenger experience. To amplify the effectiveness of such alliances, we develop a fare-setting model that jointly optimizes fares and discounts across a multimodal network. We capture commuters' travel decisions with a discrete choice model, resulting in a large-scale, mixed-integer, non-convex optimization problem. To solve this challenging problem, we develop a two-stage decomposition with the pricing decisions in the first stage and a mixed-integer linear optimization of fare discounts and passengers' travel decisions in the second stage. To solve the decomposition, we develop a new solution approach combining tailored coordinate descent, parsimonious second-stage evaluations, and interpolations using special ordered sets. This approach, enhanced by acceleration techniques based on slanted traversal, randomization and warm-start, significantly outperforms algorithmic benchmarks. Different alliance priorities result in qualitatively different fare designs: flat fares decrease the total vehicle-miles traveled, while geographically-informed discounts improve passenger happiness. The model responds appropriately to equity-oriented and passenger-centric priorities, improving system utilization and lowering prices for low-income and long-distance commuters. Our profit allocation mechanism improves outcomes for both types of operators, thus incentivizing profit-oriented MOD operators to adopt transit priorities.
}%


\maketitle

%
\vspace{-20pt}

\section{Introduction}
Cities face critical challenges in the quest to improve urban mobility. Prior to the pandemic, congestion was steadily rising, translating to \$160 billion in annual costs to U.S. cities and record-breaking contributions to greenhouse gas emissions \citep{schrank2015}. Recent declines in transit ridership demonstrate the inability of transit's static infrastructure to accommodate rapidly evolving commuting patterns \citep{economist}. Private ride-sharing apps from Transportation Network Companies (TNCs) like Uber and Lyft have challenged this fixed-infrastructure status quo. TNCs transported 2.6 billion passengers in 2017, more than doubling the ride-sharing market since 2012 \citep{s2018}. The majority of urban TNC patrons admit that they would have otherwise walked, biked, taken public transit, or not made the trip, coinciding with tens of millions in annual transit revenue losses, worsening congestion, higher emissions, lower navigability of cities, and reduced accessibility to affordable public transportation options \citep{gr2018, s2018}. 

Mobility-on-demand (MOD) services have the potential to service {\it transit deserts}---low-density areas disconnected from public transit. However, cost presents a key barrier: while all public transit modes operate at a loss, MOD services administered by transit agencies incur the highest average per-trip costs (\$23.10 vs. the next-highest \$11.19 for commuter rail) \citep{ktp2016}. High labor needs, outdated technology, and coordination difficulties lead to inefficient, expensive operations. Notably, average TNC trip costs \$13, a full \$10 less than agency-sponsored MOD trips. Outsourcing all 223 million on-demand transit trips to TNCs could hypothetically save billions of dollars for US transit agencies \citep{ktp2016}. Thus, pricing alliances between TNCs and transit agencies have the potential to improve service quality and coverage, while reducing costs and decreasing citywide vehicle-miles travelled (VMT).

\subsection{Pricing Alliances in the Real World}\label{S:realworld}
TNCs have the infrastructure to provide more cost-effective MOD services supplementing fixed-route transit. Microtransit platforms like BRIDJ and Via---differentiated from TNCs due to fleets comprising minivans or shuttles as opposed to sedans---have oriented their business model toward complementing transit. 
The Federal Transit Administration's MOD Sandbox program has provided millions in funding to transit agencies in cities like Dallas, San Francisco, and Los Angeles for on-demand pilots filling service gaps in their service regions.
Rather than incurring the high fixed costs of designing a complementary MOD system from scratch, transit agencies could outsource MOD services to TNCs that are well-established, highly connected, and widely trusted; or to microtransit providers that more closely align with transit agency goals. This section formalizes such alliances within a rigorous conceptual framework. We define a {\it pricing alliance} as a cooperative pricing scheme between a transit agency and an MOD operator with independently operated infrastructures serving overlapping or adjacent regions, to improve each operator's own objectives while also improving system-wide benefits through integration. Indeed, real-world pricing alliance pilots have shown great promise in improving service quality, affordability and ridership. 

Pricing alliances are characterized by the intended service populations and the relationship of the MOD operator's system to the fixed-route transit network. Service population may be a targeted demographics (e.g., persons with limited mobility, low-income people, or senior citizens), or residents of a particular geographic area (e.g., a transit desert or the area near a transit hub). The MOD infrastructure may complement, substitute, or extend fixed-route options. After establishing the nature of the alliance, a joint pricing scheme is selected. Well-designed fares incentivize passenger travel choices that benefit the entire system.  Table \ref{T:survey} surveys the following key traits of the recent pricing alliances:
\begin{itemize}
    \item[--] {\it MOD operator}: This can be a TNC like Uber or Lyft, or a microtransit provider like Via.
    \item[--] {\it Service population}: While most alliances are open to everyone, some recent alliances have specifically strived to serve people with limited mobility, seniors, or essential workers.
    \item[--] {\it Fare structure}: Some alliances charge a flat fare and/or a variable fare based on distance traveled. Selectively applied, interpretable discount structures for jointly offered routes can encourage multimodal travel and engineer outcomes desired by the operators.
    \item[--] {\it Route structure}: Alliances often require specific trip geography: point-to-point (PTP) (trips occur within a given geographic region), zone-based (partitions a larger region into small zones and requires intra-zonal trips), and hub-based (at least one trip endpoint is anchored at specified locations). Figure \ref{F:zoning} in Appendix \ref{A:real-alliances} shows zone- and hub-based route structures.
    \item[--] {\it Integration into transit network:} MOD can be complementary (another mode option to improve service), substitutive (replaces existing fixed-route transit), first-/last-mile (FLM) (connects travelers to the fixed-route network), and extension (serves transit desert regions).
\end{itemize}
\begin{table}[htbp!]
\centering
\caption{Recent pricing alliances. \scriptsize{PTP: point-to-point. Flat fare: same price for all passengers. Mode: fares vary by travel mode. Distance: fares increase with distance. \citep{dart, rtc, seattle, theride, newmo, njvia, stlouis, marta, indygo, tfl, fflecsi, loopmunster, bridja, bridjs} $^*$ Not operated by a transit agency. $^{**}$ Specially marketed for senior citizens. $^{***}$ Available to everyone, but only during transit closures. $^{****}$ Transported essential workers at the beginning of the COVID-19 pandemic. {\color{black} $^{*****}$ 1 Euro surcharge on top of regular bus tickets is charged.}}}
\label{T:survey}
\resizebox{\textwidth}{!}{%
\begin{tabular}{llllllll}
\toprule[1pt]
\bf Program             & \bf City            &\bf Transit agency    & \bf MOD Op. &\bf Service population      & \bf Fares                    &\bf Routes       &\bf Integration  \\ \hline
\textit{GoLink}             & Dallas, TX      & DART              & Via, Uber                & Everyone          & Mode & Zone         & Complementary, FLM   \\ 
\textit{RTC On-demand Pilot}  & Las Vegas, NV   & RTC               & Lyft               & Paratransit       & Distance       & PTP & Substitutive   \\ 
\textit{Via to Transit}      & Seattle, WA     & King County Metro & Via                & Everyone          & Flat fare                & Hub          & Complementary, FLM   \\ 
\textit{The RIDE Flex}      & Boston, MA      & MBTA              & Uber, Lyft         & Paratransit       & Flat fare                & PTP & Substitutive   \\ 
\textit{NewMo} Pilot             & Newton, MA      &  City of Newton$^*$                 & Via                & Everyone, seniors$^{**}$ & Flat fare                & Hub          & Extension, FLM   \\ 
No program name   &   Jersey City, NJ      &    NJ TRANSIT               &      Via              &       Everyone            &       Distance                  &    Zone          &        Complementary, extension      \\ 
No program name &  St. Louis, MO     &     St. Louis Metro            &   Lyft            &   Everyone              &  Distance                    &    Hub         &    Complementary, FLM          \\ 
\textit{MARTAConnect}  & Atlanta, GA       &  MARTA             &   Uber, Lyft            &  Everyone (closures)$^{***}$               &   Distance                     &     PTP        &    Extension          \\ 
\textit{IndyGo + Uber} &   Indianapolis, IN    &  IndyGo               &  Uber             &    Essential workers$^{****}$             &      Flat fare                  &    PTP         &   Substitutive           \\ 
{\color{black}\textit{LOOPMuenster}} & {\color{black}Muenster, Germany} & {\color{black}Stadtwerke Muenster} & {\color{black}door2door} & {\color{black}Everyone} & {\color{black}Flat$^{*****}$} & {\color{black}Zone} & {\color{black}Complementary} \\ 
{\color{black}\textit{fflecsi}} & {\color{black}Wales, UK} & {\color{black}Transport for Wales} & {\color{black}Stagecoach} & {\color{black}Everyone} & {\color{black}Multiple} & {\color{black}Zone} & {\color{black}Extension} \\ 
{\color{black}\textit{GoSutton}} & {\color{black}Sutton, UK} & {\color{black}Transport for London} & {\color{black}ViaVan} & {\color{black}Everyone} & {\color{black}Flat} & {\color{black}PTP} & {\color{black}Complementary} \\
{\color{black}\textit{Slide Ealing}} & {\color{black}Ealing, UK} & {\color{black}Transport for London} & {\color{black}MOIA} & {\color{black}Everyone} & {\color{black}Flat} & {\color{black}PTP} & {\color{black}Complementary} \\
{\color{black}No program name} & {\color{black}Singapore} & {\color{black}Tower Transit Singapore} & {\color{black}BRIDJ} & {\color{black}Everyone} & {\color{black}Distance} & {\color{black}Hub} & {\color{black}Substitutive} \\
{\color{black}No program name} & {\color{black}Adelaide, Australia} & {\color{black}Torrens Transit} & {\color{black}BRIDJ} & {\color{black}Everyone} & {\color{black}Distance} & {\color{black}Hub} & {\color{black}Substitutive} \\\bottomrule[1pt]
\end{tabular}
}
\end{table}

{\color{black} As seen in Table \ref{T:survey}, real-world pricing alliances differ in terms of key attributes. Yet, two common patterns hold, reflecting a public sector twist on MOD operations. First, most real-world alliances use a dedicated vehicle fleet designated specifically for serving alliance passengers. Second, all real-world alliances have transparent, easy-to-interpret static prices. Critically, once set, fares are agnostic to variations in supply and demand---a key divergence from private-sector two-sided markets operated by the MODs. This transparent fare structure makes it easier for passengers to interpret prices and budget for their travels, marking the alliance as a mode of public transit.

Another component of alliance design is the establishment of a data-sharing policy to enable truthful and successful collaboration. The extent of operations coordination may range from totally separate (e.g., coupon codes to redeem within the MOD operator's privately operated marketplace) to fully integrated (a tailored mobile application for passengers to request rides specifically within the alliance), and anywhere in between (a separate tab or account type in the MOD operator's mobile application). Through these mechanisms, it is straightforward for the transit operator to have visibility into passenger requests and flows. It is also important to establish visibility into driver availability. With a dedicated driver fleet, this point is moot; however, when the MOD operator governs driver allocation from a pooled set of contractors that might drive in both public and private sectors, a guaranteed capacity must be regulated in an alliance contract. A dishonest representation of supply-side information could result in degraded service level and higher costs to the transit operator. We elaborate on how to enforce truthful information-sharing in Section \ref{S:profit_allocation}.

}

\subsection{Literature Review}

This work is at the intersection of FLM system design and operations, integrated multimodal transportation system design, and horizontal cooperation among competing transportation operators. 

\textbf{FLM system design and operations}: Research on demand responsive connector systems develops analytical models to evaluate service quality and determine first-mile system parameters, specifying optimal zone sizes and headways, transition points between regions best served by fixed-route vs. flexible services, and best practices for inter-zone transfer coordination {\color{black}\citep{lq2010,cq2013,ks2014,lsq2014,ls2017,kls2019, grahn2023optimizing}.
\cite{chopra2023mobility} develop an analytical model for system design and operations comparing fixed-schedule, on-demand, and hybrid services, focusing on the same mile problem.} The tactical question of how to operate a first-mile system is also well-studied. The Dial-A-Ride Problem (DARP) encompasses the vehicle routing problem faced by transit agencies, given a set of trip requests and a vehicle fleet \citep{hsklpt2018, mbc2017}. The Integrated DARP (IDARP) designs vehicle routes and schedules to meet trip requests, allowing transfers with fixed-route timetabled service \citep{pah2017}. Closely related to IDARP is the problem of carpooler matching and trip integration with transit timetables \citep{sasg2017}. Finally, many studies design strategies for routing and scheduling \citep{w2019}, pricing \citep{cw2018}, and trip request acceptance \citep{acl2019} for FLM systems.

{\color{black}\textbf{Integrated multimodal transportation system design}: Our work is related to the literature on optimal design and operation of transportation systems that leverages passenger choice models.} Past research has modeled decision-making travelers with preferences. One-to-one and many-to-one assignment problems among travelers and suppliers have been addressed with preference-based stable matchings to prevent participants from leaving ride-sharing systems \citep{wae2018} or transit systems \citep{rc2019}. Passenger decisions are often captured by discrete choice models. \cite{bny2020} jointly determine multimodal transit frequencies and prices to minimize wait times, subject to passenger choice. \cite{cvbm2017} optimize airline scheduling, fleet assignment, and fares while capturing the effects of competing high-speed rail service, taking passengers' mode choices into consideration. \cite{wvj2022} develop fixed-route transit timetables to maximize welfare, subject to competition with ride-sourcing companies, and congestion effects from passengers' mode switching. \cite{wjv2022} optimize a network of vertiports for supporting urban aerial mobility, with passenger mode choices described by two alternative models, including a multinomial logit model. \cite{bhls2021} tackle a welfare-maximizing system design and pricing problem for centrally coordinated multimodal transport networks with price-dependent demand, and formulate it using mixed-integer convex optimization. In contrast, we tackle a multi-objective alliance design problem with a practically suitable pricing scheme. It enables transparent price communication to passengers, but also prevents its convexification, thus heightening the computational challenge.

\textbf{Horizontal Cooperation}: Finally, we review cooperation models among competing operators. Literature on horizontal cooperation in logistics and airline scheduling is particularly mature \citep{cdf2007, gr2016, wgs2010, hcv2013}. \cite{cks2017} design a liner shipping alliance with endogenous linear demand for a homogeneous product; shipping companies first trade physical capacity on respective networks, and then compete to sell substitutable products in an overlapping market. {\color{black}Our work also involves joint products over a shared network subject to passenger choice}, but the allied operators jointly offer those products rather than exchanging capacity to compete. \cite{afls2019} formulate an urban transportation network flow game, using exogenous passenger and cost information to coordinate a single-fare payment among competing operators. \cite{bl2019i, bl2019ii} design mechanisms for the first-mile problem incorporating personalized passenger requirements. \cite{stz2021} investigate incentive mechanisms to inspire commuters to use public transportation, by modeling commuters, transit agency, ride-sharing platform, municipal government, and local private enterprises as stakeholders.

\cite{lc2022} use a Bayesian game to investigate data sharing between competing transit agencies to improve selfish outcomes in frequency setting, subject to user equilibrium passenger flows. {\color{black} They conclude that sharing more data does not necessarily improve service, and complementary services do not guarantee a grand coalition. They recommend that, in scenarios where truthful data-sharing is beneficial, government officials can regulate compulsory data-sharing.} The most similar study to ours in this branch of literature is by \cite{sl2022}, who formulate a transport choice game in which operators cooperatively price their pooled resources, subject to a multinomial logit model of passenger choice. They design a market share allocation rule that ensures a stable grand coalition. In contrast to our study, theirs assumes that each operator offers homogeneous products with a single price to travelers with unspecified origins and destinations, thus ignoring network effects. {\color{black}Recently, \cite{zhou2023impact} developed a Nash bargaining game between a metro operator and a ride-hailing service, and analyzed the solutions for a stylized linear city. In contrast, we focus on real-world large-scale pricing alliances. \cite{ding2023mechanism} develop an auction-based resource allocation and pricing framework to match passengers to transportation providers. Customers are asked to submit bids including the spatio-temporal attributes of their travel needs and their willingness to pay. Instead of auctions, our emphasis is on joint fare-setting and discount optimization which is directly in line with current real-world alliances. \cite{wang2023quantifying} estimate the benefits of integrating multiple rideshare platforms with each other.}

In summary, most existing studies individually model either operator or passenger incentives when designing integrated, multimodal urban transportation systems; to our knowledge, studies incorporating both strategic operators and passengers provide only general high-level intervention recommendations and rules of thumb. Our work differs in that we provide a prescriptive and strategic design framework to build pricing alliances at scale and in full operational detail.

\subsection{Contributions}
We propose a prescriptive pricing alliance to enable collaboration between transit agencies and ride-sharing operators. We formulate a fare-setting model to maximize total benefits across the integrated network. We demonstrate how our framework could help operators navigate competing alliance objectives: (1) enhancing access to high-quality transportation options for underserved populations, (2) lowering emissions and congestion from single-occupancy vehicle trips, and (3) maintaining the financial well-being of participating operators and incentivizing profit-oriented operators to participate. In doing so, a key technical challenge lies in capturing interdependencies between fares and passenger choice. Our model integrates a discrete choice model of passengers' mode and route decisions based on prices and non-pricing attributes like travel times.

From a technical standpoint, our fare-setting model is a large-scale, mixed-integer, non-convex optimization problem---a challenging class of problems. Our first technical contribution is to design a two-stage decomposition in which the first-stage pricing decisions parameterize second-stage fare discounts and passengers' travel choices. The second stage becomes a more tractable mixed-integer linear optimization problem that can be solved with commercial solvers. To solve the first-stage model, we develop a new solution approach combining tailored coordinate descent, parsimonious second-stage evaluations, and interpolations using Special Ordered Sets of type 2 (SOS2) \citep{sos2}. We also develop acceleration techniques based on slanted coordinate traversal and search direction randomization. This solution approach---our second technical contribution---is applicable to any two-stage formulation with a low-dimensional, convex, continuous first-stage solution space and any computationally expensive black-box second stage. This solution approach is found to significantly improve outcomes, for passengers and operators, compared to those obtained with state-of-the-art benchmarks based on Bayesian Optimization \citep{m2012}.

From a practical standpoint, when applied to a large-scale case study in the Greater Boston Area, our model sets fares in realistic ranges and with interpretable connections to alliance goals. An alliance with a greater focus on minimizing total VMT prefers flat rather than distance-varying fares to increase system utilization by long-distance commuters. On the other hand, alliances with a greater emphasis on increasing transit access set discounts with greater geographic variation to make alliance routes more attractive to heterogeneous populations. {\color{black}The clear alignment between operator goals and passenger choices achieved by our fare structures illustrates the value of modeling passenger choice.} Analysis of our results also shows that the model is appropriately responsive to equity-oriented objectives: it sets lower fares for, and increases utilization by, low-income and long-distance commuters. When compared to non-cooperative pricing, our fares and our tailored profit allocation mechanism together incentivize profit-oriented MOD operators not only to participate in the alliance but also to adopt the transit operator's priorities. {\color{black} Finally, in profit-oriented alliances, all operators benefit from collaborative pricing---even those excluded from the coalition.}

Section \ref{S:formulation} presents the allied fare-setting model formulation, its two-stage decomposition enabling tractable solutions, and our profit allocation mechanism. Section \ref{S:reformulation} describes our parsimonious SOS2-based coordinate descent approach. Its computational performance is presented in Section \ref{S:computational_analysis}. We present practical insights in Section \ref{S:results} and conclude in Section \ref{S:conclusions}. {\color{black}Appendices \ref{A:real-alliances}--\ref{app:practical} contain proofs, algorithmic details, additional computational experiments, and parameter calibrations.}
\section{Pricing Alliance Design Problem (PADP)}\label{S:formulation}

In the PADP, the alliance---i.e., the jointly acting operators---sets a fare structure to optimize joint operator priorities over the integrated multimodal network, subject to the passengers' mode and route choice decisions. The individual operators must then decide whether or not to participate in the alliance based on the optimized fares and a profit allocation mechanism.

\subsection{Problem Characteristics and Modeling Assumptions}\label{S:assumptions}

{\color{black}We specify our characterization of system-wide benefits, passenger decision-making model, 
need for a static interpretable fare structure, our resource allocation assumptions, and cost structure.}

\paragraph{System-wide benefits.}
The alliance cooperatively sets fares over an integrated network to maximize overall societal benefits to travelers, operators, and the rest of society \citep{d2012}. Consistent with the motivation of this work, we assume that transit agency's own objective is identical to that of the alliance. We characterize an operator's benefit as its profit and a passenger's benefit as its total utility across available travel options. High passenger utility corresponds to the availability of many high-quality travel options. There are many ways to capture the system's impact on the rest of society, defined as everyone except the travelers and operators. Most people taking alternative travel options choose to drive personal vehicles, creating negative externalities, such as air pollution. Because a pricing alliance involves no change in permanent infrastructure but rather better utilization of the existing infrastructure, its key benefits to the rest of society come from single-occupancy VMT reduction. We ultimately compute system-wide benefits as a weighted sum of operator profits, passenger utility, and a penalty for the outside-option VMT. The weights are determined by the alliance's relative priorities and can be varied to evaluate trade-offs.

\paragraph{Passenger choice model.}
We model travelers as rational agents making travel decisions according to a {\it multinomial logit} (MNL) choice model. MNL choice probabilities are proportional to each option's exponentiated utility, also known as its {\it attractiveness}.
MNL allows embedding a closed form of the passengers' decision-making process in the alliance fare-setting model, but it also presents limitations related to the independence of irrelevant alternatives (IIA) property. Some have circumvented such inaccuracies by using the {\it general attraction} model (GAM), of which MNL is a special case \citep{grs2015}. The GAM formulates choice probability as a function not only the attractiveness of the available options, but also the shadow attractiveness of unavailable options. In practice, researchers have set the shadow attractiveness values to zero, in the absence of reliable data to estimate these parameters \citep{wvj2020}. Others leverage the nested MNL,
of which the MNL is {\it also} a special case. 
In \cite{bny2020}, passengers select travel mode in the first level, and a route under that mode in the second level. We populate route choice sets with the fastest route in each travel mode (transit, MOD, or transit-MOD hybrid), including the option to drive, referred to in the literature as the \textit{outside option} or the \textit{no purchase alternative}. Thus, our choice model is equivalent to a GAM with zero shadow attractiveness, or a nested MNL with second-level choice sets containing one route each.

\paragraph{Fare-setting.}
{\color{black}Many MOD operators set time-varying fares on their independently operated networks.} In particular, TNCs implement fare multipliers to manage two-sided markets between drivers and riders \citep{ckw2017}. In a pricing alliance, however, the MOD operator is a contractor to the transit agency and consequently must agree to set time- and demand-homogeneous fares over the allied network. This facilitates transparent communication with passengers who can easily anticipate public sector prices, and also allows the transit operator to set a budget for the alliance with higher confidence. {\color{black} Further, we model binary discount options over route categories to allow the alliance to fine-tune passenger traffic toward optimizing system-wide benefits. Consistent with practice, we elect to use binary, rather than continuous, discounts to maintain an interpretable and equitable rules-based fare structure (e.g., ``30\% off of route types A, B, and C,'' vs. ``22.37\% off of route type A, 3.42\% off of B, 35.91\% off of C'')}.

\paragraph{Resource allocation.} {\color{black} We assume that the transit operator is capable of serving all demand redistribution resulting from the newly set fares, and that transit capacity reallocation is therefore unnecessary to consider in the pricing alliance design process. This assumption is consistent with practice in the context of most North American public transit systems which are currently suffering from low load factors. This phenomenon is especially severe in and near transit deserts, which in turn incentivizes the pricing alliance design proposed in this paper. \cite{sl2022} make a similar assumption: they assume a constant marginal cost due to market shares that do not change significantly. In summary, we assume that the load difference on the transit network when transitioning from non-cooperative to allied fare-setting will not impose a large enough change in network utilization to necessitate resource reallocation decisions. Our practical case studies (Section \ref{S:results}) validate this assumption, showing that the potential pricing alliances indeed do not pose a risk of over-saturating the transit infrastructure. On the MOD side, the alliance trips are served by a dedicated, right-sized vehicle fleet. Appendix \ref{A:wait_time} demonstrates how we right-size the fleet. It also models the relationship between travel times and passenger volumes for a given fleet size, and presents detailed results validating our key findings under volume-dependent travel times.}

\color{black} 
\paragraph{Cost structure.} 
Operating costs may include fixed infrastructure costs, variable labor costs, and variable transportation costs. Fixed infrastructure costs are not considered since capacity reallocation is not part of the PADP. Due to the fixed fleet size, variable labor costs depend only on the time horizon considered. We capture variable transportation cost as an operator-specific marginal cost per unit distance per passenger, based on vehicle fuel economy, fuel cost, and vehicle occupancy.

\color{black} 

\subsection{Exact Formulation}\label{S:exact_formulation}

We now provide notation for our allied fare-setting model. Passengers select from a set of routes, $\mathcal{R}$, serviced by a set of operators, $\mathcal{O}$, comprising a transit operator and an MOD operator, so that $|\mathcal{O}|=2$. Appendix \ref{S:coalition} relaxes this assumption by focusing on a case involving a second MOD operator. A {\it route} is a sequence of trip legs, each served by some operator's infrastructure. {\color{black}We capture flat and distance-based fares by defining the non-discounted static price of route $r \in \mathcal{R}$ as:
\begin{equation}
    \sum_{k \in \mathcal{O}_r} (\beta_k + \mu_k\Delta_{rk})
\end{equation}
where $\mathcal{O}_r \subseteq \mathcal{O}$ is the set of operators serving route $r$;  $\Delta_{rk}$ is the distance of route $r$ covered by operator $k$; and $\beta_k$ and $\mu_k$ are respectively the base fare and markup per unit distance traveled over operator $k$'s sub-network. Collectively, the base fares and markups of all operators are called the {\it fare parameters} $\bm{f}=(\bm{\beta}, \bm{\mu})$, which are decision variables in our model. Fare parameters are bounded above by values that are determined by local legislative or operational requirements (denoted by $B$ for base fares, $M$ for markups). Each operator $k \in \mathcal{O}$ incurs a marginal cost $c_k$ per unit distance traveled per passenger; operators will set markups accordingly, so that $\mu_k \geq c_k.$

Besides the fare parameters, the operators jointly select a set of discounted routes. Only a subset of routes, $\mathcal{R}^D \subseteq \mathcal{R}$, may be discount-eligible.} Rather than deciding whether or not each individual route should receive a discount, the discount-eligible routes may be grouped into {\it discount activation} categories. Routes in the same discount activation category may share common geographic components specified by the alliance. Grouping routes into categories allows passengers to easily interpret which routes are discounted from a map or a simple set of rules. Example definitions for discount activation categories include all routes anchored on a particular hub location, or all routes with origins and destinations contained in specified regions. {\color{black}Let $\mathcal{R}_a^D \subset \mathcal{R}$ be the set of routes in discount activation category $a \in \mathcal{A}.$ The sets $\mathcal{R}_a^D$ partition $\mathcal{R}^D,$ i.e., $\mathcal{R}^D=\cup_{a \in \mathcal{A}} \mathcal{R}_a^D$ and $\mathcal{R}_a^D \cap \mathcal{R}_b^D = \emptyset$ for $a \neq b \in \mathcal{A}$. Let $x_a \in \{0,1\}$ denote the decision variable activating discounts on all routes in $\mathcal{R}_a^D$. This assumption is not restrictive; the absence of discount activation categories is the same as putting each route in its own category: $|\mathcal{R}_a^D|=1, \forall a \in \mathcal{A}$. Note that relaxing this assumption might result in discount rules that are difficult to communicate to passengers in large-scale systems.

The routes selected for a discount receive a discount multiplier $\Lambda$, which has an allowable range of $[0, L]$ with $L \leq 1.$ The customer-facing price, $p_r$, of route $r \in \mathcal{R}$, is given as:
\begin{equation}\label{E:customerfacingprice}
    p_r = \begin{cases}
    (1-\Lambda x_a) \cdot \left(\sum_{k \in \mathcal{O}_r} (\beta_k + \mu_k\Delta_{rk})\right) &\text{if }\: r \in \mathcal{R}_a^D, a \in \mathcal{A}\\
    \sum_{k \in \mathcal{O}_r} (\beta_k + \mu_k\Delta_{rk}) &\text{if }\: r  \in \mathcal{R}\setminus \mathcal{R}^D
    \end{cases} 
\end{equation}
}
We consider a set $\mathcal{N}$ of passenger types. Each {\it passenger type} is a unique combination of origin, destination, and preference profile (i.e., their utility coefficients). $\mathcal{R}_i$ is the set of routes available to passenger type $i \in \mathcal{N}$. Some passengers are more averse to expensive travel options, whereas others are more sensitive to travel time, constituting different preference profiles. {\color{black}The number of passengers of each type depends on the time of day. There are $N_{it}$ passengers of type $i \in \mathcal{N}$ at time $t \in \mathcal{T}$. The utility to a passenger of type $i \in \mathcal{N}$ of route $r \in \mathcal{R}_i$ is $u_{irt} + \alpha_i p_r$, where ${u}_{irt}$ is from non-monetary route attributes (including travel time) at time $t \in \mathcal{T}$ and $\alpha_i < 0$ is the utility per unit price. The MNL-based market share of route $r \in \mathcal{R}_i$ for passenger type $i \in \mathcal{N}$ at time $t \in \mathcal{T}$ is:
\begin{equation}
    s_{irt} = \frac{e^{u_{irt} + \alpha_i p_r}}{e^{u_{i0t}} + \sum_{s \in \mathcal{R}_i} e^{u_{ist} + \alpha_i p_s}},
\end{equation}
where the outside option (not contained in set $\mathcal{R}$) has a utility $u_{i0t}$ and a market share given by:
\begin{equation}
    s_{i0t} = \frac{e^{u_{i0t}}}{e^{u_{i0t}} + \sum_{s \in \mathcal{R}_i} e^{u_{ist} + \alpha_i p_s}}.
\end{equation}
For now, we assume exogenous travel times to maintain tractability, but evaluate the effects of relaxing this assumption in Appendix \ref{A:wait_time}. Appendix \ref{A:wait_time} shows that under volume-dependent travel times, the pricing alliance design is minimally impacted, thus validating our assumption.
}

Finally, the operators' relative priorities over the system-wide performance metrics are captured by non-negative weights: $\pi^{PX}, \pi^{PR}, \pi^{VM}$, respectively, corresponding to passenger benefits, operator benefits, and the benefits from negative externality reduction.  Table \ref{T:notation} in Appendix \ref{App:notation} summarizes all notation. Model PADP-FS \eqref{O:obj}-\eqref{C:last} provides the exact formulation for the PADP fare-setting model. 
{\color{black}
    \begin{align}
           &(\text{PADP-FS}) \qquad \nonumber \\
           &\max_{\mathbf{p}, \mathbf{s}, \mathbf{x}, \bm{\beta}, \bm{\mu}, \Lambda} \quad \sum_{t \in \mathcal{T}} \sum_{i \in \mathcal{N}} N_{it}  \Bigg(\pi^{PX}  \big(u_{i0t} + \sum_{r \in \mathcal{R}_i} (u_{irt} + \alpha_i p_r) \big) +  \pi^{PR}  \sum_{r \in \mathcal{R}_i} \big(p_r - \sum_{k \in \mathcal{O}_r}c_k \Delta_{rk} \big)  s_{irt} - \pi^{VM} \Delta_i^0 s_{i0t} \Bigg) \label{O:obj} 
    \end{align}
\begin{align}
    \text{s.t. } \quad &p_r = \sum_{k \in \mathcal{O}_r} (\beta_k + \mu_k\Delta_{rk}) &&r \in \mathcal{R} \setminus \mathcal{R}^D \label{C:prices1} \\
    &p_r = (1 - \Lambda x_a) \cdot \left(\sum_{k \in \mathcal{O}_r} (\beta_k + \mu_k\Delta_{rk}) \right) &&a \in \mathcal{A}, r \in \mathcal{R}_a^D \label{C:prices2} \\
    &s_{irt} = \frac{e^{u_{irt} + \alpha_i p_r}}{e^{u_{i0t}} + \sum_{s \in \mathcal{R}_i}e^{u_{ist} + \alpha_i p_s} } &&i \in \mathcal{N}, r \in \mathcal{R}_i, t \in \mathcal{T} \label{C:badshare} \\
    &s_{i0t} = \frac{e^{u_{i0t}}}{e^{u_{i0t}} + \sum_{s \in \mathcal{R}_i}e^{u_{ist} + \alpha_i p_s} } &&i \in \mathcal{N}, t \in \mathcal{T} \label{C:badshare0}\\
    &0 \leq \beta_k \leq B,\, c_k\leq \mu_k \leq M &&k \in \mathcal{O} \label{C:first} \\    &0 \leq \Lambda  \leq L\label{C:lambda_last}\\
    &x_a \in \{0,1\} &&a \in \mathcal{A}\label{C:last}
    \end{align}
}
The PADP-FS jointly sets fares and discounts across the integrated network to maximize a weighted sum of passenger benefits, operator profits, and the negative of the vehicle miles traveled. Discounts are applied on selected routes (Constraints \eqref{C:prices1} and \eqref{C:prices2}) and utility-maximizing passengers make route selections according to an MNL model (Constraints \eqref{C:badshare} and \eqref{C:badshare0}). Fare parameters and the discount multipliers obey bounds (Constraints \eqref{C:first}-\eqref{C:lambda_last}). Discount activation decisions are binary (Constraints \eqref{C:last}). {\color{black} In summary, the PADP-FS model sets fares and activates discounts over the integrated network to improve system-wide outcomes, subject to fare-responsive passengers.}

\subsection{Two-stage Decomposition}\label{S:decomposition}

The PADP-FS model is a non-convex mixed-integer nonlinear optimization problem (MINLOP). There are no commercial solvers that accommodate non-convex MINLOPs, and no open-source solvers accept non-convex MINLOPs at practically large scale. Therefore, we propose a different solution approach. We decompose the formulation to tractably obtain high-quality solutions for practically sized problems (with tens of thousands of variables and hundreds of thousands of constraints in our case study). 
{\color{black} The key reason behind the success of the decomposition approach is its separation of the two main sources of model complexity, namely, non-convexity and the presence of discrete variables. The first stage of the decomposition directly handles the non-convexity but relegates all discrete variables to the second stage. By letting first-stage pricing decisions parameterize second-stage discount activations and passenger choice, the first stage solution space can be characterized as a small-dimensional, continuous, convex region.
The second stage can be formulated as a more tractable mixed integer linear optimization problem (MILOP). There exist many commercial solvers that provide certifiably optimal solutions to MILOPs of practically large scale.}

Let $\mathcal{F} :=  [0,B]^2 \times [c_k, M]^2$ and $\mathcal{L} := [0,L]$ be the domains of allowable fare parameters and discount multipliers. We parameterize the second-stage problem by $(\bm{\widehat{f}}, \widehat{\Lambda})=([\widehat{\beta}_k]_{k \in \mathcal{O}},[\widehat{\mu}_k]_{k \in \mathcal{O}}, \widehat{\Lambda}) \in \mathcal{F} \times \mathcal{L}$ with feasible region $\mathcal{S}(\widehat{\bm{f}}, \widehat{\Lambda})$. We reformulate the choice model constraints \eqref{C:badshare} and \eqref{C:badshare0} using sales-based linear programming \citep{grs2015}. 
\begin{align}
    &e^{u_{i0t}} \cdot s_{irt} = e^{u_{irt} + \alpha_i p_r} \cdot s_{i0t}  &&i \in \mathcal{N}, r \in \mathcal{R}_i, t \in \mathcal{T} \label{E:proportional}\\
    &s_{i0t} + \sum_{r \in \mathcal{R}_i} s_{irt} = 1 &&i \in \mathcal{N}, t \in \mathcal{T}\label{E:validprob} \\
    &s_{i0t} \geq 0 &&i \in \mathcal{N}, t \in \mathcal{T} \label{E:validprob2} \\
    &s_{irt} \geq 0 &&i \in \mathcal{N}, r \in \mathcal{R}_i, t \in \mathcal{T}\label{E:validprob3}
\end{align}
Equation \eqref{E:proportional} ensures that the market share of each route is proportional to its attractiveness. Constraints \eqref{E:validprob}, \eqref{E:validprob2}, and \eqref{E:validprob3} ensure that the market shares are non-negative and sum to 1. When the fare parameters and discount multipliers are fixed, Constraints \eqref{E:proportional} can be linearized by replacing $p_r$ with $\sum_{k \in \mathcal{O}_r} (\widehat{\beta}_k + \widehat{\mu}_k \Delta_{rk})$ for $r \in \mathcal{R} \setminus \mathcal{R}^D$  (see Constraints \ref{E:share_prod_triv}) and using big-$M$ inequalities for $r \in \mathcal{R}^D$. Let $\mathcal{N}_a \subset \mathcal{N}$ be the set of passenger types with at least one route option corresponding to discount activation category $a \in \mathcal{A}$, and define $M_{irt}^s = (e^{u_{i0t}})/(e^{u_{irt} + \alpha_i \sum_{k \in \mathcal{O}_r}(\widehat{\beta}_k + \widehat{\mu}_k \Delta_{rk})})  \geq 0$.
\begin{equation}
    \begin{cases}
    s_{i0t} \leq \frac{e^{u_{i0t}}}{e^{u_{irt} + \alpha_i \sum_{k \in \mathcal{O}_r}(\widehat{\beta}_k + \widehat{\mu}_k \Delta_{rk})}} s_{irt} \\
    s_{i0t} \geq \frac{e^{u_{i0t}}}{e^{u_{irt} + \alpha_i \sum_{k \in \mathcal{O}_r}(\widehat{\beta}_k + \widehat{\mu}_k \Delta_{rk})}} s_{irt} - M_{irt}^s x_{a} \\
    s_{i0t} \leq \frac{e^{u_{i0t}}}{e^{u_{irt} + \alpha_i (1 - \widehat{\Lambda}) \sum_{k \in \mathcal{O}_r}(\widehat{\beta}_k + \widehat{\mu}_k \Delta_{rk})}} s_{irt} + M_{irt}^s (1 - x_{a}) \\
    s_{i0t} \geq \frac{e^{u_{i0t}}}{e^{u_{irt} + \alpha_i (1 - \widehat{\Lambda}) \sum_{k \in \mathcal{O}_r}(\widehat{\beta}_k + \widehat{\mu}_k \Delta_{rk})}} s_{irt} 
    \end{cases} \qquad a \in \mathcal{A}, i \in \mathcal{N}_a, r \in \mathcal{R}_{i}\cap \mathcal{R}_a^D, t \in \mathcal{T}  \label{E:share_prod}
\end{equation}

To handle bilinearity in the profit terms in objective \eqref{O:obj}, {\color{black} we define a new decision variable $w_{irt} = \left(p_r - \sum_{k \in \mathcal{O}_r}c_k\Delta_{rk}\right) s_{irt}$. }Similarly as with equation \eqref{E:proportional}, the value of $w_{irt}$ can be set with linear constraints (see Constraints \ref{E:disc_prod_triv}) for $r \in \mathcal{R} \setminus \mathcal{R}^D$ by replacing $p_r$ with the full route price, and using Big-M inequalities (with $M_{ir}^w$ defined as $M_{ir}^w = \widehat{\Lambda}  \sum_{k \in \mathcal{O}_r} (\widehat\beta_k + \widehat\mu_k\Delta_{rk}) \geq 0$) for $r \in \mathcal{R}^D$).
\color{black}
\begin{equation}
    \begin{cases}
    w_{irt} \leq \sum_{k \in \mathcal{O}_r} (\widehat\beta_k + (\widehat\mu_k-c_k)\Delta_{rk})  s_{irt} \\
    w_{irt} \geq \sum_{k \in \mathcal{O}_r} (\widehat\beta_k + (\widehat\mu_k-c_k)\Delta_{rk})   s_{irt} - M_{ir}^w  x_{a} \\
    w_{irt} \leq  \sum_{k \in \mathcal{O}_r} \Big((1 - \widehat{\Lambda}) (\widehat\beta_k + \widehat\mu_k\Delta_{rk}) - c_k\Delta_{rk} \Big) s_{irt} + M_{ir}^w  (1 - x_{a}) \\
    w_{irt} \geq  \sum_{k \in \mathcal{O}_r} \Big((1 - \widehat{\Lambda}) (\widehat\beta_k + \widehat\mu_k\Delta_{rk}) - c_k\Delta_{rk} \Big)   s_{irt} 
    \end{cases} \qquad a \in \mathcal{A}, i \in \mathcal{N}_a, r \in \mathcal{R}_{i} \cap \mathcal{R}_{a}^D, t \in \mathcal{T} \label{E:disc_prod}
\end{equation}
\color{black}

In response to fare parameters set in the first stage, the second-stage problem activates discounts that optimize system-wide performance metrics, subject to passenger choice. The optimal value of the second-stage problem is denoted by $W$ in Equation \eqref{E:welfare} with $\mathcal{S}(\bm{\widehat{f}}, \widehat{\Lambda})$ defined as below.

{\small
\begin{equation}\label{E:welfare}
    W(\bm{\widehat{f}}, \widehat{\Lambda}) := \max_{(\bm{x}, \bm{s}, \bm{w}, \bm{p}) \in \mathcal{S}(\bm{\widehat{f}}, \widehat{\Lambda})}\sum_{t \in \mathcal{T}} \sum_{i \in \mathcal{N}} N_{it}  \left(\pi^{PX}  \big(u_{i0t} + \sum_{r \in \mathcal{R}_i} (u_{irt} + \alpha_i p_r) \big)  + \pi^{PR}  \sum_{r \in \mathcal{R}_i}  w_{irt} - \pi^{VM}  \Delta_i^0 s_{i0t} \right)
\end{equation}
}

\begin{align*}
  \mathcal{S}(\bm{\widehat{f}},\widehat{\Lambda}) \, \equiv \Big\{(\bm{x}, \bm{s}, \bm{w},  \bm{p}) \in \{0,1\}^{|\mathcal{A}|} \times \mathbb{R}_+^{|\mathcal{T}|\sum_{i \in \mathcal{N}}(|\mathcal{R}_i| + 1)} \times \mathbb{R}^{|\mathcal{T}|\sum_{i \in \mathcal{N}} |\mathcal{R}_i \cap \mathcal{R}^D|} \times  \mathbb{R}^{|\mathcal{R}|} \, : \, \qquad \qquad \qquad \qquad \qquad
  \end{align*}
  \vspace{-12mm}
  \begin{align}
  &\text{Constraints \eqref{E:share_prod} - \eqref{E:disc_prod}} \nonumber \\
    &s_{i0t} + \sum_{r \in \mathcal{R}_i} s_{irt} = 1 &&i \in \mathcal{N}, t \in \mathcal{T} \label{E:sum}\\
    &s_{i0t} = \frac{e^{u_{i0t}}}{e^{u_{irt} + \alpha_i \sum_{k \in \mathcal{O}_r}(\widehat{\beta}_k + \Delta_{rk} \widehat{\mu}_k)}}  s_{irt} &&i \in \mathcal{N}, r \in \mathcal{R}_i \setminus \mathcal{R}^D, t \in \mathcal{T} \label{E:share_prod_triv}\\
    &w_{irt} = \sum_{k \in \mathcal{O}_r} (\widehat\beta_k + (\widehat\mu_k-c_k)\Delta_{rk}) s_{irt} &&i \in \mathcal{N}, r \in \mathcal{R}_i \setminus \mathcal{R}^D, t \in \mathcal{T} \label{E:disc_prod_triv} \\
    &p_r = \sum_{k \in \mathcal{O}_r} (\widehat\beta_k + \widehat\mu_k\Delta_{rk}) &&r \in \mathcal{R} \setminus \mathcal{R}^D \label{E:price_nodisc}\\
    &p_r = (1 - \widehat{\Lambda} x_{a}) \cdot \left(\sum_{k \in \mathcal{O}_r} (\widehat\beta_k + \widehat\mu_k\Delta_{rk})\right)  &&a \in \mathcal{A}, r \in \mathcal{R}_a^D \label{E:price_disc} &&\Big\}
  \end{align}

Now we define the PADP-FS2SD model, the two-stage decomposition of the PADP-FS model.
    \begin{align}
        (\text{PADP-FS2SD}) \qquad \max_{\bm{f},\Lambda} \quad &W(\bm{f}, \Lambda)
        \quad \quad \text{s.t.} \quad \bm{f} \in \mathcal{F}, \Lambda \in \mathcal{L}
    \end{align}

\begin{lemma}\label{L:equiv}
Formulations PADP-FS and PADP-FS2SD are equivalent.
\end{lemma}

\begin{lemma}\label{L:feas}
Fix $(\widehat{\bm{f}}, \widehat{\Lambda}) \in \mathcal{F} \times \mathcal{L}.$ Then
$\mathcal{S}(\widehat{\bm{f}}, \widehat{\Lambda}) \neq \emptyset.$
\end{lemma}

Lemma \ref{L:equiv} shows that the two-stage decomposition is equivalent to the full formulation. Lemma \ref{L:feas} shows the second-stage feasibility for any first-stage solution (proofs in Appendices \ref{A:equivalent}-\ref{A:feas}).

\subsection{Profit Allocation Mechanism}\label{S:profit_allocation}

When considering a pricing alliance, an operator assesses whether the cooperative regime would improve their prioritized system-wide metrics over the non-cooperative regime. 
The MOD operator is solely profit-maximizing, while the transit agency may maximize a linear combination of multiple system-wide metrics. By entering a pricing alliance, the transit agency (denoted as TR) is guaranteed to fare no worse than that under the non-cooperative regime. However, the profit-maximizing MOD operator's participation depends on whether the alliance participation increases its profit.

We now design a profit allocation mechanism that guarantees the MOD operator's alliance participation. Let $\bm{f}^{nc}$ and $(\bm{f}^a, \Lambda^a)$ denote the non-cooperative equilibrium fare parameters and allied optimal fare parameters. Let $\Theta_k(\bm{f}^{nc})$ be the profit of operator $k$ in the non-cooperative regime, and $\Theta(\bm{f}^a, \Lambda^a; \bm{\pi})$ be the total profit of the alliance, for transit operator priority weights $\bm{\pi}$. 

\begin{lemma}\label{L:profit_allocation}
Let $\delta = \Theta(\bm{f}^a, \Lambda^a; \bm{\pi}) - \sum_{k \in \mathcal{O}} \Theta_k(\bm{f}^{nc})$ be the surplus allied profit compared to the total non-cooperative profit. Define $\Phi_k : \mathbb{R}^{|\mathcal{O}|}_+ \times \mathbb{R}_+ \to \mathbb{R}_+$ to be the profit allocation to operator $k \in \mathcal{O}$:
\begin{align}
    \Phi_{TR} ((\Theta_k(\bm{f}^{nc}))_{k \in \mathcal{O}}, \Theta(\bm{f}^a, \Lambda^a; \bm{\pi})) &= \begin{cases}
    \Theta_{TR}(\bm{f}^{nc}) + \frac{\delta}{|\mathcal{O}|} \cdot &\text{if } \delta \geq 0, \\
    \Theta(\bm{f}^a, \Lambda^a; \bm{\pi}) - \sum_{k \in \mathcal{O} \setminus \{TR\}}\Theta_k(\bm{f}^{nc}) &\text{otherwise,}
    \end{cases} \\
    \Phi_k ((\Theta_l(\bm{f}^{nc}))_{l \in \mathcal{O}}, \Theta(\bm{f}^a, \Lambda^a; \bm{\pi})) &= 
    \Theta_{k}(\bm{f}^{nc}) + \left(\tfrac{\delta}{|\mathcal{O}|}\right)^+, \qquad\forall k \in \mathcal{O}\setminus \{TR\}.
\end{align}
\begin{enumerate}[(a)]
    \item The MOD operator(s) will enter the pricing alliance with payment rule $\Phi$. 
    \item When $\delta \geq 0$, the mechanism satisfies Pareto efficiency, symmetry, the core property, scale invariance, and independence of irrelevant alternatives.
    \item {\color{black}  If $\pi^{PR} > 0, \pi^{PX} = 0,$ and $\pi^{VM} = 0,$ then $\delta \geq 0.$ If $\pi^{PR} = 0$ and $\pi^{PX} + \pi^{VM} > 0,$ with $c_k=0$ for each $k \in \mathcal{O}$, then $\delta < 0.$}
\end{enumerate}
\end{lemma}

The proof of Lemma \ref{L:profit_allocation} is in Appendix \ref{App:profit_allocation}.
Because the alliance's priorities may include benefits to passengers and/or to the rest of society (as reduced VMT), the alliance may earn less profit than the operators' combined profit in the non-cooperative regime. The payment rule in Lemma \ref{L:profit_allocation} assumes that the transit operator chooses to guarantee that the profit-oriented MOD operator earns at least as much as it would have earned outside of the alliance. We assume that the MOD operator will participate if its non-cooperative and allied profits are equal. The operators split the surplus evenly if the alliance accrues strictly more profit than that in the non-cooperative regime.

{\color{black} This raises the question of exactly when the alliance can generate surplus profit. Lemma \ref{L:profit_allocation}c shows that a purely profit-oriented alliance will naturally guarantee non-negative surplus, and a purely altruistic one will guarantee maximum reduction in profit compared to the corresponding non-cooperative setting. As profit is de-prioritized by the transit operator (and the alliance) relative to other objectives, the alliance will earn less profit because lower route prices improve passenger utility and reduce outside-option VMT. Consequentially, the MOD operators earn more profit in the non-cooperative regime as they increase their own route prices to monetize some passenger surplus into profits. Together, these two factors combine to shrink $\delta$ as $\pi^{PR}$ decreases and as $\pi^{VM}$ and/or $\pi^{PX}$ increase. In the absence of analytic expressions for $\Theta$ and $\Theta_k,$ the exact decision boundary between these two outcomes (namely, $\delta\geq 0$ and $\delta<0$) over the $\bm{\pi}$ space cannot be provided in closed form. However, we illustrate this transition empirically in Section \ref{S:value_added}, showing the value of cooperation in a real-world case study.

A key assumption in Lemma \ref{L:profit_allocation} is that the MOD operator truthfully shares information regarding system capacity---namely, a guaranteed fleet size. Such information is straightforward to share when the operators jointly establish a dedicated fleet of employee drivers to fulfill alliance service requests. However, in a two-sided TNC market, an MOD operator could potentially choose to withhold or misrepresent information about driver availability. Without regulation, they might surreptitiously allocate drivers to more profitable private demand pools. A smaller fleet inflates passenger travel times and leads to lower passenger volumes as well as lower total profit (see Appendix \ref{A:wait_time} for an experiment on endogenous travel times). The MOD operator may induce these dynamics under the guarantee that they will receive their non-cooperative profit under the profit allocation mechanism. Under certain data-sharing arrangements---including mutually managed mobile applications---the transit operator has direct oversight into available driver supply. Under others---including coupon vouchers in the MOD operator's marketplace---a dedicated driver fleet hired by the transit operator and managed by the MOD operator may be necessary to ensure truthful supply representation. In any case, \cite{lc2022} recommend that government officials should regulate compulsory data-sharing in transit coopetitions when it is beneficial. As transit agencies are closely associated with local governments, compulsory data-sharing is possible but should be enforced as a last resort. Considerable care should be taken to establish data-sharing practices  that enable either direct management or visibility by all operators into the available fleet when pricing alliance contracts are established, to enable rapid gains. It will ensure that the resources allocated to public transit are protected, and that the profit allocation mechanism is sound. 
}

\color{black}
\section{Solution Approach}\label{S:reformulation}

Our two-stage decomposition ensures that the second-stage is a mixed-integer linear optimization problem, and the first-stage is a low-dimensional decision problem over a convex space with a non-closed-form objective function. Without an analytic closed form, the first-stage objective function's gradients are inaccessible, which eliminates gradient-based approaches. Bayesian Optimization is applicable and leveraged in recent urban transportation studies on MOD systems \citep{lbds2019}, but does not provide clear convergence criteria. The nonlinear interdependencies between stages eliminate Benders decomposition. Finally, incompatibility with the simpler centralized welfare-maximization problem structure also eliminates convexification strategies \citep{bhls2021}.

Low-dimensional, convex first-stage space is amenable to {\it coordinate descent}, which takes turns fixing all fare parameters except one and greedily optimizing along the free dimension. But a one-dimensional search is also difficult for a search space comprising a continuous spectrum of optimal MILOP solutions. While the second-stage problem solves fast enough to be a useful tool (see Section \ref{S:computational_analysis}: needing $<5$ seconds on average), it is slow enough to warrant a judicious selection of first-stage valuation points. Our tailored coordinate descent scans each search direction by solving an auxiliary model that approximates the objective function in that search direction with Special Ordered Sets of type 2 (SOS2) \citep{sos2}. It terminates when no improvements are found along any coordinate direction. This leads to a new procedure, \textit{SOS2 Coordinate Descent (SOS2-CD)}, discussed in Section \ref{S:SOS2_CD}. To further improve final solution quality, we develop three acceleration strategies that build upon SOS2-CD. Section \ref{S:SOS2_CD_variants} describes the overall algorithm, which leverages multiple SOS2-CD runs, acceleration strategies, initialization procedures, and time limits.


\subsection{SOS2 Coordinate Descent (SOS2-CD)}\label{S:SOS2_CD}

The SOS2-CD (Algorithm \ref{A:SOS2_CD}) is a new coordinate descent (CD) procedure. Similar to basic CD, it starts from an initial point (Step \ref{L:startingpoint}) and successively optimizes a multivariate function along a single dimension at a time, holding other dimensions fixed. Unlike basic CD, it produces the next candidate solution with an SOS2 interpolation procedure, which generates an auxiliary objective function optimized in the one-dimensional optimization step (Step \ref{L:interpolateLast}).
First, in Step \ref{L:searchDir}, the subroutine {\sc Search Directions} provides a comprehensive ordered list of search directions (including multidimensional and/or randomized---discussed in Section \ref{S:SOS2_CD_variants}). Then, in Step \ref{L:interpolateFirst}, the subroutine {\sc Generate Anchors} returns evenly spaced SOS2 anchors along the search direction. The default is to cycle through coordinate axes when $random$ and $multidim$ are both set to $false$. Appendix \ref{App:sos2} presents subroutines {\sc Search Directions} and {\sc Generate Anchors} in detail. Step \ref{L:interpolateLast} picks the solution that maximizes the auxiliary objective function over the given search direction by computing an optimal second-stage solution for each anchor by solving the corresponding second-stage models and using these anchor solutions for interpolation. Step \ref{L:update} computes the true value of $W$ at the interpolated candidate solution and updates the current solution if the objective has improved. The algorithm iterates until convergence. 

\begin{algorithm}
\caption{SOS2 Coordinate Descent for maximization of $W$}\label{A:SOS2_CD}
\begin{algorithmic}[1]

\State {\sc args} $\bm{y}^{(0)}:$ Initial solution in $\mathcal{Y}$; $D$: No. of SOS2 anchors; $random$: Boolean, whether to randomize search directions; $multidim$: Boolean, whether to use multidimensional slanted search

\Procedure{SOS2 Coordinate Descent}{$\bm{y}^{(0)}, D, random, multidim$}

\State $objPrev \leftarrow -\infty$; $objCur \leftarrow W(\bm{y}^{(0)})$; $k \leftarrow 0$ \label{L:startingpoint}

\While{$objCur - objPrev > \epsilon$}

\State $k \leftarrow k + 1;\, objPrev \leftarrow objCur;\, \bm{y}^{(k)} \leftarrow \bm{y}^{(k-1)}$

\For{ $i \in $ {\sc Search Directions} ($random$, $multidim$) } \label{L:searchDir}

\State $(\bar{\bm{y}}^d)_{d \in 1..D} \leftarrow $ {\sc Generate Anchors}$(\bm{y}^{(k)}, i, D)$ \label{L:interpolateFirst}

\State  Interpolate new solution: $(\bm{\lambda}, \bm{z}) = SOS2^*((\bar{\bm{y}}^d)_{d \in 1..D})$ and set $\bm{y}^* = \sum_{d =1}^D \lambda_d \bar{\bm{y}}^d$\label{L:interpolateLast}

\State $\bm{y}^{(k)} \leftarrow \bm{y}^*$ if $W(\bm{y}^*) > objCur$ else $\bm{y}^{(k)}$ \label{L:update}

\EndFor

\State $objCur \leftarrow W(\bm{y}^{(k)})$

\EndWhile

\State \Return $\bm{y}^{(k)}$

\EndProcedure
\end{algorithmic}
\end{algorithm}

The rest of this subsection presents Step \ref{L:interpolateLast} of Algorithm \ref{A:SOS2_CD}, which selects a candidate pivot solution with SOS2 interpolation. First, we introduce relevant notation.
Let $\mathcal{Y} := \mathcal{F} \times \mathcal{L}$ be the first-stage solution space, and let $\bm{y} := (\bm{f}, \Lambda)\in \mathcal{Y}$ be a first-stage solution. 
Define $\mathcal{S}(\bm{y})$ and $\mathcal{S}^*(\bm{y}) \subseteq \mathcal{S}(\bm{y})$, respectively, as the set of feasible and optimal second-stage decisions for a first-stage solution $\bm{y}$:
{\small
\[ \mathcal{S}^*(\bm{y})  := \arg\max_{(\bm{x}, \bm{s}, \bm{w}, \bm{p}) \in \mathcal{S}(\bm{y})} \, \sum_{t \in \mathcal{T}} \sum_{i \in \mathcal{N}} N_{it}  \left(\pi^{PX}  \big(u_{i0t} + \sum_{r \in \mathcal{R}_i} (u_{irt} + \alpha_i  p_r) \big) + \pi^{PR}  \sum_{r \in \mathcal{R}_i} w_{irt} - \pi^{VM}  \Delta_i^0 s_{i0t} \right). \]
}

\begin{figure}[htbp!]
    \centering
    \includegraphics[width=0.3\textwidth]{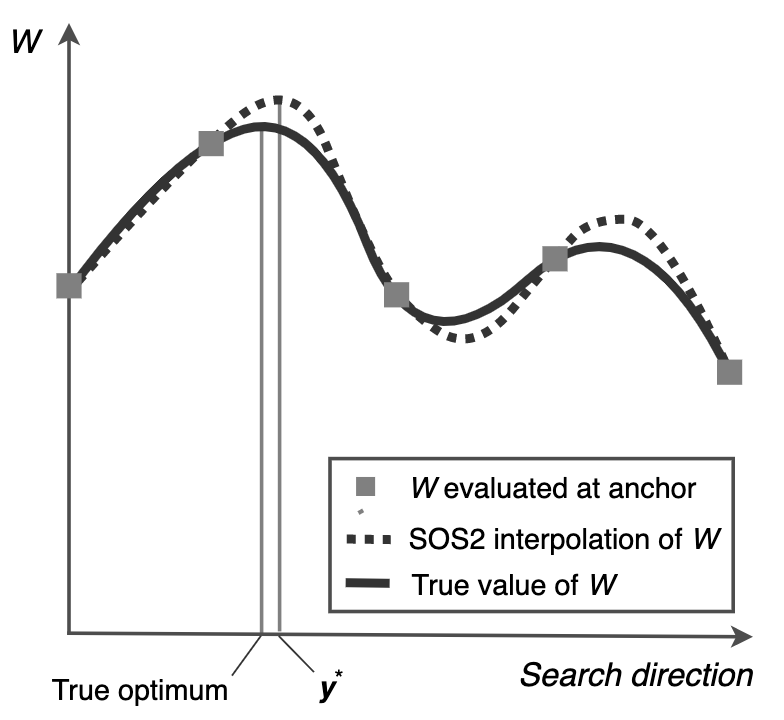}
    \caption{SOS2 interpolation of $W$ value and the selection of next candidate first-stage solution $\bm{y}^*$.}
    \label{F:sos2_schematic}
\end{figure}

Note that the number of anchors, $D$, is a design choice, navigating a trade-off between interpolation accuracy and computation speed. Larger $D$ improves accuracy, but reduces speed. For each $d \in \{1,\cdots,D\}$, denote the anchor as $\bar{\bm{y}}^d := (\bm{f}^d, \Lambda^d)$ and its corresponding second-stage solution as $(\bm{x}^d, \bm{s}^d, \bm{w}^d, \bm{p}^d) \in \mathcal{S}^*(\bar{\bm{y}}^d)$.
For a given $D$, the SOS2 model \citep{sos2} is specified as $ SOS2(D) \equiv \Big\{ \, (\bm{\lambda}, \bm{z}) \in \mathbb{R}_+^{D} \times \{0,1\}^{D-1} \, : \, \sum_{d = 1}^D \lambda_d = 1, \sum_{d =1}^{D-1} z_d = 1, \lambda_1 \leq z_1, \lambda_d \leq z_{d-1} + z_{d}\: \forall d \in \{2, \dots, D-1\}, \lambda_D \leq z_{D-1} \Big\}$. Here, the $\bm{\lambda}$ variables are the convex combination weights for anchors $(\bm{f}^d,\Lambda^d)$, and the binary $z_d$ variable indicates whether the segment between anchors $d$ and $d+1$ is selected for interpolation. The optimal SOS2 interpolation coefficients given a set of anchors are: 
\begin{align}
   SOS2^*&((\bar{\bm{y}}^d)_{d \in \{1,\cdots,D\}}) :=  \argmax_{(\bm{\lambda}, \bm{z}) \in SOS2(D)} \quad \sum_{t \in \mathcal{T}} \sum_{i \in \mathcal{N}} N_{it} \Big(\pi^{PX}  \big(u_{i0t} + \sum_{r \in \mathcal{R}_i} (u_{irt} + \alpha_i \sum_{d =1}^D p_r^d \lambda_d)  \big) +\nonumber \\
   &\qquad\qquad \pi^{PR} \sum_{r \in \mathcal{R}_i}\Big( \sum_{d =1}^D p_r^d  \lambda_d  - \sum_{k \in \mathcal{O}_r} c_k\Delta_{rk}\Big)\cdot\Big(\sum_{d =1}^D s_{irt}^d  \lambda_d\Big) - \pi^{VM}  \Delta_i^0 \sum_{d =1}^Ds_{i0t}^d \lambda_d \Big) \label{sos2_model}
\end{align}
Given $(\bm{\lambda}, \bm{z}) \in SOS2^*((\bar{\bm{y}}^d)_{d \in \{1,\cdots,D\}}),$ the approximately optimal solution in this search direction is interpolated as $\bm{y} = \sum_{d =1}^D \lambda_d \bar{\bm{y}}^d.$ Figure \ref{F:sos2_schematic} visualizes the selection of the next candidate solution using SOS2 variables. The objective function $W$ is evaluated at every anchor and approximated between anchors via interpolation. The next candidate solution is selected where the approximation of $W$ is maximized. As seen in Expression \eqref{sos2_model}, we interpolate the price and market share variables, $\bm{p}$ and $\bm{s}$, using piecewise linear functions, leading to a piecewise quadratic approximation for $W$. Since the objective function's nonlinearities are quadratic in nature due to the multiplicative revenue terms, we effectively capture them with this SOS2 approximation. Fortunately, it can be solved almost instantly to global optimality with commercial solvers, because $D$ is small by design.


\subsection{Final Algorithm}\label{S:SOS2_CD_variants}

We now present three acceleration strategies to escape local optima and improve solution quality. The first one relaxes the assumption of deterministic search direction order. Because the search direction ordering is arbitrary, we can randomize it at each iteration by setting the $random$ argument to {\tt TRUE}. Second, we exploit the fact that the SOS2 approximation is valid along \textit{any search direction} through the first-stage problem's solution space, not just those parallel to coordinate axes. Natural search direction candidates are those which jointly vary an operator's base fare and markup while holding other parameters constant. Given a pair of dimensions, this strategy randomly selects the spanning dimension, and then selects the line's slope in this 2D plane uniformly at random from the set of lines that intersect the current solution and span the selected dimensions. Finally, it drops anchors at evenly spaced points along the chosen line and obtains the next candidate solution maximizing the auxiliary objective function. When the $multidim$ argument is set to {\tt TRUE}, the list of search directions contains three items: (1) transit parameters, (2) MOD parameters, and (3) the discount multiplier. Whenever an operator's fare parameters are selected as the search direction, we sample a slanted line with the aforementioned procedure. By considering SOS2-CD over such slanted directions, we unlock directions that navigate trade-offs between high base fares and low markups vs. low base fares and high markups, which would be unavailable with single-coordinate search directions. Appendix \ref{App:sos2} fully specifies the subroutine {\sc Search Directions}.

Finally, we mitigate SOS2-CD's sensitivity to random initializations by leveraging warm-starts. The outcome of a single run of Algorithm \ref{A:SOS2_CD} depends on the initial solution $\bm{y}^{(0)}$. The overall algorithm (Algorithm \ref{A:SOS2_CD_augmented}) repeats SOS2-CD until a computational time budget is reached. Each run of SOS2-CD (Algorithm \ref{A:SOS2_CD}) is called a {\it trajectory.} The best fare parameters found across all trajectories are returned. Convergence to higher quality solutions is more likely with intelligent initializations. We warm-start the algorithm by first obtaining a few samples with a specified $warmStartProcedure$, and selecting the best starting points from them. The $warmStartProcedure$ might simply be uniform sampling from the solution space, or it can consist of searching the space in a more principled way, such as with Bayesian Optimization.

Algorithm \ref{A:SOS2_CD_augmented} presents the overall solution algorithm, with $\tau^{WS}$ and $\tau$ as the time limits on the warm-start and SOS2-CD procedures, respectively. The Boolean arguments $random$ and $multidim$ indicate whether randomized and/or multidimensional search directions are used. The procedure for generating informed initializations is specified by the argument $warmStartProcedure$.

\begin{algorithm}
\caption{Timed SOS2-CD with warm-start initialization}\label{A:SOS2_CD_augmented}
\begin{algorithmic}[1]

\State {\sc args} $\tau^{WS}$: Warm-start time limit (sec); $\tau$: SOS2-CD time limit (sec); $random$: Boolean, whether to randomize search directions; $multidim$: Boolean, whether to use multidimensional slanted search; $warmStartProcedure$: Initialization procedure; $D$: Number of SOS2 anchors

\Procedure{Timed SOS2-CD}{$\tau^{WS}$, $\tau$, $random$, $multidim$, $warmStartPocedure$, $D$}

\State $objCur \leftarrow - \infty$; $\mathcal{Y}^0 \leftarrow \emptyset$; $T^{WS} \leftarrow \tau^{WS}$; $T \gets \tau$; draw $\bm{y}^* \in \mathcal{Y}$ uniformly at random

\While{$T^{WS} > 0$} {\tt \quad // generate warm-start solutions }

\State Draw $\bm{y} \in \mathcal{Y}$ with $warmStartPocedure$

\State Subtract from $T^{WS}$ the time to run $warmStartPocedure$ and to compute $W(\bm{y})$

\State \textbf{if }{$T^{WS} \geq 0$} \textbf{then:}  Insert $\bm{y}$ into set $\mathcal{Y}^0$ 

\EndWhile

\While{$T > 0$} {\tt \quad // execute SOS2-CD}

\State \textbf{if }{$\mathcal{Y}^0\neq \emptyset$} \textbf{then:} $\bm{y}^{0} \leftarrow \arg\max_{\bm{y} \in \mathcal{Y}^0} \{W(\bm{y})\}$; Remove $\bm{y}^{(0)}$ from $\mathcal{Y}^0$
 
\State \textbf{else:} Draw $\bm{y}^{(0)}$ from $\mathcal{Y}$ uniformly at random

\State ${\bm{y}}\leftarrow$  {\sc SOS2 Coordinate Descent} $(\bm{y}^{(0)}, D, random, multidim)$

\State Subtract from $T$ the time to run {\sc SOS2 Coordinate Descent} and to compute $W({\bm{y}})$

\State \textbf{if }{$W({\bm{y}}) > objCur$ and $T \geq 0$} \textbf{then:} $\bm{y}^* \leftarrow {\bm{y}}$; $objCur \leftarrow W({\bm{y}})$ 

\EndWhile
\State \Return $\bm{y}^*$

\EndProcedure
\end{algorithmic}
\end{algorithm}
\vspace{-5mm}
\subsection{SOS2-CD and Multiobjective Optimization}
SOS2-CD procedure (Algorithm \ref{A:SOS2_CD}) assumes the second-stage optimization problem to be feasible over the entire first-stage decision space. Infeasibility of any anchor's associated subproblem would invalidate line \ref{L:interpolateFirst} of the procedure. Lemma \ref{L:feas} proves that $W$ in Equation \eqref{E:welfare} can be obtained for all $(\bm{f}, \Lambda) \in \mathcal{F} \times \mathcal{L}$. However, recall that the subproblem's objective function \eqref{E:welfare} is a linear scalarization of an inherently multi-criteria decision problem, optimizing operator profits, passenger happiness, and vehicle miles traveled. While it would be preferable to generate a Pareto-efficient curve of the model, the incorporation of an $\epsilon$-constraint would violate the anchor feasibility assumption. For example, an $\epsilon$-constraint that requires profit to exceed a certain percentage of the optimal profit is impossible to satisfy at a null anchor (zero fare parameters), which yields an infeasible anchor subproblem. To circumvent this limitation, we compile an approximated Pareto frontier from incumbent solutions generated through the SOS2-CD procedure. Appendix \ref{App:pareto}  depicts these outcomes and compares them with optimal solutions to the scalarized multiobjective formulation.
\color{black}
\section{Computational Results}\label{S:computational_analysis}

We now discuss the effectiveness of our new SOS2-CD approach through computational experiments, {\color{black} using a synthetic example in Section \ref{S:syntheticComp},} and then for a large-scale case study (see Appendix \ref{App:case_study} for details) of the Greater Boston Area in Sections \ref{S:1hour} and \ref{S:moreTime}. All optimization models are solved with Gurobi v9.0 and the JuMP package in Julia v1.4.

\color{black}
\subsection{Illustration of SOS2 Coordinate Descent over a Synthetic Example}\label{S:syntheticComp}

We begin with a small synthetic example, where a single passenger type chooses from three routes: a transit (TR), an on-demand (MOD), and a hybrid (HYB) route, and also an outside option (i.e., driving). Two allied operators jointly set distance-based markups to maximize total profit. They also decide whether to discount the hybrid route at 25\% (so that $\Lambda_{HYB} = 0.25$ and $\Lambda_{TR} = \Lambda_{MOD} = 0$). Each operator's marginal cost of transportation $c_k$ is set to 0. This example can be directly formulated in terms of three unknowns: the two markups $\bm{\mu}$ and the binary discount decision $x$.
\begin{align}\label{E:synthetic}
    \max_{\bm{\mu} \in [c_k, M ]^2,\, x \in \{0,1\}} \quad & \frac{\sum_{r \in \mathcal{R}}\left(\sum_{k \in \mathcal{O}_r}\left( (1-\Lambda_r x)\Delta_{rk} \mu_k -  \Delta_{rk}c_k\right)\right) \cdot e^{u_r + \alpha (1-\Lambda_r x) \sum_{k \in \mathcal{O}_r} \Delta_{rk} \mu_k} }{e^{u_0} + \sum_{s \in \mathcal{R}} e^{ u_s + \alpha (1-\Lambda_s x) \sum_{k \in \mathcal{O}_s} \Delta_{sk} \mu_k} }
\end{align}
Because of its simplicity, one can find the optimal solution in this example by brute force. We proceed by enumerating all markup pairs between $c_k=\$0$ and $M=\$5$ at a \$0.01 granularity, and evaluate the profit for each discount decision ($x=0$ or $x=1$). Table \ref{T:synthetic_attr} summarizes travel times (including wait times) and distance attributes. 
\begin{table}[htbp!]
    \centering
    \footnotesize
    \caption{Route attributes in synthetic example.}
    \begin{tabular}{ccccc}
    \toprule[1pt]
        Route & Travel Time  & MOD Distance  &  TR Distance   & Driving Distance  \\ 
        & (minutes) & (miles) & (miles) & (miles) \\\hline
        TR  & 75& 0&20 & 0\\
        MOD  & 55& 25& 0& 0\\
        HYB & 80 & 5&20 &0\\
        Outside & 40 & 0 & 0 & 25 \\\bottomrule[1pt]
    \end{tabular}
    \label{T:synthetic_attr}
\end{table}

\begin{figure}[htbp!]
    \centering
    \begin{subfigure}[t]{0.4\textwidth}
        \centering
        \includegraphics[width=\textwidth]{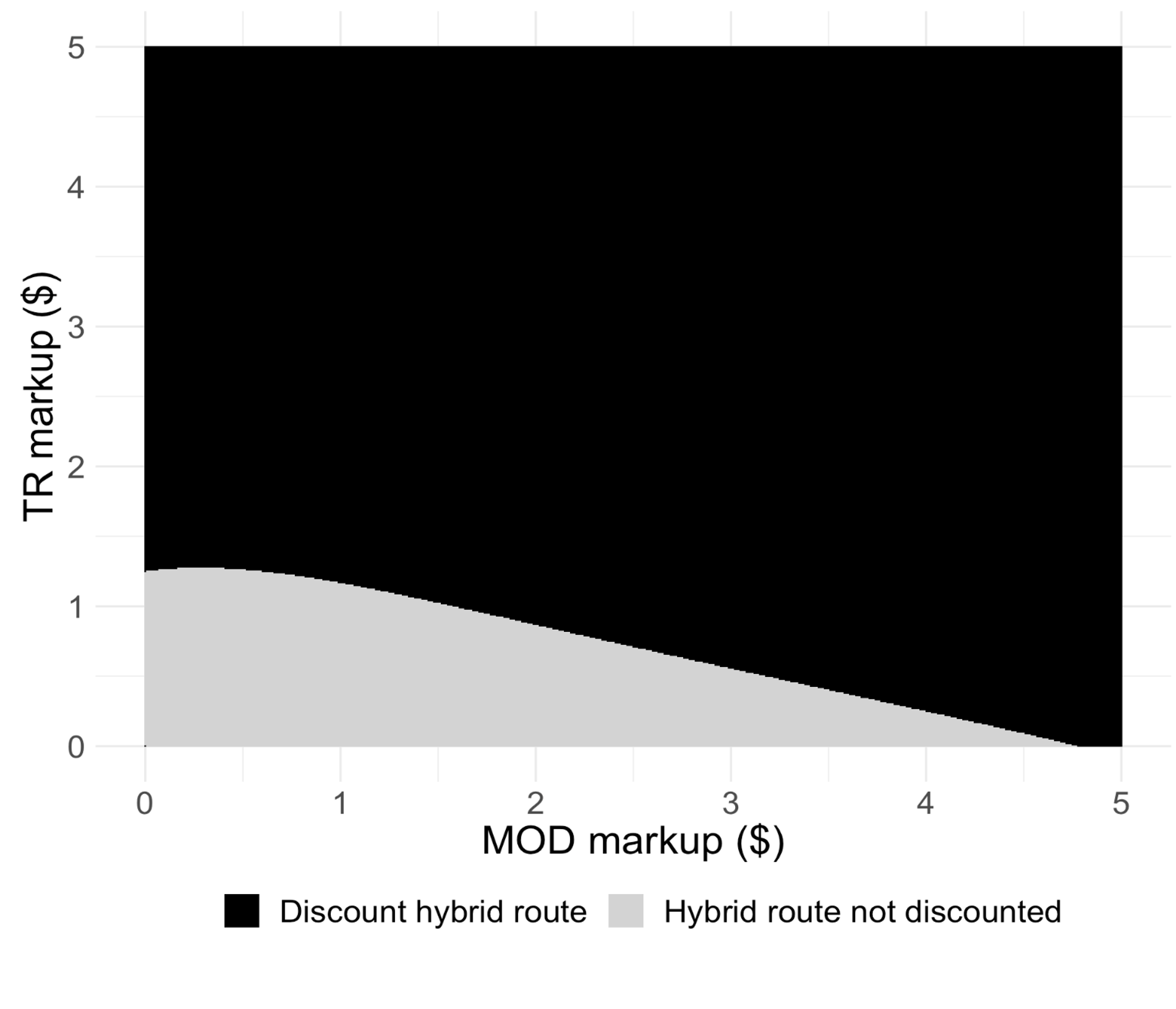}
        \caption{\footnotesize Optimal discount per markup pair.}
        \label{F:synthetic_discount}
    \end{subfigure}
    \begin{subfigure}[t]{0.4\textwidth}
        \centering
        \includegraphics[width=\textwidth]{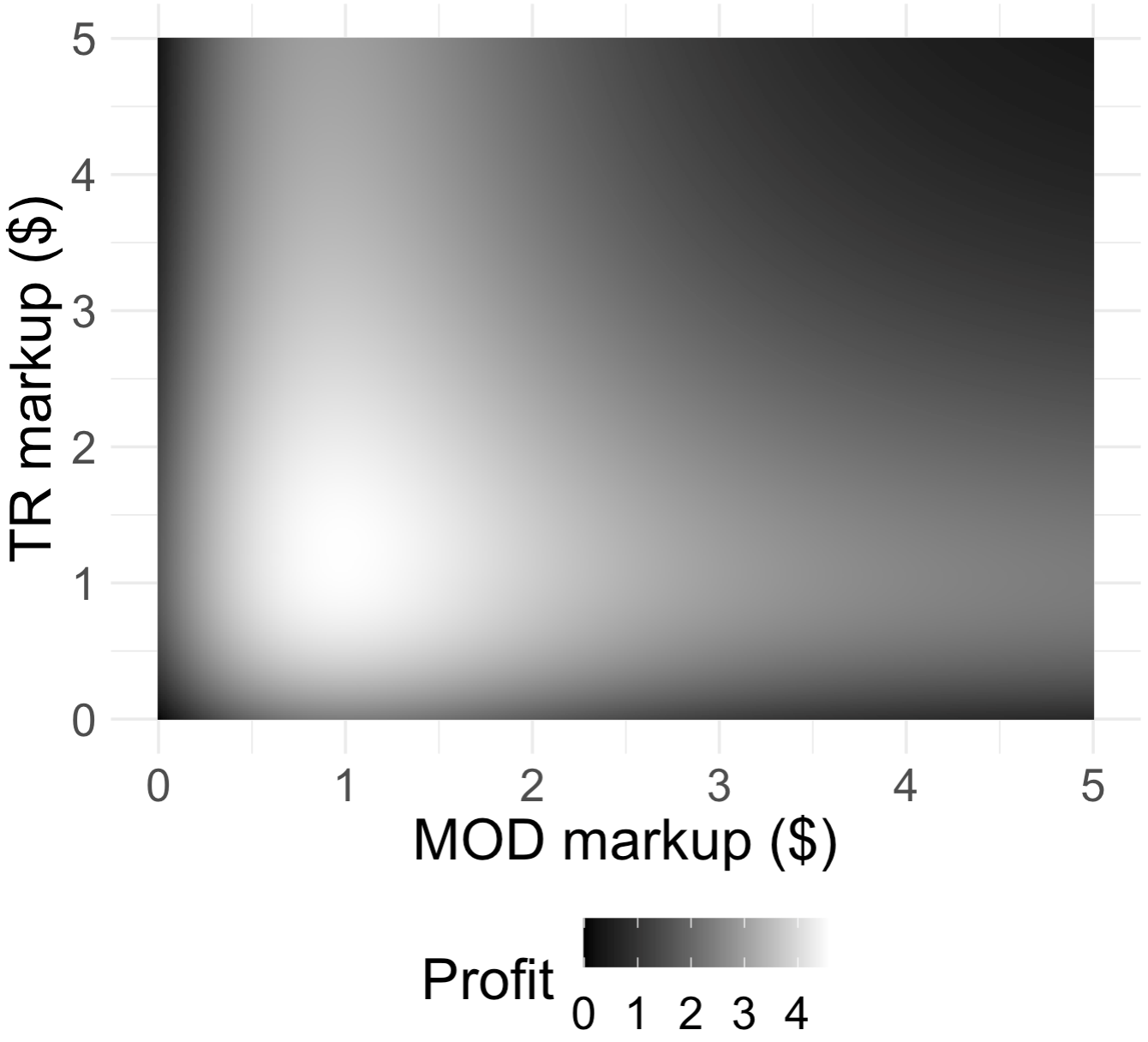}
        \caption{\footnotesize Total profit per markup pair.}
        \label{F:synthetic_objective}
    \end{subfigure}
    \caption{Feasible region characteristics for synthetic example \eqref{E:synthetic}.}
    \label{F:synthetic_heatmaps}
\end{figure}
\vspace{-5mm}
Figure \ref{F:synthetic_heatmaps} describes the example's feasible region, including the optimal discount decisions (Figure \ref{F:synthetic_discount}) and corresponding profit (Figure \ref{F:synthetic_objective}) for each markup combination. Figure \ref{F:synthetic_discount} shows that the hybrid route is discounted whenever the transit markup is moderate to high or both markups are high. The optimal total profit in Figure \ref{F:synthetic_objective} occurs when both operators' markups are moderately low. The optimal solution for the brute force problem is $\mu_{TR} = \$1.27$ and $\mu_{MOD} = \$0.99$ with $x = 1$. 

We will execute SOS2-CD over this feasible region and compare the outcomes with the brute force solution. For given markups $\widehat{\bm{\mu}}$, the second-stage objective function will be denoted as $W(\widehat{\bm{\mu}}).$ The full second-stage formulation is a binary linear optimization model in 11 variables. The first-stage problem is to select the profit maximizing markup pair.
\begin{align}
    \max_{\bm{\mu}} \quad &W(\bm{\mu}) \quad \quad \quad
    \text{s.t.} \quad \mu_k \in \left[c_k, M\right], \qquad \forall k \in \mathcal{O} \label{E:synth1s} 
\end{align}

\begin{figure}[htbp!]
    \centering
    \includegraphics[width=\textwidth]{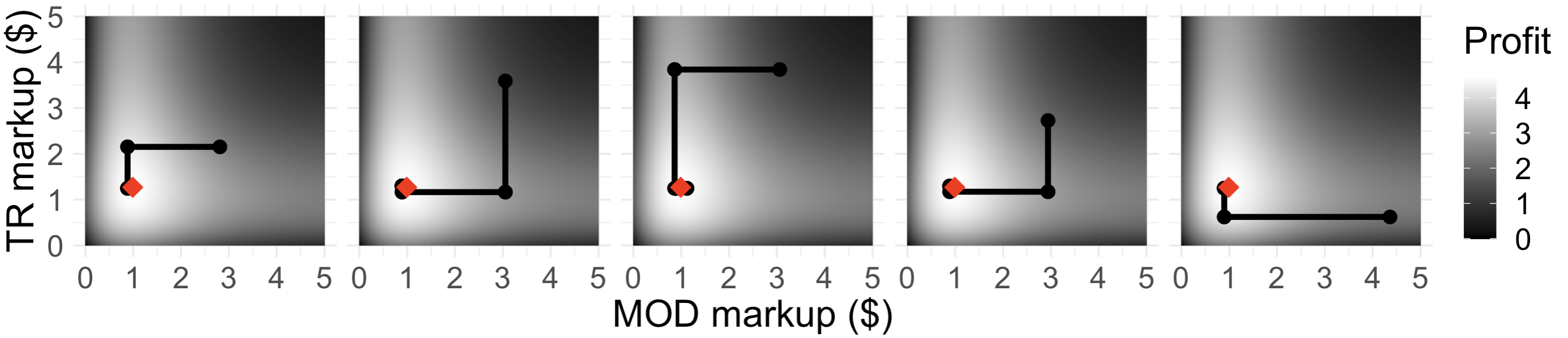}
    \caption{Five SOS2-CD trajectories for  synthetic example \eqref{E:synthetic}.}
    \label{F:synthetic_trajectories}
\end{figure}

We execute SOS2-CD to solve model \eqref{E:synth1s} with 100 random initializations. Figure \ref{F:synthetic_trajectories} depicts five of those SOS2-CD trajectories compared to the optimal solution. Across the 100 runs, on an average, SOS2-CD solution achieves 99.64\% of the brute force optimal profit, with a range of  99.51\%-99.85\%. In general, each trajectory efficiently navigates the search space starting from its random initial solution, always terminating within three iterations. These outcomes empirically confirm the quality of the SOS2-CD search process and the associated terminal solutions in this small synthetic example. Next, we move on to the full-scale case study.

\color{black} 
\subsection{Comparisons under 1-Hour Computational Time Budget}\label{S:1hour}

We now demonstrate the superior computational performance of our approach (Algorithm \ref{A:SOS2_CD_augmented}) over a realistic case study of an entire metropolitan area. Table \ref{T:ws} compares different versions of our approach, with multidimensional search (SOS2-CD-MD), randomized search directions (SOS2-CD-R), both (SOS2-CD-MD-R) and neither (SOS2-CD). None of these use intelligent warm-starts. They are compared to two algorithmic benchmarks---Brute-Force Coordinate Descent (BF-CD) and Bayesian Optimization (BO). BF-CD differs from SOS2-CD in how it conducts each iteration of coordinate descent. BF-CD uses a much higher number of ``anchors'' along the search direction, solves a second-stage model at each anchor, and it selects the anchor with the highest value of the second-stage objective function (instead of the SOS2-based interpolation) as the new candidate solution. It allows a drastic computation time increase for a more accurate evaluation of the points along a search direction. We implement BF-CD using a 1\% granularity for the discount multiplier and \$0.01 for base fares and markups. BO is a global optimization method for black-box functions (like our $W$ function) that are expensive to evaluate and may not have gradients \citep{m2012}. BO imposes upon $W$ a prior belief about the space of possible objective values based on the candidate solutions considered so far. The posterior distribution decides which candidate solution to evaluate next, so that a sequential search explores unseen regions in the decision space and exploits regions that are more likely to host global optima based on prior beliefs. See Appendix \ref{A:BO} for details.

We also tested time-limited SOS2-CD-MD-R with BO warm-starts, with varying time limit allocations to the warm-starts. Essentially, we execute Algorithm \ref{A:SOS2_CD_augmented} where the $warmStartProcedure$ is Bayesian Optimization. Warm-start trials have names ending in BO-$TL$, where $TL$ is the BO warm-start time limit in minutes. Table \ref{T:ws} presents results across 50 trials each with a 1-hour limit. All outcomes are expressed in surplus USD over the average performance of 1-hour BO benchmark.

First, we see that all four variations of our approach, even without warm-starts, significantly outperform the BO benchmark in terms of the average (by \$13.5K-\$19.3K) and best-case (by \$5.1K-\$5.8K) performance. Moreover, our approaches with multidimensional search (SOS2-CD-MD-R and SOS2-CD-MD), beat the BO benchmark even on the worst-case performance across the 50 trials (by \$7.5K-\$22.3K). Note that the average and worst-case performances of the approaches with either acceleration strategy (MD or R or both) were superior to those of the basic SOS2-CD approach. BF-CD never terminated within the one-hour time limit; in fact, it could not even evaluate one full set of anchors in all but 7 cases. Furthermore, our approaches with warm-starts perform even better than those without. In particular, a 40-minute BO warm-start drastically outperforms the benchmarks in the worst case and provides the best average-case performance, while a 20-minute BO warm-start has the strongest best-case performance. In summary, all our approaches significantly beat benchmarks, and all three acceleration strategies (random search, slanted search and warm-start) enhanced the performance of our basic SOS2-CD solution approach.

\begin{table}[htbp!]
\centering
\footnotesize
\caption{Objective function statistics with 1-hour time limits and 50 trials each, expressed in terms of surplus compared to the average performance of the 1-hour BO. $^*$ Average BO performance = \$3,634,074.}
\label{T:ws}
\resizebox{\textwidth}{!}{
\begin{tabular}{ccccccccc}
\toprule[1pt]
& & & & & & \multicolumn{3}{c}{Objective (Thousand \$)} \\
Algorithm & $random$ & $multidim$ & $warmStartProcedure$ & $\tau^{WS}$ & $\tau$ & Min & Avg. & Max  \\\hline
BO & - & - & - & - & - & $-28.0$ & $0.0^*$ & $18.6$  \\
SOS2-CD-MD-R-BO-$50$ & Yes & Yes & BO & 50 & 10 & $-11.9$ & $15.4$ & $24.1$\\
SOS2-CD-MD-R-BO-$40$ & Yes & Yes & BO & 40 & 20 & $\bm{12.6}$ & $\bm{20.3}$ & $24.1$\\
SOS2-CD-MD-R-BO-$30$ & Yes & Yes & BO & 30 & 30 & $0.0$ & $19.7$ & $24.0$\\
SOS2-CD-MD-R-BO-$20$ & Yes & Yes & BO & 20 & 40 & $-4.3$ & $19.1$ & $\bm{24.6}$ \\
SOS2-CD-MD-R-BO-$10$ & Yes & Yes & BO & 10 & 50 & $-9.8$ & $20.0 $ & $24.2$\\
SOS2-CD-MD-R & Yes & Yes & - & - & 60 & $-20.5$ &  $19.0$ & $23.7$ \\\hline
BF-CD & - & - & - & - & 60 & - & - & $4.8$ \\
SOS2-CD & No & No & - & - & 60 & $-103.8$ & $13.5$ & $24.1$ \\
SOS2-CD-MD & No & Yes & - & - & 60 &  ${-5.7}$ &  ${19.3}$ & $24.1$\\
SOS2-CD-R & Yes & No & - & - & 60 & $-77.6$ &  $17.1$ & ${24.4}$ \\ 
\bottomrule[1pt]
\end{tabular}}
\end{table}

\subsection{Comparisons under Higher Computational Time Budgets}\label{S:moreTime}

Table \ref{T:comp} compares performances under three time budgets---1 hour, 6 hours, and 12 hours, and provides statistics on the number of trajectories. All versions of our approach under all time budgets outperform the BO benchmark on average, and especially so for the versions with acceleration strategies. The larger time budgets allow accelerated SOS2-CD to offer robust performance. In general, the trajectories of SOS2-CD with multidimensional search (i.e., SOS2-CD-MD and SOS2-CD-MD-R) converge more quickly, allowing more trajectories to be computed in a given time limit. BF-CD is extremely slow and did not terminate before the 12-hour limit in any of our runs. We report the performance statistics for BF-CD corresponding to the best solutions found within the computational time budgets, prior to termination. While the best-case runs of BF-CD provide a slight edge over all benchmarks (of merely \$200 USD), the average and worst-case performance is significantly worse than our methods. BF-CD is more thorough for a single random initialization, but computationally too intensive to properly explore the search region. Note that warm-starts did not provide much additional value for longer time budgets and hence warm-start approaches are omitted from Table \ref{T:comp}. Appendix \ref{App:component_performance} further analyzes the solution times and SOS2 optimality gaps.

\begin{table}[htbp!]
\caption{Objective function statistics with varying time budgets and 50 trials each, expressed in terms of surplus compared to the average performance of the 1-hour BO. $^*$ Average BO performance = \$3,634,074.}
\label{T:comp}
\centering
\footnotesize
\begin{tabular}{llcccccc}
\toprule[1pt]
Time & Algorithm & \multicolumn{3}{c}{Trajectories} & \multicolumn{3}{c}{Objective (Thousand USD)} \\
limit &  & Min & Avg. & Max & Min & Avg. & Max \\ \hline
\multirow{6}{*}{1 hour} & BO  & - & - & - & $-28.0$ & $0^*$ & $18.6$  \\
& BF-CD & 0 & 0 & 0 & - & - & $4.8$ \\
& SOS2-CD & $2$ & $3.8$ & $5$ & $-103.8$ & $13.5$ & $24.1$ \\
& SOS2-CD-MD &  $2$ & $4.0$ & $7$ &  $\bm{-5.7}$ &  $\bm{19.3}$ & $24.1$\\
& SOS2-CD-R &  $2$ & $3.8$ & $5$ & $-77.6$ &  $17.1$ & $\bm{24.4}$ \\
& SOS2-CD-MD-R & $2$ & $4.2$ & $6$ & $-20.5$ &  $19.0$ & $23.7$ \\ \hline 
\multirow{6}{*}{6 hours} & BO &  - & - &- & $7.7$ & $17.9$ & $23.9$  \\
& BF-CD &  0 & 0 & 0 & $-2132.7$ & $-791.9$ & $\bm{24.6}$ \\
& SOS2-CD & $16$ & $20.1$ & $25$ & $-2.7$ & $22.3$ & $24.3$\\
& SOS2-CD-MD & $15$ & $21.3$ & $28$ & $22.0$ & $23.3$ & $24.1$\\
& SOS2-CD-R & $13$ & $19.7$ & $24$ & $\bm{22.4}$ & $\bm{23.8}$ & $24.4$\\
& SOS2-CD-MD-R & $15$ & $21.4$ & $26$ & $21.1$ & $23.3$ & $24.3$ \\\hline 
\multirow{6}{*}{12 hours} & BO &  - & - &- & $18.7$ & $21.9$  & $23.9$ \\
& BF-CD &  0 & 0 & 0 & $-2003.4$ & $-112.7$ & $\bm{24.6}$ \\
& SOS2-CD &  $32$ & $39.6$ & $50$ & $22.0$ & $23.7$ & $24.3$\\
& SOS2-CD-MD & $32$ & $42.4$ & $54$ & $22.8$ & $23.7$ & $24.2$\\
& SOS2-CD-R &  $26$ & $38.2$ & $47$ & $\bm{23.5}$ & $\bm{24.0}$ & $24.4$\\
& SOS2-CD-MD-R & $31$ & $42.7$ & $52$ & $22.2$ & $23.7$ & $24.4$\\
\bottomrule[1pt]
\end{tabular}
\end{table}
\section{Insights from Practical Case Study}\label{S:results}

To inform policymaking, we obtained practical results with the PADP model over a Greater Boston Area case study, for a potential pricing alliance between the Massachusetts Bay Transit Authority (MBTA) and a TNC like Uber or Lyft. MBTA subsidizes Uber and Lyft trips as part of their on-demand paratransit program called The RIDE Flex.
We model an alliance with a wider passenger scope aligned with MBTA goals outlined in a recent report \citep{mbta2019}, where MBTA identified 14 towns (called ``urban gateways'') adjacent to the commuter rail network whose residents had the greatest likelihood of utilizing---and benefiting from---targeted transit expansion efforts. We identify these 14 towns as the {\it service region} of the potential pricing alliance (see Figure \ref{F:pplace} in Appendix \ref{App:case_study}). Our case study integrates many datasets describing travel characteristics in the Greater Boston Area for a single time period---weekday morning commute at 6-10 am {\color{black} (see Appendix \ref{app:timePeriods} for an experiment assessing the impact of multiple time periods on optimal fares)}. We consider passenger travel patterns for those who commute from the service region to the inner city (Boston and Cambridge), or those who commute locally within the service region. See Appendix \ref{App:case_study} for a detailed presentation of the case study setup including data sources and data processing.

Section \ref{S:validation} confirms that our model yields interpretable outputs with prices in realistic ranges. An equity-oriented case study in Section \ref{S:equity} shows the value of accurately capturing passenger preferences. Section \ref{S:value_added} demonstrates the value of cooperative pricing, {\color{black}illustrating how transit agencies might negotiate alliance priorities with profit-oriented MOD operators. Appendix \ref{S:coalition} investigates coalition-forming behavior when there are multiple MOD operators.} All results over the Greater Boston Area case study are obtained with the SOS2-CD-MD-R approach and a 12-hour time limit.

\subsection{Model Validation}\label{S:validation}

\begin{figure}[htbp!]
\centering
\begin{subfigure}[b]{0.8\textwidth}
    \centering
    \includegraphics[width=\textwidth]{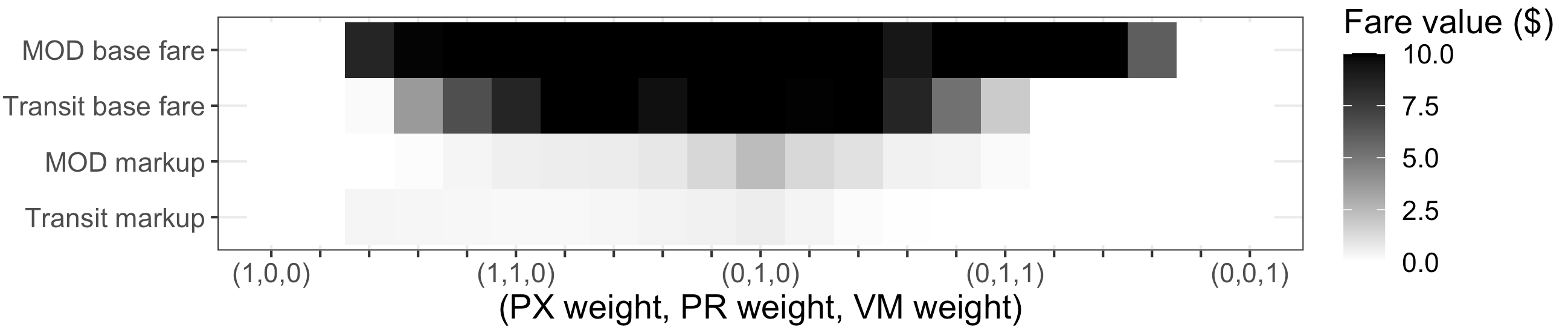}
    \caption{\footnotesize Optimal fare parameter values.}
    \label{F:fares}
\end{subfigure}
\begin{subfigure}[b]{0.8\textwidth}
    \centering
    \includegraphics[width=0.9\textwidth]{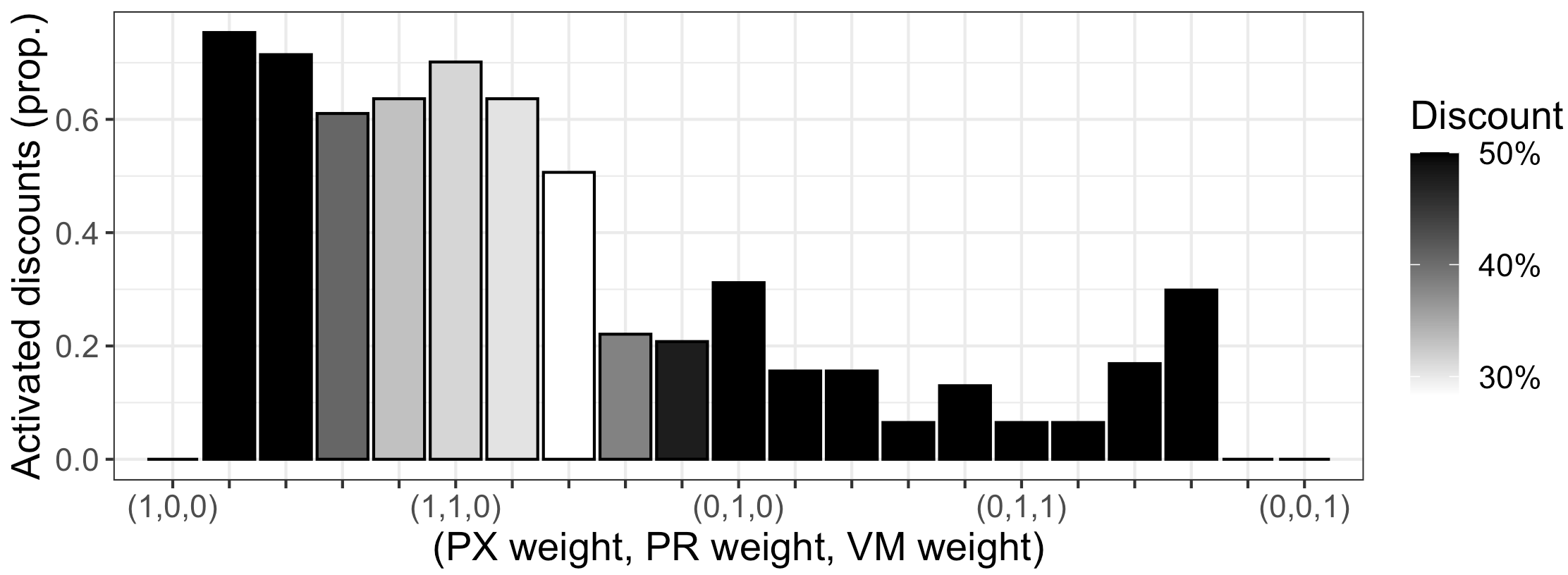}
    \caption{\footnotesize Optimal discount characteristics.}
    \label{F:discounts}
\end{subfigure}
\caption{\footnotesize Optimal fares across varying alliance priority regimes.}\label{F:testsuite_fares}
\end{figure}

We now demonstrate how the allied fare-setting model sets route prices in realistic and practically reasonable ranges. Further, we find that the optimal fares intuitively reflect various portfolios of alliance priorities. We vary the objective function coefficients $\bm{\pi}$, i.e., relative weights among the three performance metrics. In particular, we focus on regimes with varying combinations of priorities between profit and passenger utility (i.e., setting $\pi^{VM}=0$), as well as between profit and VMT (i.e., setting $\pi^{PX}=0$). We do not emphasize regimes that completely exclude profit as a priority, because they intuitively result in zero fares and are not interesting from an analysis standpoint. Thus, all experiments have $\pi^{PR}>0$. We also do not analyze regimes that vary all three weights (explained later in this section). Table \ref{T:testsuite} summarizes route prices, system utilization, and performance metrics across tested priority regimes. Figure \ref{F:testsuite_fares} depicts optimal fare parameters and discount multipliers, demonstrating the different fare-setting strategies across regimes.

\begin{table}[htbp!]
    \centering
    \footnotesize
        \caption{Aggregate metrics for different operator priority regimes. Performance metrics are normalized against best possible values. System utilization (util.) is the alliance's total market share, i.e., the percentage of travelers electing to travel on a transit, MOD, or hybrid option instead of driving a single-occupancy vehicle.}
    \label{T:testsuite}
    \begin{tabular}{ccc|ccc|ccc|c}
    \toprule[1pt]
    \multicolumn{3}{c|}{Objective weights} & \multicolumn{3}{c|}{Route price (\$)} & \multicolumn{3}{c|}{Performance metrics} & System  \\    
    $\pi^{PX}$ & $\pi^{PR}$ & $\pi^{VM}$ & Min. & Mean & Max.  & PX  & PR  & VM  & util. \%  \\\hline 
    1.0 & 0 & 0 & \$0.00 & \$0.00& \$0.00& 100.00\% &	0.00\% &	100.00\% &	50.32\%  \\
    1.0 & 0.2 & 0 & \$0.00 &	\$0.04 &	\$0.54 &99.26\% &	3.19\% &	100.59\% & 50.32\% \\
    1.0 & 0.4 & 0 & \$0.18 &	\$4.24 &	\$12.29 &80.15\% &	61.62\% &	113.01\% & 45.02\%\\
    1.0 & 0.6 & 0 &\$3.63 &\$7.88	& \$23.62 & 70.50\% &	78.38\% &	117.89\% &40.20\%  \\
    1.0 & 0.8 & 0 & \$6.63 &	\$10.43 &	\$25.89 & 64.61\% &	85.66\%	 & 120.39\%	 & 37.33\% \\
    1.0 & 1.0 & 0 & \$7.04 &\$12.00	& \$28.14 & 59.99\% &	90.04\% &	122.35\% & 35.93\% \\
    0.8 & 1.0 & 0 & \$7.17 &\$13.14 &\$29.66 &56.35\% &	92.79\% &	123.81\% &34.95\%  \\
    0.6 & 1.0 & 0 & \$7.47 &\$15.20 &	\$31.64 & 51.94\% &	95.39\% &	125.33\% & 32.95\%\\
    0.4 & 1.0 & 0 &\$6.91	& \$16.27 &	\$38.89 & 46.81\% &	97.57\% &	127.13\% & 31.82\%\\
    0.2 & 1.0 & 0 & \$8.59 &	\$18.59 &	\$45.27 & 40.03\% &	99.33\% &	129.17\% & 29.95\%\\
    0 & 1.0 & 0 & \$10.00 &	\$21.58 &	\$60.37& 30.91\% &	100.00\% &	131.65\% & 27.65\%\\
    0 & 1.0 & 0.2 & \$9.71 &	\$18.28&	\$42.20 &45.90\% &	97.77\%&	127.10\%	 & 29.95\%\\
    0 & 1.0 & 0.4 & \$8.52 &	\$16.41 &	\$36.63	 &59.59\%&	89.37\%&	121.46\%& 30.97\%\\
    0 & 1.0 & 0.6 &\$8.50 &	\$13.02 &	\$25.58 & 70.54\%&	76.54\%&	116.32\% &33.99\% \\
    0 & 1.0 & 0.8 & \$5.22 &	\$11.13 &	\$21.65 & 80.62\%&	56.33\%&	110.38\% & 35.17\%\\
    0 & 1.0 & 1.0 &\$1.83	& \$8.39&	\$14.59 & 90.11\%&	29.53\%	&104.39\% & 37.83\%\\
    0 & 0.8 & 1.0 & \$0.00 &	\$6.77 &	\$10.00& 95.40\%&	10.83\%&	100.87\% & 39.48\%\\
    0 & 0.6 & 1.0 & \$0.00 &	\$6.67 &	\$9.96 & 	95.58\%	&10.45\%&	100.81\% & 39.54\% \\
    0 & 0.4& 1.0 & \$0.00	& \$3.70	 & \$6.00 & 97.63\%	&6.53\%	&100.43\% &  44.04\% \\
    0 & 0.2 & 1.0 & \$0.00 &	\$0.00 &	\$0.00 & 100.00\%&	0.00\%&	100.00\% & 50.32\%\\
    0 & 0 & 1.0 & \$0.00	& \$0.00 &	\$0.00 &100.00\%&	0.00\%&	100.00\% & 50.32\%\\
    \bottomrule[1pt]
    \end{tabular}
\end{table}

A few representative solutions from Table \ref{T:testsuite} are presented in Table \ref{T:benchmark} alongside real-world fares, and the corresponding ridership values obtained by our model for the real-world fares. While all regimes have slightly lower ridership than that under real fares, all benchmarks achieve non-negligible improvements in system-wide metrics. In particular, the PR, PR+PX, and PR+VM regimes respectively achieve objective value increases of 47.4\%, 5.6\%, and 1.8\% respectively.

\begin{table}[htbp!]
    \centering
    \footnotesize
       \caption{\footnotesize Summary statistics of representative priority regimes. PR, PR+PX, and PR+VM regimes respectively have objective weights $(\pi^{PR}, \pi^{PX}, \pi^{VM})$ equal to $(1, 0, 0), (1, 1, 0), $ and $(1, 0, 1)$. ``\% routes discounted'' provides the proportion of routes in the system with activated discounts. ``Number of travelers'' is the total alliance passenger count originating within the alliance service region.}\label{T:benchmark}
    \begin{tabular}{l | c c c c c c c}
    \toprule[1pt]
    {\bf Priority regime} & \multicolumn{2}{c}{\bf Base fares (\$)} &\multicolumn{2}{c}{\bf Markups (\$/mile)} & {\bf Discount} & {\bf \% routes} & {\bf Number of} \\
    & {\bf MOD} & {\bf Transit} & {\bf MOD} & {\bf Transit} & {\bf Multiplier} & {\bf discounted} & {\bf travelers} \\\hline 
    Real fares & \$4.53 & \$4.50 & \$1.07 & \$0.16 & 0\% & 0\% &  25,816 \\
    PR & \$10.00 & \$10.00 & \$2.37 & \$0.63 & 50\% & 31.20\%  & 18,309 \\
    PR+PX &  \$10.00 &          \$8.50 &      \$0.56 &         \$0.25 & 31\% & 70.13\% & 23,793 \\
    PR+VM & \$10.00 &	\$1.83 &	\$0.16	& \$0.00 & 50\%  & 6.50\% & 25,051 \\
    \bottomrule[1pt]
    \end{tabular}
\end{table}

As shown in Table \ref{T:testsuite}, each set of priorities induces interpretable optimal prices and passenger decisions. The minimum, mean, and maximum real-world route prices in the service region are respectively \$4, \$10.23, and \$54.46, while those given by our model are in the range \$0-\$60.37; thus optimal fares are set at the correct order of magnitude in all regimes. Intuitively, route prices are the highest for the profit maximizing regime, and they decrease gradually as the importance of VMT or passenger utility increases. System-wide performance metrics are normalized against the best possible values across tested regimes, naturally achieved by each metric's corresponding single-objective optimization. The lowest profit is achieved in regimes that solely maximize passenger utility or minimize VMT, because very low prices achieve very low profit, but increase passenger happiness and entice more passengers away from single-occupancy vehicles. Analogously, highest fares achieve the highest profit, with more passengers electing to travel outside the system, and lowering overall passenger utility. System-wide outcomes vary smoothly with gradually changing alliance priorities. Even under single-objective profit maximization, route prices remain in the ballpark of real-world fares. While base fares and discount multiplier reach their upper limits, both optimal markups stay in the interior of the allowable range. These intuitive observations confirm that our fare-setting model is suitable for generating trustworthy qualitative insights.

At a quick glance, minimizing VMT and maximizing passenger utility seem to achieve similar outcomes in Table \ref{T:testsuite}---lower prices and higher system utilization. But Figure \ref{F:testsuite_fares} illustrates how each objective yields qualitatively very different designs. As VMT minimization increases in importance (moving from the middle to the right in Figure \ref{F:fares}), the markup is zeroed out, equalizing fares across longer and shorter routes. The elimination of a distance-based markup entices more longer-distance commuters to travel on the allied network, thus lowering VMT. When maximizing passenger utility, a more nuanced fare structure emerges to address heterogeneous passenger preferences. All pricing levers are employed: base fares, markups, discount multipliers, and discount activations. Higher markups and base fares are coupled with more numerous discounts across hybrid options, illustrated in the left halves of Figures \ref{F:fares} and \ref{F:discounts}. Thus, prioritization of each objective (PX versus VM) results in similar system-wide performance metrics by qualitatively different means. To extract the corresponding fare designs, we consider case studies prioritizing at most one of PX and VM at a time, with varying weights for profit. Section \ref{S:equity} further investigates geographic factors. {\color{black} Table \ref{T:testsuite_extended} in Appendix \ref{app:objTableExtension} elaborates on results that prioritize both passenger utility and VMT.} 

Finally, note that the variation in system utilization due to allied fare-setting is small compared to the integrated network's total loads. Before the pandemic, MBTA's peak hour bus load factor was already below 75\% \citep{loadfactors}. MBTA commuter rail transported around 120k passengers on an average weekday in 2018 \citep{mbtaCommuterRailLoad}, with 81.2\% of inbound ridership on peak trains \citep{commuterRailCounts}, yielding approximately 49k travelers on commuter rail during the AM rush. Moreover, approximately 116k daily TNC rides were destined for Boston in 2018, while another 15.3k originated in the alliance service region every day \citep{bostonTNCtraffic}. In contrast, the ridership numbers in Table \ref{T:benchmark} show that the alliance's ridership under real fares is a small proportion of the entire integrated network and that the system is capable of accommodating all demand redistribution as a result of allied fare-setting. In fact, Table \ref{T:benchmark} shows that the aggregate alliance ridership is slightly lower than that under real fares. Thus, drops on other routes will compensate for the slight ridership increases that may happen on certain routes under our proposed pricing alliance, ensuring that the system-wide transit load factors and MOD detour times are expected to remain largely unchanged as a result of the pricing alliance. We conclude that the linked resource reallocation problem need not be considered, when looking for rapid gains through pricing alliance formation.

\subsection{Equitable Access through a Refined Income-Aware Model Specification}\label{S:equity}

Our approach captures passengers' preferences and travel decisions when designing fares, a critical step to satisfying passenger needs. However, different groups of passengers have different preferences. Ignoring such differences can lead to inequitable and socially undesirable outcomes. After all, equity is a key driver for integrating on-demand services into public transport options. Acknowledging this challenge, we further refine our choice model and quantify the impacts of this nuanced model specification on system-wide metrics compared to an aggregated, average-case choice model.

Towns targeted for transit expansion by the MBTA in our case study have wide-ranging median household incomes, translating into varying price sensitivities. Affluent travelers' route choices are less susceptible to fare changes than those of low-income travelers. To partly account for such passenger heterogeneity, we compute a ratio of each town's median household income to the average of the median household incomes across the entire service region \citep{CensusBureau2018}. We use this income ratio to scale passengers' price sensitivities.
\vspace{-5mm}
\begin{table}[htbp!]
    \centering
    \footnotesize
        \caption{\footnotesize Allied system's daily morning rush ridership with (``Refined'') and without (``Base'') the choice model refinement, their percentage difference (Diff.), real-world ridership (Real), median household income (HHI), and average distance (Dist.) of alliance routes originating in the corresponding town and ending in the inner city.
    }\label{T:util}
    \resizebox{\textwidth}{!}{
    \begin{tabular}{lcc|ccc|ccc|ccc|c}
    \toprule[1pt]
    &\bf HHI & \bf Dist. &\multicolumn{3}{c|}{\bf PR+PX} & \multicolumn{3}{c|}{\bf PR} & \multicolumn{3}{c|}{\bf PR+VM} &  \\
    \bf Town & \bf \$K & \bf Miles & \bf Base & \bf Refined & \bf Diff. & \bf Base & \bf Refined & \bf Diff. & \bf Base & \bf Refined & \bf Diff. & \bf Real \\ \hline  
    Lawrence & 41.6 & 32.4 & 	1585 &	1982 & 125\%& 918 &	1881 & 205\% & 1331&	2065 & 155\% & 1566\\
    Lowell & 52.0 & 30.8 & 2000 &	2525 & 126\% & 1285 &	1565 & 122\% & 1852&	2633 & 142\% & 2115 \\
    Lynn & 54.6 & 13.5 &3539	& 3769 & 107\%	& 2283 &	2560 & 112\% & 3615&	3188 & 88\% & 4023 \\
    Brockton & 55.1 & 24.6 &	3498 &	3709 & 106\% & 1949&	3364 & 173\% & 3416&	3919 & 115\% & 3872 \\
    Salem & 65.6 & 18.9 &  1885 &	2043 & 108\% & 1189&	1382 & 116\% & 1961&	1786 & 91\% & 2171 \\
    Haverhill & 67.6 & 39.1 &	1605 &	1746 & 109\% & 1068&	1214 & 114\% & 1630&	1809 & 111\% & 1801 \\
    Framingham & 79.1 & 26.1 & 1203 &	1327 & 110\% & 	780	 &   936 & 120\% & 1093&	1117 & 102\% & 1198 \\
    Waltham & 85.7 & 11.8 & 1974 &	1873 & 95\% & 1597&	1713 & 107\% & 2101&	1946 & 93\% & 2249 \\
    Woburn & 88.7 & 14.3& 1795 &	1828 & 102\% & 1479&	1633 & 110\% & 2047&	1910 & 93\% & 2199 \\
    Stoneham & 94.8 & 11.6 & 	761	& 785 & 103\%	 & 	610	 &   691 & 113\% & 880	 &   830 & 94\% & 948 \\
    Wakefield & 95.3 & 13.3 & 1004 &	1007 & 100\%		 & 783	  &  889 & 114\% & 1140&	1076 & 94\% & 1228 \\
    Melrose & 103.7 & 9.8 &810	& 821	& 101\% &637	   & 711 & 112\%  & 889	 &   842 & 95\% & 949 \\
    Burlington & 105.4 & 17.6 & 1091 &	1144 & 105\% & 975	   & 1053 & 108\% & 1264&	1193 & 94\% & 1340\\
    Reading & 112.6 & 16.6 & 	 846 &	871	& 103\% & 703	   & 787 & 112\% & 972	 &   927 & 95\% & 1038 \\
    Winchester & 159.5 &10.8 & 197 &	207	& 105\% & 193	   & 198 & 103\% & 218	  &  210 & 96\% & 227 \\
    \bottomrule[1pt]
    \end{tabular}}
\end{table}

We compare the allied system ridership across priority regimes and towns under fares corresponding to base (i.e., income-agnostic) and refined (i.e., income-aware) choice models. We first calculate fares using the fare-setting model incorporating income-agnostic and income-aware choice models separately, and then evaluate both fare designs by calculating (and reporting in Table \ref{T:util}) ridership using only the income-aware choice model. Due to the use of the income-aware model, the system ridership in the PR+PX regime increases by 13\% on average for the towns with below-average median HHI compared to a less than 2\% average increase for the towns with above-average median HHI. While ridership increases across the board due to generally lower fares, the greatest increases occur in Lawrence and Lowell---the two towns with the lowest median household incomes. The lower middle-income bracket (Lynn, Brockton, Salem, Haverhill, Framingham) sees the next-highest ridership increase. Similar trends are observed in the PR regime: ridership increases across all towns as higher profit can be achieved with lower fares and higher volume. The above-average income towns gain 10\% in average ridership, while the below-average income towns gain 37\%.

In contrast to the 7\% and 23\% ridership gains in PR+PX and PR regimes, ridership grows by less than 4\% in the PR+VM regime. But interestingly, this regime has significantly larger ridership increases for longer distance passengers. In particular, the five towns farthest from the inner city (Lawrence, Lowell, Brockton, Haverhill, Framingham) are the only ones with ridership increases. Their average increase is $\sim$25\% while the remaining 10 towns see an average 7\% ridership drop. Thus, the PR+VM regime increases access to the allied network for commuters who are farther from the inner city. Overall, compared to the real-world fares, the prioritized system-wide objectives improved by 19.16\% for PR+PX, 0.35\% for PR and 0.99\% for PR+VM regime.

These results underline the importance of capturing passenger preference heterogeneity to amplify our model's practical impact. A choice model that reflects passenger preferences more accurately improves system-wide outcomes and especially maximizes passenger benefits. We conclude that our income-aware refined model improves transportation equity for passengers---as compared to the income-agnostic aggregate model---making it a valuable tool for transit agencies to incorporate into strategic decision-making.

\subsection{Quantifying the Value of Cooperation}\label{S:value_added}

To quantify the value of operator cooperation for operators and passengers, we solve the non-cooperative fare-setting model for all allied priorities in Table \ref{T:testsuite}. In each experiment, the transit operator's priorities are identical to alliance priorities, whereas the MOD operator always maximizes profit. Thus, the prospective alliance is assumed to adopt the transit operator's priorities.

\begin{figure}[htbp!]
    \centering
    \includegraphics[width=0.8\textwidth]{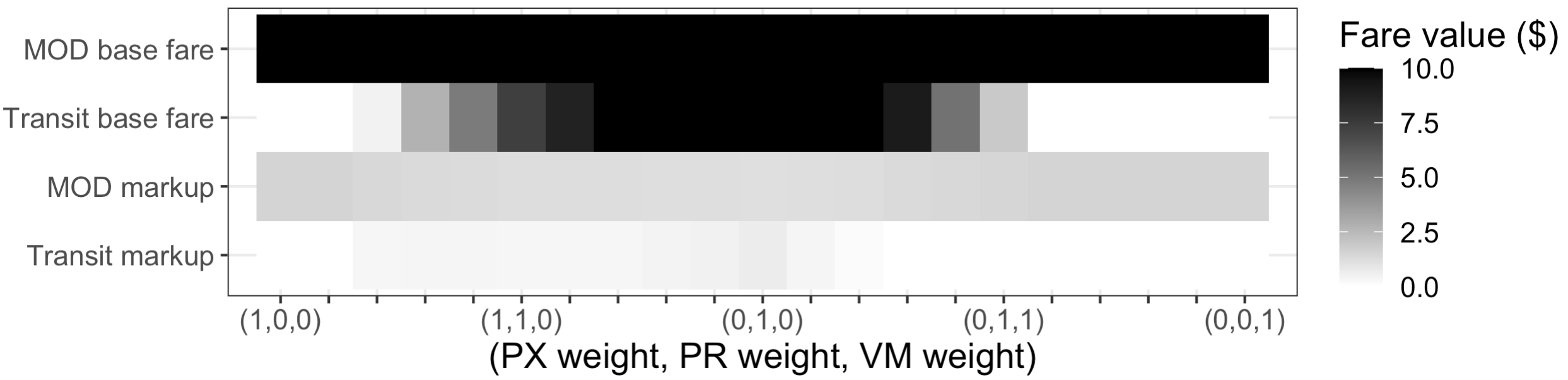}
    \caption{\footnotesize Non-cooperative fare parameters for different transit operator priorities. }
    \label{F:noncoop_fares}
\end{figure}

\begin{figure}[htbp!]
    \centering
    \includegraphics[width=0.8\textwidth]{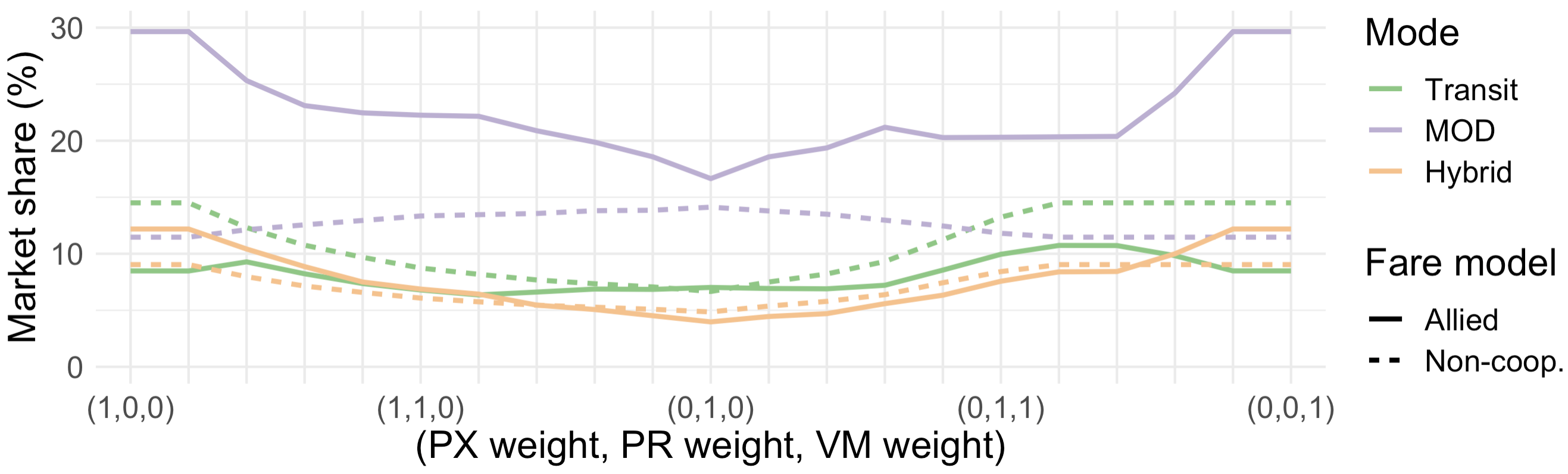}
    \caption{Market shares by mode across priority regimes for allied and non-cooperative fare-setting models.}
    \label{F:mode_choices}
\end{figure}

Figure \ref{F:noncoop_fares} lists the non-cooperative equilibrium fares. Discount multipliers are not applicable in the non-cooperative setting and hence excluded. Table \ref{T:noncooperative} compares allied and non-cooperative outcomes. {\color{black}We interpret the absolute profit outcomes by reporting the surplus allied profit ($\delta$), each operator's allied profit allocation,} \% increase in the alliance objective compared to that under the non-cooperative setting (transit obj. \% inc.), and the percentage increase in MOD profit due to the alliance (MOD obj. \% inc.). The profit allocations are determined by the profit allocation mechanism in Section \ref{S:profit_allocation}. {\color{black}The surplus allied profit is negative whenever the alliance accrues less total profit than in the non-cooperative regime, potentially resulting in a negative profit allocation for the transit operator if the MOD non-cooperative profit exceeds total allied profit.} The non-cooperative system utilization and average route prices are also provided. Finally, Figure \ref{F:mode_choices} shows passenger mode choices across all tested regimes for both allied and non-cooperative fares.

\begin{table}[htbp!]
    \centering
    \footnotesize
        \caption{\footnotesize Allied vs. non-cooperative outcomes. Transit operator's priorities represented as $\bm{\pi}_{TR}$. Profit amounts in thousands of dollars (\$K).}
    \label{T:noncooperative}
    \resizebox{\textwidth}{!}{\begin{tabular}{ccc|ccc|ccccc|cc}
         \toprule[1pt] 
          \multicolumn{3}{c|}{Transit obj. weights} & MOD non-  & Transit non- & Total allied & $\delta$ & MOD allied & Transit allied & Transit obj. &MOD obj. & System & Route price\\
         $\pi_{TR}^{PX}$ & $\pi_{TR}^{PR}$ & $\pi_{TR}^{VM}$ & coop. prof.  & coop. prof.  & prof.  & & prof. alloc. & prof. alloc. & \% inc. &  \% inc.   & util. \% & (\$) Mean \\ \hline
     1.0 & 0 & 0 &501.90	&0.00	   &0.00	    &-501.90	&501.90	&-501.90 & 9.60\%  & 0.00\%  & 35.01\% & \$12.83\\
    1.0 & 0.2 & 0 &501.90	&0.00	   &0.00	    &-501.90	&501.90	&-501.90&6.87\%  & 0.00\% &35.01\% & \$12.83\\
    1.0 & 0.4 & 0 &487.91	&1,935.20	&2,336.80	&-86.31	   &487.91	&1,848.89& 4.93\% & 0.00\%  &32.41\% & \$13.45  \\
    1.0 & 0.6 & 0 &479.04	&2,488.03	&2,816.71	&-150.36	&479.04	&2,337.68&4.80\% & 	0.00\%  &30.46\% & \$14.78\\
    1.0 & 0.8 & 0 &473.18	&2,734.56	&3,097.45	&-110.29	&473.18	&2,624.27& 5.45\% & 0.00\%& 29.21\% &\$15.83 \\
    1.0 & 1.0 & 0 &468.14	&2,871.10	&3,292.77	&-46.46	   &468.14	&2,824.63& 7.57\%& 0.00\%  &28.16\% &\$16.96\\
    0.8 & 1.0 & 0 &465.06	&2,968.75	&3,406.00	&-27.81	   &465.06	&2,940.94& 27.71\%& 0.00\% & 27.39\% &\$17.77 \\
    0.6 & 1.0 & 0 &462.29	&3,043.50	&3,484.80	&-20.99	   &462.29	&3,022.50&5.82\% & 0.00\%&26.70\% &\$18.54\\
    0.4 & 1.0 & 0 &460.45	&3,123.17	&3,562.44	&-21.18	   &460.45	&3,101.99&1.31\% & 0.00\% &26.43\% & \$18.66 \\
    0.2 & 1.0 & 0 &458.83	&3,168.79	&3,637.17	&9.56	    &463.61	&3,173.57& 0.54\% & 1.04\% &26.00\% &\$18.97 \\
    0 & 1.0 & 0 &456.73	&3,188.96	&3,657.81	&12.12	    &462.79	&3,195.02& 0.33\%& 1.33\%  &25.65\% & \$19.24\\
    0 & 1.0 & 0.2 &461.34	&3,088.11	&3,590.19	&40.74	    &481.71	&3,108.48&1.55\% & 4.42\% &26.66\% &\$18.50\\
    0 & 1.0 & 0.4  &465.54	&2,839.69	&3,305.07	&-0.16	    &465.54	&2,839.52&0.52\% & 0.00\% &27.50\% &\$18.07 \\
    0 & 1.0 & 0.6 &471.77	&2,356.26	&2,868.99	&40.96	    &492.25	&2,376.74& 0.37\%& 4.34\%&28.67\% &	\$17.40 \\
    0 & 1.0 & 0.8 &482.97	&1,619.43	&2,056.25	&-46.15	   &482.97	&1,573.28&0.51\% & 0.00\%&31.13\% &\$15.33 \\
    0 & 1.0 & 1.0 &494.38	&684.88	    &1,085.64	&-93.61	    &494.38	&591.27& 0.54\%& 0.00\%& 	33.48\% & \$13.76\\
    0 & 0.8 & 1.0 &501.90	&0.00	   &396.31	    &-105.59	&501.90	&-105.59&0.60\%  & 0.00\% & 35.01\% &\$12.83\\
    0 & 0.6 & 1.0 &501.90	&0.00	   &382.17	    &-119.73	&501.90	&-119.73& 0.73\% & 0.00\% & 35.01\% &\$12.83\\
    0 & 0.4& 1.0 &501.90	&0.00	   &238.92	    &-262.97	&501.90	&-262.97&0.92\%  & 0.00\%& 35.01\% &\$12.83 \\
    0 & 0.2 & 1.0 &501.90	&0.00	   &0.00	    &-501.90	&501.90 &-501.90& 1.36\%& 0.00\% &35.01\% &	\$12.83 \\
    0 & 0 & 1.0 &501.90	&0.00	   &0.00	    &-501.90	&501.90	&-501.90& 1.90\%& 0.00\%& 35.01\% &\$12.83\\
         \bottomrule[1pt]
    \end{tabular}}
\end{table}

\vspace{-5mm}
Figure \ref{F:noncoop_fares} illustrates that the MOD operator's profit-maximizing strategy remains relatively constant, regardless of the transit operator's priorities. Still, when the transit operator prioritizes VMT or passenger benefits, the MOD operator selects higher fare parameters (mainly through higher markups) than in the corresponding allied setting. Thus, the average route price (Route price (\$) Mean) columns in Tables \ref{T:testsuite} and \ref{T:noncooperative} show that the non-cooperative average route prices are higher than average allied route prices in every regime except for the one where the transit operator only prioritizes profit. {\color{black}In these scenarios, the MOD operator is exploiting the additional room in passenger budgets afforded by lower transit fares, to further maximize their own profit.}

When prioritizing passenger benefits or VMT, the alliance sets MOD fare parameters lower than in the non-cooperative scenario (Figure \ref{F:fares} vs. \ref{F:noncoop_fares}). Figure \ref{F:mode_choices} shows that MOD-only market shares decrease and hybrid market shares increase in the non-cooperative scenario as passenger benefits or VMT are increasingly prioritized. We conclude that, although fewer passengers utilize MOD-only options in non-cooperative scenarios where the transit operator is VMT- or passenger-oriented, more passengers select hybrid options due to the very low (or free) transit fares observed in Figure \ref{F:noncoop_fares}. Thus, the MOD operator earns more profit in those non-cooperative scenarios in which the transit operator is more altruistic. This is reflected in the profit allocation mechanism: observe in Table \ref{T:noncooperative} that the MOD operator earns strictly more profit in almost all regimes where the transit operator is not solely a profit maximizer, even though the system as a whole generates strictly less total profit, as seen in comparison with Table \ref{T:testsuite}. The profit allocation mechanism ensures that the MOD operator receives their non-cooperative earnings, despite the lower MOD fares that the alliance sets to achieve lower VMT or higher passenger benefits. This in turn reduces the transit's profit allocation as high VMT or low passenger benefits are increasingly penalized. On the other hand, the transit operator always strictly improves its objective of optimizing total system-wide performance, however it chooses to define it. As a result, the MOD operator interestingly finds it in its interest to adopt transit's priorities as transit increasingly diverges from profit maximization. In other words, {\it the profit-maximizing MOD operator would not prefer a profit-maximizing alliance}. In fact, the MOD operator would benefit most from total altruism on the transit side (an exclusive focus on either passenger utility or VMT reduction). Passengers win due to the strictly lower prices and higher system utilization that result from such alliances. {\color{black}We contrast this conclusion with the observation that strictly positive surplus profit is observed in profit-maximizing alliances (as supported by Lemma \ref{L:profit_allocation}c and the corresponding discussion in Section \ref{S:profit_allocation}). After all, operators are not only strategically setting fares, but also deciding on their respective---or collective---priorities.}

The transit operator must ultimately set the ceiling in terms of the price they are willing to pay for the alliance benefits. Table \ref{T:noncooperative} shows that the transit agency runs a deficit to appease the MOD operator if its profit emphasis is too low. A deficit alone might not be enough to dissuade the transit agency from participating in the alliance: as we have noted, every public transit mode operates at a loss, and especially so for the on-demand options \citep{ktp2016}. To weigh the financial implications of an alliance, the agency may compare the magnitude of the loss to the cost of the analogous MOD system operated by transit on their own in the absence of outsourcing through an alliance. Oftentimes, transit agencies also receive grants to fund pricing alliances \citep{sandbox}, and the daily deficit rate can be compared to the grant award amount and intended duration. To enable a daily deficit rate that keeps pace with financial resources over time, the transit agency may propose to reset overall performance metric goals to induce different optimal fares, or to adjust the geographic and/or temporal scope of the alliance. In the end, each transit agency is expected to choose the trade-off point between its financial, passenger-focused, and environmental goals that most closely aligns with their overall policy and various practical and financial constraints. Regardless, our allied and non-cooperative fare-setting models together with our profit allocation mechanism jointly provide a toolkit usable by transit agencies to weigh these trade-offs as they evaluate a potential pricing alliance.
\section{Conclusions, Limitations, and Future Directions}\label{S:conclusions}

We contribute a pricing alliance design framework to enable incentive-aligned collaboration between transit agencies and MOD operators. Our allied fare-setting model captures the interdependent decisions of passengers and operators. We accomplish large scale by developing a tractable two-stage fare-setting formulation equivalent to the original mixed-integer non-convex optimization problem, and solve it with a tailored SOS2 coordinate descent approach. Our approach consistently and significantly outperforms benchmarks, enabling additional daily system-wide benefits worth tens of thousands dollars. Practically speaking, our framework aligns profit-oriented MOD operators with transit goals of passenger utility and single-occupancy VMT reduction. In other words, cooperative pricing results in win-win-win outcomes for passengers, MOD operators, and transit agencies. Finally, by tuning passenger route choice models, the alliance can prioritize lower fares and higher utilization for low-income or long-distance commuters, thus improving passenger equity.

Consistent with pricing alliance practice, we assumed that the MOD operator serves as a contractor to the transit agency and agrees to set static fares for trips in the integrated system. While this ensures transparent communication of public sector prices, it may be possible to set fare schemes that allow MOD operators to maintain dynamic prices, perhaps through the transit operator subsidizing passenger trip costs up to a fixed dollar amount. Additionally, we model average-case travel times over a fixed set of route options for the MOD operator's portion of the network to represent average-case operations, which simplifies their typically dynamic routing scheme. In other words, we consider the case where the alliance sets one permanent fare scheme that is optimized for average-case performance. Future research may consider optimizing for the worst-case performance, and/or integrating dynamic routing into the fare-setting model, requiring the integration of two complex problem classes. {\color{black} Finally, we did not incorporate joint resource reallocation into the pricing scheme, due to the observation that the transit operator is capable of accommodating all demand redistributions attributed to changing prices, and the MOD fleet is right-sized accordingly. This assumption works well for contemporary North American transit systems because of the prevalence of low ridership, empty seats and low load factors. Future research could consider relaxing this assumption to generalize the analysis beyond the North American context by jointly modeling the capacity allocation and pricing decisions.}

\ACKNOWLEDGMENT{%
This material is based upon work supported by the National Science Foundation under Grant no. 1122374.
}

%
%
%


\bibliographystyle{informs2014trsc} 
\bibliography{ref.bib} 

\newpage 
\begin{APPENDICES}
{\color{black}
\section{Real-World Pricing Alliances}\label{A:real-alliances}
}

Here we provide examples of real-world pricing alliance maps in Figure \ref{F:zoning}. Figure \ref{F:plano} illustrates a zone-based route structure, while Figure \ref{F:newton} illustrates a hub-based route structure.

\begin{figure}[htbp!]
    \centering
     \begin{subfigure}[b]{0.4\textwidth}
         \centering
         \includegraphics[width=\textwidth]{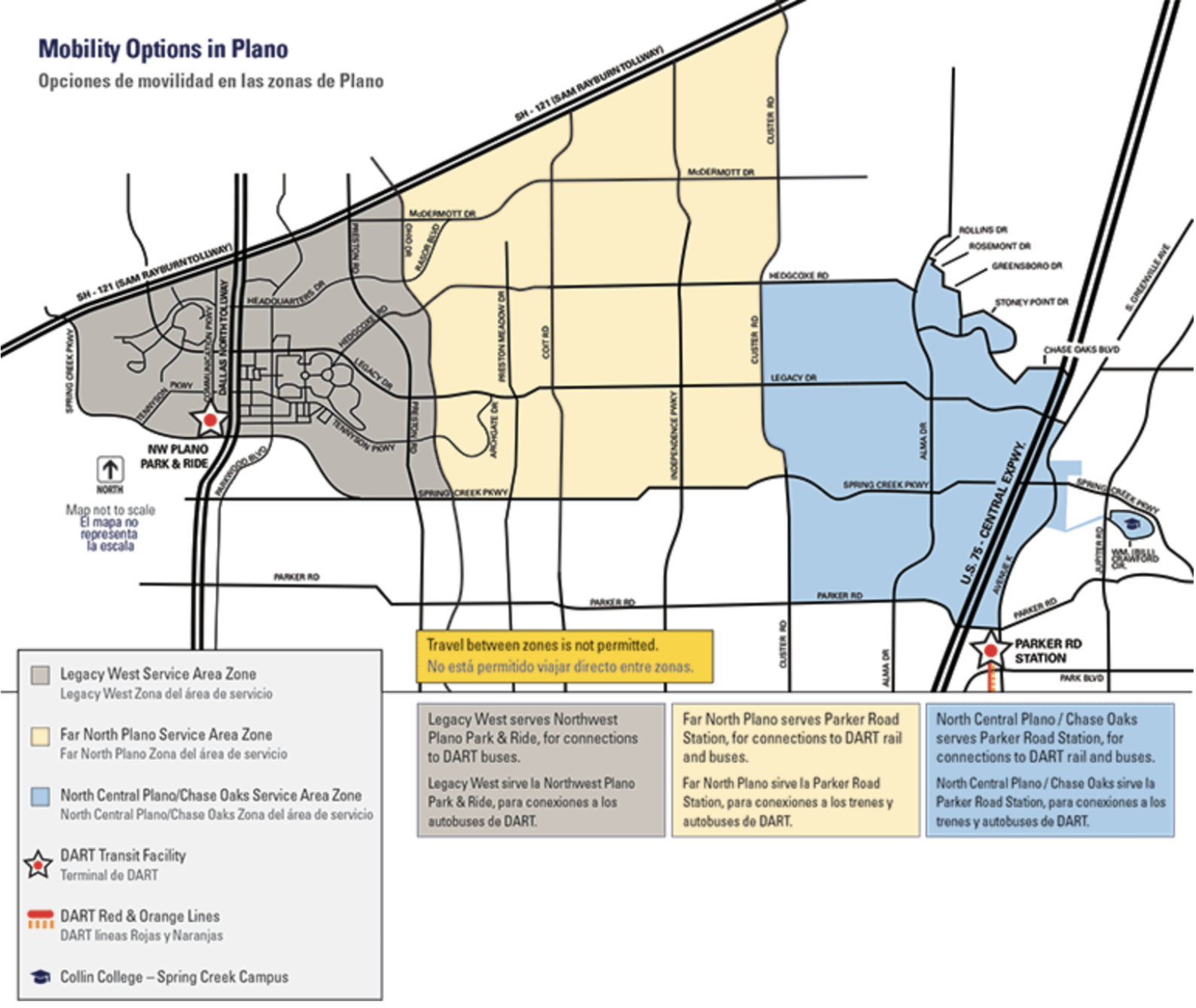}
         \caption{ GoLink in Dallas, TX (Plano region) utilizes a zone-based route structure.}
         \label{F:plano}
     \end{subfigure}
     \qquad \qquad
     \begin{subfigure}[b]{0.45\textwidth}
         \centering
         \includegraphics[width=\textwidth]{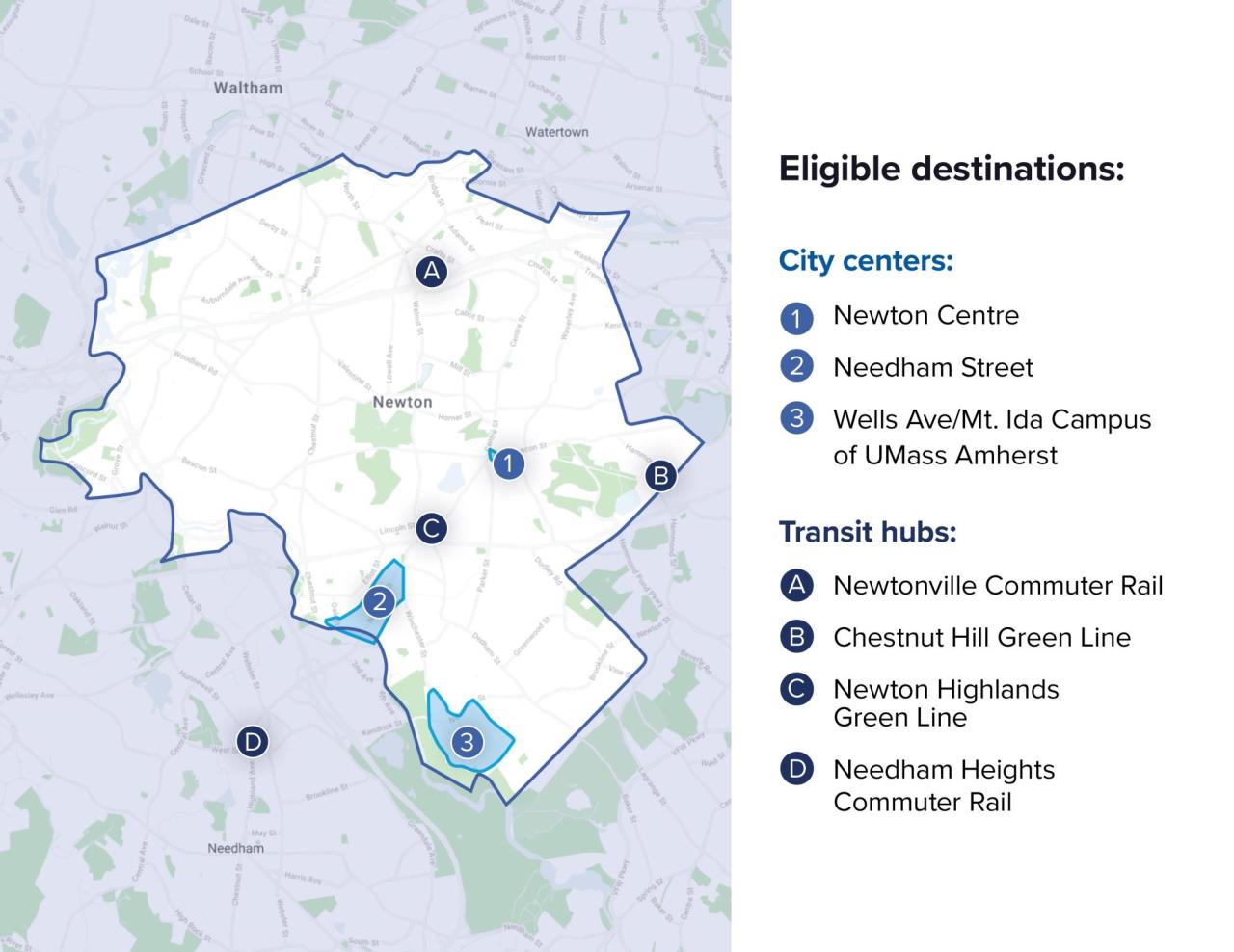}
         \caption{ The {\it NewMo} Pilot in Newton, MA utilized a hub-based route structure.}
         \label{F:newton}
     \end{subfigure}
    \caption{ Route structures of recent pricing alliances. In Plano \citep{dart}, passengers spend up to \$3 to travel anywhere within a color-blocked region. In the {\it NewMo} Pilot \citep{newmo}, passengers could travel anywhere within town limits for \$2, as long as either trip endpoint was one of seven hubs. {\it NewMo} now allows passengers to travel anywhere within Newton.}
    \label{F:zoning}
\end{figure}

\section{Proofs of Lemmas}\label{A:proofs}

\subsection{Proof of Lemma \ref{L:equiv}.}\label{A:equivalent}

Given an optimal solution from each formulation, we construct a feasible solution for the other with the same objective value. We denote the non-discounted route price $r \in \mathcal{R}$ as $\sigma_r(\bm{f}) = \sum_{k \in \mathcal{O}_r} (\beta_k + \mu_k\Delta_{rk})$, and define $\gamma_{irt}(\bm{f}, \Lambda)= \frac{\exp(u_{i0t})}{\exp(u_{irt} + \alpha_i (1 - \Lambda) \sum_{k \in \mathcal{O}_r}(\beta_k + \mu_k \Delta_{rk}))}$ as the ratio of the outside option's attractiveness for passenger $i \in \mathcal{N}$ to the attractiveness of route $r \in \mathcal{R}_i$, as a function of the fare parameters and percent discount.

\paragraph{PADP-FS2SD $\to$ PADP-FS.} 
Let $(\bm{f}^{2SD}, \Lambda^{2SD}, \bm{x}^{2SD}, \bm{s}^{2SD}, \bm{w}^{2SD}, \bm{p}^{2SD}) $ be an optimal solution to the PADP-FS2SD model with objective value $z^{2SD}$. We construct the following solution to the PADP-FS model.
\begin{equation}
    \bm{f}^{FS} = \bm{f}^{2SD}, \quad
    \Lambda^{FS} = \Lambda ^{2SD}, \quad
    \bm{x}^{FS} = \bm{x}^{2SD}, \quad
    \bm{s}^{FS} = \bm{s}^{2SD}, \quad
    \bm{p}^{FS} = \bm{p}^{2SD}
\label{form2_sol}
\end{equation}
We demonstrate feasibility of solution \eqref{form2_sol} to the PADP-FS model. $\bm{x}^{FS}$ are binary by construction. Constraints \eqref{C:prices1} and \eqref{C:prices2} hold by construction. Constraints \eqref{C:badshare} and \eqref{C:badshare0} are feasible due to Constraints \eqref{E:share_prod}, \eqref{E:sum}, and \eqref{E:share_prod_triv}. Specifically, for a given $i \in \mathcal{N}$, Constraints \eqref{E:share_prod} and \eqref{E:share_prod_triv} ensure that $s_{irt}^{FS} \propto\exp(u_{irt} + \alpha_i \cdot p_r^{FS}) $ for each $r \in \mathcal{R}_i, t \in \mathcal{T}$, and $s_{i0t}^{FS} \propto \exp(u_{i0t}) $ for each $ t \in \mathcal{T}$. Constraint \eqref{E:sum} forces each $s_{irt}$ and $s_{i0t}$ to be normalized by $\exp(u_{i0t}) + \sum_{r \in \mathcal{R}_i} \exp(u_{irt} + \alpha_i \cdot p_r^{FS})$, and the result follows. Constraints \eqref{C:first} - \eqref{C:lambda_last} hold by definitions of $\mathcal{F}$ and $\mathcal{L}.$
Finally, to show that solution \eqref{form2_sol} has the same objective value as $z^{2SD},$ we must only prove that $\left(p_r^{FS}- \sum_{k \in \mathcal{O}_r} c_k \Delta_{rk} \right) \cdot s_{irt}^{FS} = w_{irt}^{2SD}$ for all $i \in \mathcal{N}$, $r \in \mathcal{R}_i,$ and $t \in \mathcal{T},$ which follows directly from Constraints \eqref{E:disc_prod} and \eqref{E:disc_prod_triv}.

\paragraph{PADP-FS $\to$ PADP-FS2SD.} Let $(\bm{f}^{FS}, \Lambda^{FS}, \bm{p}^{FS}, \bm{s}^{FS}, \bm{x}^{FS}) $ be an optimal solution to the PADP-FS model with objective value $z^{FS}$. We construct the following solution to the PADP-FS2SD model.
\begin{equation}
\begin{cases}
    \bm{f}^{2SD} = \bm{f}^{FS}, \quad
    \Lambda^{2SD} = \Lambda ^{FS}, \quad
    \bm{x}^{2SD} = \bm{x}^{FS}, \quad
    \bm{s}^{2SD} = \bm{s}^{FS}, \quad
    \bm{p}^{2SD} = \bm{p}^{FS} \\
    w_{irt}^{2SD} = \left(p_r^{FS} - \sum_{k \in \mathcal{O}_r} c_k \Delta_{rk} \right)\cdot s_{irt}^{FS} \quad \forall i \in \mathcal{N}, r \in \mathcal{R}_i, t \in \mathcal{T}
\end{cases} \label{form3_sol}
\end{equation}

First we demonstrate the feasibility of \eqref{form3_sol} to PADP-FS2SD. Variables $\bm{s}^{2SD}$, $\bm{w}^{2SD}$ and $\bm{p}^{2SD}$ are non-negative and $\bm{x}^{2SD}$ are binary by construction. For Constraints (\ref{E:share_prod}) and \eqref{E:share_prod_triv}, fix some $i \in \mathcal{N}.$ For $r \in \mathcal{R}_i \setminus \mathcal{R}^D$ and $t \in \mathcal{T},$ and by expressions \eqref{C:badshare} and \eqref{C:badshare0},
    \begin{align*}
        \gamma_{irt}(\bm{f}^{2SD} , 0) \cdot s_{irt}^{2SD} &=  \gamma_{irt}(\bm{f}^{FS} , 0) \cdot s_{irt}^{FS} = \frac{\exp(u_{i0t})}{\exp(u_{irt} + \alpha_i \cdot p_r^{FS})} \cdot \frac{\exp(u_{irt} + \alpha_i \cdot p_r^{FS})}{\exp(u_{i0t}) + \sum_{s \in \mathcal{R}_i} \exp(u_{ist} + \alpha_i \cdot p_s^{FS})} \\
        &= \frac{\exp(u_{i0t})}{\exp(u_{i0t}) + \sum_{s \in \mathcal{R}_i} \exp(u_{ist} + \alpha_i \cdot p_s^{FS})} = s_{i0t}^{FS}= s_{i0t}^{2SD}.
    \end{align*}
    Thus, Constraints \eqref{E:share_prod_triv} are satisfied. Now, consider any $a \in \mathcal{A}$ with $r \in \mathcal{R}_i \cap \mathcal{R}_a^D.$ If $x_a^{2SD} = 0$, then we can use the exact same argument as above to confirm that $\gamma_{irt}(\bm{f}^{2SD} , 0) \cdot s_{irt}^{2SD} = s_{i0t}^{2SD}$ again holds.
    Since
    \begin{align*}
        &\gamma_{irt}(\bm{f}^{2SD}, \Lambda^{2SD}) \cdot s_{irt}^{2SD} \leq \gamma_{irt}(\bm{f}^{2SD}, 0) \cdot s_{irt}^{2SD} = s_{i0t}^{2SD} \leq \gamma_{irt}(\bm{f}^{2SD}, 0)\\
        &\leq\gamma_{irt}(\bm{f}^{2SD}, 0) + \gamma_{irt}(\bm{f}^{2SD}, \Lambda^{2SD}) \cdot s_{irt}^{2SD}= \gamma_{irt}(\bm{f}^{2SD}, \Lambda^{2SD}) \cdot s_{irt}^{2SD} + M_{irt}^s \cdot (1 - x_a^{2SD}),
    \end{align*}
    so, Constraints \eqref{E:share_prod} hold. Similar logic applies when $x_a^{2SD} = 1.$
For Constraints (\ref{E:disc_prod}) and \eqref{E:disc_prod_triv}, fix any $i \in \mathcal{N}$. Constraints \eqref{E:disc_prod_triv} hold by construction for all $r \in \mathcal{R}_i \setminus \mathcal{R}^D$ and $t \in \mathcal{T}.$ 
    For $a \in \mathcal{A}$ with $r \in \mathcal{R}_i \cap \mathcal{R}_a^D$ and $t \in \mathcal{T},$ we have $w_{irt}^{2SD} = \left(p_r^{2SD}- \sum_{k \in \mathcal{O}_r} c_k \Delta_{rk}\right) \cdot s_{irt}^{2SD} \leq \sigma_r(\bm{f}^{2SD}) \cdot s_{irt}^{2SD}$ regardless of the value of $x_a^{2SD}.$ If $x_a^{2SD} = 0,$  
    \begin{align*}
        &\left((1 - \Lambda^{2SD}) \cdot \sigma_r(\bm{f}^{2SD})- \sum_{k \in \mathcal{O}_r} c_k \Delta_{rk}\right) \cdot s_{irt}^{2SD} \leq w_{irt}^{2SD}=\left(\sigma_r(\bm{f}^{2SD}) - \sum_{k \in \mathcal{O}_r} c_k \Delta_{rk}\right) \cdot s_{irt}^{2SD} \\
        &= \left((1 - \Lambda^{2SD}) \cdot \sigma_r(\bm{f}^{2SD}) - \sum_{k \in \mathcal{O}_r} c_k \Delta_{rk}\right) \cdot s_{irt}^{2SD} + \Lambda^{2SD} \cdot \sigma_r(\bm{f}^{2SD}) \cdot s_{irt}^{2SD} \\
        &\leq ((1 - \Lambda^{2SD}) \cdot \sigma_r(\bm{f}^{2SD}) - \sum_{k \in \mathcal{O}_r} c_k \Delta_{rk}) \cdot s_{irt}^{2SD} + \Lambda^{2SD} \cdot \sigma_r(\bm{f}^{2SD})\\
        &= ((1 - \Lambda^{2SD}) \cdot \sigma_r(\bm{f}^{2SD}) - \sum_{k \in \mathcal{O}_r} c_k \Delta_{rk}) \cdot s_{irt}^{2SD} + M_{ir}^w \cdot (1 - x_a^{2SD}),
    \end{align*}
    and hence Constraints (\ref{E:disc_prod}) hold. A similar argument applies when $x_a^{2SD} = 1.$ From expressions \eqref{C:badshare} and \eqref{C:badshare0}, it is easy to see that \eqref{E:sum} holds.
Constraints \eqref{E:price_nodisc} and \eqref{E:price_disc}, as well as Constraints $\bm{f} \in \mathcal{F}$ and $\Lambda \in \mathcal{L}$ hold by construction. Now we show that solution \eqref{form3_sol} to PADP-FS2SD has the same objective value of $z^{FS}$.
\begin{align*}
    z^{2SD} &= \sum_{t \in \mathcal{T}}\sum_{i \in \mathcal{N}} N_{it} \Big[\pi^{PX}  (u_{i0t} + \sum_{r \in \mathcal{R}_i} (u_{irt} + \alpha_i p_r^{2SD}) + \pi^{PR} \sum_{r \in \mathcal{R}_i} w_{irt}^{2SD} - \pi^{VM} (\Delta^0_i s_{i0t}^{2SD}) \Big] \\
    &= \sum_{t \in \mathcal{T}}\sum_{i \in \mathcal{N}} N_{it} \Big[\pi^{PX}  (u_{i0t} + \sum_{r \in \mathcal{R}_i} (u_{irt} + \alpha_i p_r^{FS}) + \pi^{PR} \sum_{r \in \mathcal{R}_i} \left(p_r^{FS}- \sum_{k \in \mathcal{O}_r} c_k \Delta_{rk} \right) s_{irt}^{FS} - \pi^{VM} (\Delta^0_i s_{i0t}^{FS}) \Big]= z^{FS}
\end{align*}

\color{black}
\subsection{Proof of Lemma \ref{L:feas}.}\label{A:feas}

Given a feasible first-stage solution $(\widehat{\mathbf{f}}, \widehat{\Lambda}) \in \mathcal{F} \times \mathcal{L},$ we demonstrate that a solution to the second-stage problem always exists, i.e., that $\mathcal{S}(\widehat{\mathbf{f}}, \widehat{\Lambda})\neq \emptyset$. In particular, the non-discounted solution is always feasible.

\begin{equation}
\begin{cases}
    \overline{x}_a = 0 &a \in \mathcal{A} \\
     \overline{p}_r = \sum_{k \in \mathcal{O}_r} \big(\widehat{\beta}_k + \widehat{\mu}_k\Delta_{rk}  \big) &r \in \mathcal{R} \\
     \overline{s}_{i0t} = \frac{\exp(u_{i0t})}{\exp(u_{i0t}) + \sum_{r \in \mathcal{R}_i} \exp(u_{irt} + \alpha_i \cdot  \overline{p}_r)} &i \in \mathcal{N} \\
     \overline{s}_{irt} = \frac{\exp(u_{irt} + \alpha_i \cdot  \overline{p}_r)}{\exp(u_{i0t}) + \sum_{s \in \mathcal{R}_i} \exp(u_{ist} + \alpha_i \cdot  \overline{p}_s)} &i \in \mathcal{N}, r \in \mathcal{R}_i \\
     \overline{w}_{irt} = \left(\overline{p}_r- \sum_{k \in \mathcal{O}_r} c_k \Delta_{rk} \right) \cdot \overline{s}_{irt} &i \in \mathcal{N}, r \in \mathcal{R}_i
\end{cases} \label{feas_sol}
\end{equation}

We demonstrate that solution \eqref{feas_sol} is feasible to the second-stage problem, i.e., that it lies within $\mathcal{S}(\widehat{\bm{f}}, \widehat{\Lambda}).$ Throughout the proof, we continue to use the notation $\gamma_{irt}$ defined in the previous subsection.

Constraint \eqref{E:share_prod_triv}:
    Fix $i \in \mathcal{N}$, $r \in \mathcal{R}_i \setminus \mathcal{R}^D,$ and $t \in \mathcal{T}.$
    \begin{align*}
        \frac{\exp(u_{i0t})}{\exp(u_{irt}+\alpha_i \cdot \overline{p}_r)} \cdot \overline{s}_{irt}
        &=\frac{\exp(u_{i0t})}{\exp(u_{irt}+\alpha_i \cdot \overline{p}_r)} \cdot \frac{\exp(u_{irt} + \alpha_i \cdot  \overline{p}_r)}{\exp(u_{i0t}) + \sum_{s \in \mathcal{R}_i} \exp(u_{ist} + \alpha_i \cdot  \overline{p}_s)} &&\text{Value of } \overline{s}_{irt} \\
        &= \frac{\exp(u_{i0t})}{\exp(u_{i0t}) + \sum_{s \in \mathcal{R}_i} \exp(u_{ist} + \alpha_i \cdot  \overline{p}_s)}= \overline{s}_{i0t} &&\text{Value of }\overline{s}_{i0t}
    \end{align*}
    
Constraint \eqref{E:share_prod}: Fix $a \in \mathcal{A}, i \in \mathcal{N}_a, 
r \in \mathcal{R}_i \cap \mathcal{R}_a^D,$ and $t \in \mathcal{T}.$ We first note that $\overline{s}_{irt} \gamma_{irt}(\widehat{\bm{f}}, 0) = \overline{s}_{i0t}$ due to the same algebraic argument demonstrated for Constraint \eqref{E:share_prod_triv}. It follows that
$\overline{s}_{i0t} \leq \gamma_{irt}(\widehat{\bm{f}}, 0) \overline{s}_{irt}$
as well as $\overline{s}_{i0t} \geq \gamma_{irt}(\widehat{\bm{f}}, 0) \overline{s}_{irt} - M_{irt}^s \overline{x}_a$.
Further, the value of $\overline{p}_r$ yields 
$\gamma_{irt}(\widehat{\bm{f}}, \widehat{\Lambda}) = \exp(u_{i0t})/\exp(u_{irt} + \alpha_i (1 - \widehat{\Lambda}) \overline{p}_r).$ Thus, by observing $\exp(u_{irt} + \alpha_i (1 - \widehat{\Lambda}) \overline{p}_r) \geq \exp(u_{irt} + \alpha_i \overline{p}_r),$ we obtain $\gamma_{irt}(\widehat{\bm{f}}, \widehat{\Lambda}) \leq \gamma_{irt}(\widehat{\bm{f}}, 0).$ Therefore, $\overline{s}_{i0t} = \gamma_{irt}(\widehat{\bm{f}}, 0) \overline{s}_{irt} \geq \gamma_{irt}(\widehat{\bm{f}}, \widehat{\Lambda}) \overline{s}_{irt}$. We show that the final inequality holds with simple algebra:
$\overline{s}_{i0t} = \gamma_{irt}(\widehat{\bm{f}}, 0) \overline{s}_{irt}\leq \gamma_{irt}(\widehat{\bm{f}}, \widehat{\Lambda}) \overline{s}_{irt} + \gamma_{irt}(\widehat{\bm{f}}, 0) \overline{s}_{irt}\leq \gamma_{irt}(\widehat{\bm{f}}, \widehat{\Lambda}) \overline{s}_{irt} + \gamma_{irt}(\widehat{\bm{f}}, 0) = \gamma_{irt}(\widehat{\bm{f}}, \widehat{\Lambda}) \overline{s}_{irt} + M_{irt}^s (1 - \overline{x}_a)$. Here, the second inequality holds because $\overline{s}_{irt} \leq 1$ and the last equality holds because of the definition of $M_{irt}^s$.
    
Constraint \eqref{E:disc_prod}: Fix $a \in \mathcal{A}$, $i \in \mathcal{N}_a,$ $r \in \mathcal{R}_i \cap \mathcal{R}_a^D,$ and $t \in \mathcal{T}.$ We first note that $\overline{w}_{irt} = \left(\overline{p}_r - \sum_{k \in \mathcal{O}_r} c_k \Delta_{rk} \right) \cdot \overline{s}_{irt} = \left(\sum_{k \in \mathcal{O}_r} (\widehat\beta_k + (\widehat\mu_k-c_k)\Delta_{rk})\right)\cdot \overline{s}_{irt}$ due to the value of $\overline{p}_r.$ It follows that $\overline{w}_{irt} \leq \left(\sum_{k \in \mathcal{O}_r} (\widehat\beta_k + (\widehat\mu_k-c_k)\Delta_{rk})\right) \cdot \overline{s}_{irt}$ and $\overline{w}_{irt} \geq \left(\sum_{k \in \mathcal{O}_r} (\widehat\beta_k + (\widehat\mu_k-c_k)\Delta_{rk})\right) \cdot \overline{s}_{irt} - M_{ir}^w \cdot\overline{x}_a$. The discounted price is bounded above by the non-discounted price, so
    \begin{align*}
        \overline{w}_{irt} = \left(\sum_{k \in \mathcal{O}_r} (\widehat\beta_k + (\widehat\mu_k-c_k)\Delta_{rk})\right) \cdot \overline{s}_{irt} \geq \left((1 - \widehat{\Lambda}) \cdot \Big(\sum_{k \in \mathcal{O}_r} (\widehat\beta_k + \widehat\mu_k\Delta_{rk})\Big)- \sum_{k \in \mathcal{O}_r} c_k \Delta_{rk} \right) \cdot \overline{s}_{irt}.
    \end{align*}

    Finally, we leverage bounds on market shares to validate the final inequality.
    {\small\begin{align*}
        \overline{w}_{irt} &\leq \left((1 - \widehat{\Lambda}) \cdot \left(\sum_{k \in \mathcal{O}_r} (\widehat\beta_k + \widehat\mu_k\Delta_{rk})\right) - \sum_{k \in \mathcal{O}_r} c_k \Delta_{rk}\right) \cdot \overline{s}_{irt} + \widehat{\Lambda}\cdot \left(\sum_{k \in \mathcal{O}_r} (\widehat\beta_k + \widehat\mu_k\Delta_{rk})\right) \cdot \overline{s}_{irt}  &&\text{Rearranging}\\
        &\leq \left((1 - \widehat{\Lambda}) \cdot \left(\sum_{k \in \mathcal{O}_r} (\widehat\beta_k + \widehat\mu_k\Delta_{rk})\right) - \sum_{k \in \mathcal{O}_r} c_k \Delta_{rk}\right) \cdot \overline{s}_{irt} + \widehat{\Lambda}\cdot \left(\sum_{k \in \mathcal{O}_r} (\widehat\beta_k + \widehat\mu_k\Delta_{rk})\right)&&\text{Value of }\overline{s}_{irt} \leq 1 \\
        &= \left((1 - \widehat{\Lambda}) \cdot \left(\sum_{k \in \mathcal{O}_r} (\widehat\beta_k + \widehat\mu_k\Delta_{rk})\right) - \sum_{k \in \mathcal{O}_r} c_k \Delta_{rk}\right) \cdot \overline{s}_{irt} + M_{ir}^w \cdot (1 - \overline{x}_a) &&\text{Value of }\overline{x}_a = 0
    \end{align*}}
    
Constraint \eqref{E:sum}: Fix $i \in \mathcal{N}$ and $t \in \mathcal{T}.$
    \begin{align*}
        \overline{s}_{i0t} + \sum_{r \in \mathcal{R}_i} \overline{s}_{irt} &= \frac{\exp(u_{i0t})}{\exp(u_{i0t}) + \sum_{r \in \mathcal{R}_i} \exp(u_{irt} + \alpha_i \cdot  \overline{p}_r)} + \sum_{r \in \mathcal{R}_i} \frac{\exp(u_{irt} + \alpha_i \cdot  \overline{p}_r)}{\exp(u_{i0t}) + \sum_{s \in \mathcal{R}_i} \exp(u_{ist} + \alpha_i \cdot  \overline{p}_s)} \\
        &= \frac{\exp(u_{i0t}) + \sum_{r \in \mathcal{R}_i}\exp(u_{irt} + \alpha_i \cdot  \overline{p}_r)}{\exp(u_{i0t}) + \sum_{r \in \mathcal{R}_i}\exp(u_{irt} + \alpha_i \cdot  \overline{p}_r)}= 1
    \end{align*}
    
Constraint \eqref{E:disc_prod_triv} holds by construction. Constraint \eqref{E:price_nodisc} is satisfied by definition of $\overline{p}_r$ for $r \in \mathcal{R} \setminus \mathcal{R}^D.$ For Constraint \eqref{E:price_disc}, fix $a \in \mathcal{A}$ and $r \in \mathcal{R}_a^D.$ Then, $\overline{p}_r = \sum_{k \in \mathcal{O}_r} \left(\widehat{\beta}_k + \widehat{\mu}_k\Delta_{rk} )\right) = (1 - \widehat{\Lambda} \overline{x}_a) \left(\sum_{k \in \mathcal{O}_r} \left(\widehat{\beta}_k + \widehat{\mu}_k\Delta_{rk} \right)\right)$.
The second equality holds because the value of $\overline{x}_a=0$.
Variables $\overline{x}_a = 0$ are binary for each $a \in \mathcal{A}$, and the market shares $\mathbf{s}$ are non-negative due to the non-negative range of the exponential function.

\color{black} 
\subsection{Proof of Lemma \ref{L:profit_allocation}.}\label{App:profit_allocation}

\paragraph{Proof of Lemma \ref{L:profit_allocation}a:}
The MOD operator participates in the alliance whenever they can earn at least as much profit by cooperating with the transit agency as they can otherwise. By construction, for $k\in\mathcal{O}\setminus\{TR\}$, $\Phi_k((\Theta_l(\bm{f}^{nc}))_{l \in \mathcal{O}}, \Theta(\bm{f}^a, \Lambda^a,\bm{\pi})) = \Theta_k(\bm{f}^{nc}) + \big(\tfrac{\delta}{|\mathcal{O}|}\big)^+ \geq \Theta_k(\bm{f}^{nc}), $ and so the result follows.

\paragraph{Proof of Lemma \ref{L:profit_allocation}b:}
When the allied profit equals/exceeds combined non-cooperative profit, we show that the allocation is a Nash bargaining solution, which is a classic payment rule guaranteeing properties of Pareto efficiency, symmetry, scale invariance, and independence of irrelevant alternatives. We frame the alliance profit allocation problem as a collective bargaining problem. Let the disagreement outcome be the operator profits resulting from non-cooperation, $\bm{d}=(\Theta_k(\bm{f}^{nc}))_{k \in \mathcal{O}}$. Let $\mathcal{X}$ be the set of all profit allocations such that each operator receives at least their disagreement outcome, and such that the sum of profit allocations does not exceed allied profit $p = \Theta(\bm{f}^{a}, \Lambda^a;\bm{\pi})$. $\mathcal{X}(\bm{d},p) := \Big\{\big(a_k\big)_{k \in \mathcal{O}}: a_k \geq d_k, \forall k \in \mathcal{O}; \, \sum_{k \in \mathcal{O}} a_k \leq p\Big\}$.

Let $\mathcal{F}$ be the set of all such allocation problems, with each $(\bm{d},p) \in \mathcal{F}$ corresponding to a different allocation problem (i.e., a different potential alliance). We seek an allocation solution function $\Phi: \mathcal{F} \to \mathcal{X}$ that allocates the available profit according to the axioms of Pareto efficiency, symmetry, scale invariance, and independence of irrelevant alternatives. It is known that Nash bargaining solutions---which satisfy the above axioms---exactly coincide with optimal solutions to optimization problem \eqref{NashBargaining}.
\begin{equation}\label{NashBargaining}
     \max_{\bm{a} \in \mathcal{X}(\bm{d},p)}  \prod_{k \in \mathcal{O}}(a_k-d_k) 
\end{equation}
Expression \eqref{E:nash} clearly solves \eqref{NashBargaining}.
\begin{equation}\label{E:nash}
     d_k + \frac{p - \sum_{k \in \mathcal{O}}d_k}{|\mathcal{O}|} = \Theta_k(\bm{f}^{nc}) + \frac{\delta}{|\mathcal{O}|} \equiv \Phi_k(\bm{d}, p), \quad \forall k \in \mathcal{O}
\end{equation}
To establish the core property, we point to Lemma (\ref{L:profit_allocation}a), as well as the fact that $\Phi_{TR}(\bm{d}, p) \geq \Theta_{TR}(\bm{f}^{nc})$ in the case that $\delta\geq0$. This proves Lemma (\ref{L:profit_allocation}a).

\color{black}
\paragraph{Proof of Lemma \ref{L:profit_allocation}c:} 

First we consider the case when $\pi^{PR} > 0, \pi^{PX} = 0,$ and $\pi^{VM} = 0.$ Consider the feasible allied solution $(\bm{f}^{nc}, 0) \in \mathcal{F} \times \mathcal{L},$ i.e., the non-cooperative fare parameters with a zero discount multiplier. This solution accrues allied profit $\sum_{k \in \mathcal{O}} \Theta_k(\bm{f}^{nc})$, providing the lower bound $\Theta(\bm{f}^a, \Lambda^a; \bm{\pi}) \geq \sum_{k \in \mathcal{O}} \Theta_k(\bm{f}^{nc}),$ which yields $\delta \geq 0.$

Next, we consider the case when $\pi^{PR} = 0$ and either $\pi^{PX} > 0$ or $\pi^{VM} > 0,$ assuming $c_k=0$ for each operator $k.$ Each passenger $i \in \mathcal{N}$ has a negative utility per unit price, i.e., $\alpha_i < 0.$ Thus, decreasing any price $p_r$ for any route $r \in \mathcal{R}_i$ results in larger utilities $u_{irt} + \alpha_i p_r$. Without any profit incentive, allied route prices will decrease as much as they are allowed---in other words, we will see $p_r = 0$ for all $r \in \mathcal{R},$ since $c_k=0.$ These minimum prices either maximize total utility (the sum of all utilities), or maximize utilization of intrasystem routes (which according to the MNL choice model, is proportional to the attractiveness of the routes, which increases with utility). In the non-cooperative regime, the individually acting MOD operator will set strictly positive prices on all of their routes to maximize their profit. Thus, $\Theta_{MOD}(\bm{f}^{nc}) > \Theta(\bm{f}^a, \Lambda^a; \bm{\mu}) = 0,$ ensuring that $\delta < 0.$

\color{black}

\section{Model Notation}\label{App:notation}

Here, Table \ref{T:notation} summarizes all notation from Section \ref{S:formulation}.
\vspace{-5mm}

\begin{table}[htbp!]
\centering
\scriptsize
\caption{\footnotesize Notation for the PADP-FS and PADP-FS2SD.}\label{T:notation}
\begin{tabular}{llp{11.5cm}}
\toprule[1pt]
{\bf Component} & {\bf Type} & {\bf Description} \\
\hline
$\mathcal{A}$ & Set & Discount activation categories \\
$\mathcal{N}$ & Set & Passenger types \\
$\mathcal{O}$ & Set & Operators \\
$\mathcal{R}$ & Set & Intrasystem routes, not including the outside option \\
$\mathcal{T}$ & Set & Discrete set of time intervals \\
$\mathcal{O}_r$ & Set & Operators who help service route $r \in \mathcal{R}$ \\
$\mathcal{R}_i$ & Set & Route options available to passengers of type $i \in \mathcal{N}$ \\
$\mathcal{R}^D$ & Set & Discount-eligible routes \\
$\mathcal{R}_a^D$ & Set & Routes in discount activation category $a \in \mathcal{A}$ \\
$\mathcal{F}$& Set & Allowable fare parameter values \\
$\mathcal{L}$& Set & Allowable discount multiplier values\\
$\mathcal{N}_a$ & Set & Passenger types with $\geq 1$ discount-eligible route in category $a\in\mathcal{A}$\\
\hline
\color{black} $N_{it}$ & \color{black}
 Param. & Number of passengers of type $i \in \mathcal{N}$ at time $t \in \mathcal{T}$ \\
\color{black} $c_k$ & \color{black} Param. & Marginal cost for operator $k \in \mathcal{O}$ per unit distance per passenger \\ 
$\Delta_i^0$ & Param. & Distance driven by passenger of type $i \in \mathcal{N}$ when the outside option is selected \\
$\Delta_{rk}$ & Param. & Distance the passenger travels with operator $k \in \mathcal{O}$ on route $r \in \mathcal{R}$ \\
$u_{irt}$ & Param. & Non-monetary utility to passenger type $i \in \mathcal{N}$ on route $r \in \mathcal{R}_i$ at time $t \in \mathcal{T}$ \\
$u_{i0t}$ & Param. & Utility accrued by a passenger of type $i \in \mathcal{N}$ by driving at time $t \in \mathcal{T}$ \\
$\alpha_i$ & Param. & Utility per unit price to a passenger of type $i \in \mathcal{N}$ \\
$B$ & Param. & Maximum allowable base fare \\
$M$ & Param. & Maximum allowable distance-based markup \\
$L$ & Param. &Maximum allowable value of discount multipliers\\
$\pi^{PX}, \pi^{PR}, \pi^{VM}$ & Param. & Relative priority weights of system-wide performance metrics\\
$M_{ir}^w, M_{irt}^s$& Parameter & Big-$M$ parameters for each $i \in \mathcal{N}, r \in \mathcal{R}_i, t \in \mathcal{T}$ \\
\hline 
$x_a$ &Variable & Binary. Whether to activate discount option $a \in \mathcal{A}$ \\
$\beta_k, \mu_k$ &Variable & Continuous. Base fare and markup of operator $k \in \mathcal{O}$. Also, $\bm{f}=(\bm{\beta}, \bm{\mu})$ \\
$p_r$ &Variable & Continuous. Customer-facing price of route $r \in \mathcal{R}$ \\
$\Lambda$ & Variable & Continuous. Discount multiplier applied to routes with activated discounts \\
$s_{irt}$ &Variable & Continuous. Share of passengers of type $i \in \mathcal{N}$ choosing route $r \in \mathcal{R}_i$ at $t \in \mathcal{T}$ \\ 
$s_{i0t}$ &Variable & Continuous. Share of passengers of type $i \in \mathcal{N}$ choosing the outside option at $t \in \mathcal{T}$\\
$w_{irt}$ & Variable & Continuous. Represents profit for passengers of type $i \in \mathcal{N}$, route $r \in \mathcal{R}_i$, time $t \in \mathcal{T}$ \\
\bottomrule[1pt]
\end{tabular}
\end{table}
\vspace{-7mm}

\color{black}
\section{SOS2 Coordinate Descent Subroutines}\label{App:sos2}

This section specifies two subroutines of SOS2 Coordinate Descent, presented in Algorithm \ref{A:SOS2_CD}. {\sc Search Directions} provides a comprehensive ordered list of search directions that can be multidimensional and/or randomized (or neither), while {\sc Generate Anchors} computes an ordered set of evenly spaced SOS2 anchors along the specified search direction. In case of a multidimensional search direction, it determines the range of slopes that will ensure that the line spans the selected dimension, and samples uniformly from that range.

\begin{algorithm}
\scriptsize 
\caption{Subroutines for SOS2 Coordinate Descent}\label{A:SOS2_CD_prelim}
\begin{algorithmic}[1]
\Procedure{Search Directions}{$random$, $multidim$}

\If{$multidim$}

\State $searchDirs \leftarrow \{(\beta_{TR}, \mu_{TR}),\: (\beta_{MOD}, \mu_{MOD}),\: \Lambda\} $ {\tt \quad // index names of }$\bm{y}$

\Else 

\State $searchDirs \leftarrow \{\beta_{TR},\: \mu_{TR},\: \beta_{MOD},\: \mu_{MOD},\: \Lambda\} $

\EndIf

\If{$random$}

\State $searchDirs \leftarrow \text{\sc shuffle}(searchDirs)$ {\tt \quad  // randomly permutes the set $searchDirs$}

\EndIf

\State \Return{$searchDirs$}

\EndProcedure

\Procedure{Generate Anchors}{$(\bm{f}, \Lambda)$, $searchDir$, $D$}

\If{$searchDir \in \{(\beta_{k'}, \mu_{k'}) \, : \, k' \in \mathcal{O}\}$} {\tt \quad // if $searchDir$ is multidimensional}

\State $\beta_k, \mu_k \leftarrow searchDir$

\If{Unif$(0, 1) < 0.5$} {\tt \quad // select spanning dimension}

\State $(x, x_{min}, x_{max}), (y, y_{min}, y_{max}) \leftarrow  (\beta_k, 0, B), (\mu_k, c_k, M)$ 

\Else 

\State $(x, x_{min}, x_{max}), (y, y_{min}, y_{max}) \leftarrow(\mu_k,  c_k, M), (\beta_k, 0, B)$

\EndIf

\If{$x == 0$} \label{slopeBegin} {\tt \quad // determine valid slopes for spanning affine lines}

\State $m_{min} \leftarrow \tfrac{ y_{min} - y}{x_{max} - x}$; \, $m_{max} \leftarrow \tfrac{y_{max} - y}{x_{max} - x}$ 

\ElsIf{$x == x_{max}$}

\State $m_{min} \leftarrow \tfrac{y_{max} - y}{x_{min} - x}$; \, $m_{max} \leftarrow \tfrac{y_{min} - y}{x_{min} - x}$

\Else 

\State $m_{min} \leftarrow \max\{\tfrac{y_{max} - y}{x_{min} - x}, \tfrac{y_{min} - y}{x_{max} - x}\}$; \, $m_{max} \leftarrow \min\{\tfrac{y_{min} - y}{x_{min} - x},  \tfrac{y_{max} - y}{x_{max} - x}\}$ \label{slopeEnd} {\tt \quad // see Figure \ref{F:slopes}}

\EndIf 

\State $m \leftarrow \text{Unif}(m_{min}, m_{max})$; \, $b \leftarrow y - m \cdot x$; \, $\delta_x = x_{max} - x_{min}$ {\tt \quad // slope of spanning affine line}

\State $anchors \leftarrow \{(x_{min} + \tfrac{i-1}{D-1}\cdot \delta_x, b + m \cdot (x_{min} + \tfrac{i-1}{D-1}\cdot \delta_x), \beta_{-k}, \mu_{-k}, \Lambda) \, : \, i \in \{1,\dots,D\} \}$

\Else 

\If{$searchDir == \Lambda$} 

\State $maxVal \leftarrow L; \, minVal \leftarrow 0$

\ElsIf{$searchDir$ in $\{\beta_k \, : \, k \in \mathcal{O}\}$}

\State $maxVal \leftarrow B; \, minVal \leftarrow 0$

\Else 

\State $maxVal \leftarrow M; \, minVal \leftarrow c_k$

\EndIf

\State $anchors \leftarrow \{(minVal + \tfrac{i-1}{D-1} \cdot (maxVal - minVal), (\bm{f}, \Lambda) \setminus \{searchDir\}) \, : \,i \in \{1,\dots,D\}\} $

\EndIf

\State Insert $(\bm{f}, \Lambda)$ into the ordered set $anchors$

\State \Return{$anchors$}

\EndProcedure

\end{algorithmic}
\end{algorithm}

Figure \ref{F:slopes} visualizes the computation of the slope range in lines \eqref{slopeBegin}-\eqref{slopeEnd} of the {\sc Generate Anchors} procedure. The most negative slope of a line passing through the current solution is determined by the maximum of the slopes of the two line segments connecting the current solution to the upper left and bottom right corners of the diagram. Similarly, the most positive slope is determined by the minimum of the slopes of the two line segments connecting the current solution to the lower left and upper right corners.

\begin{figure}[htbp!]
    \centering
    \includegraphics[width=0.4\textwidth]{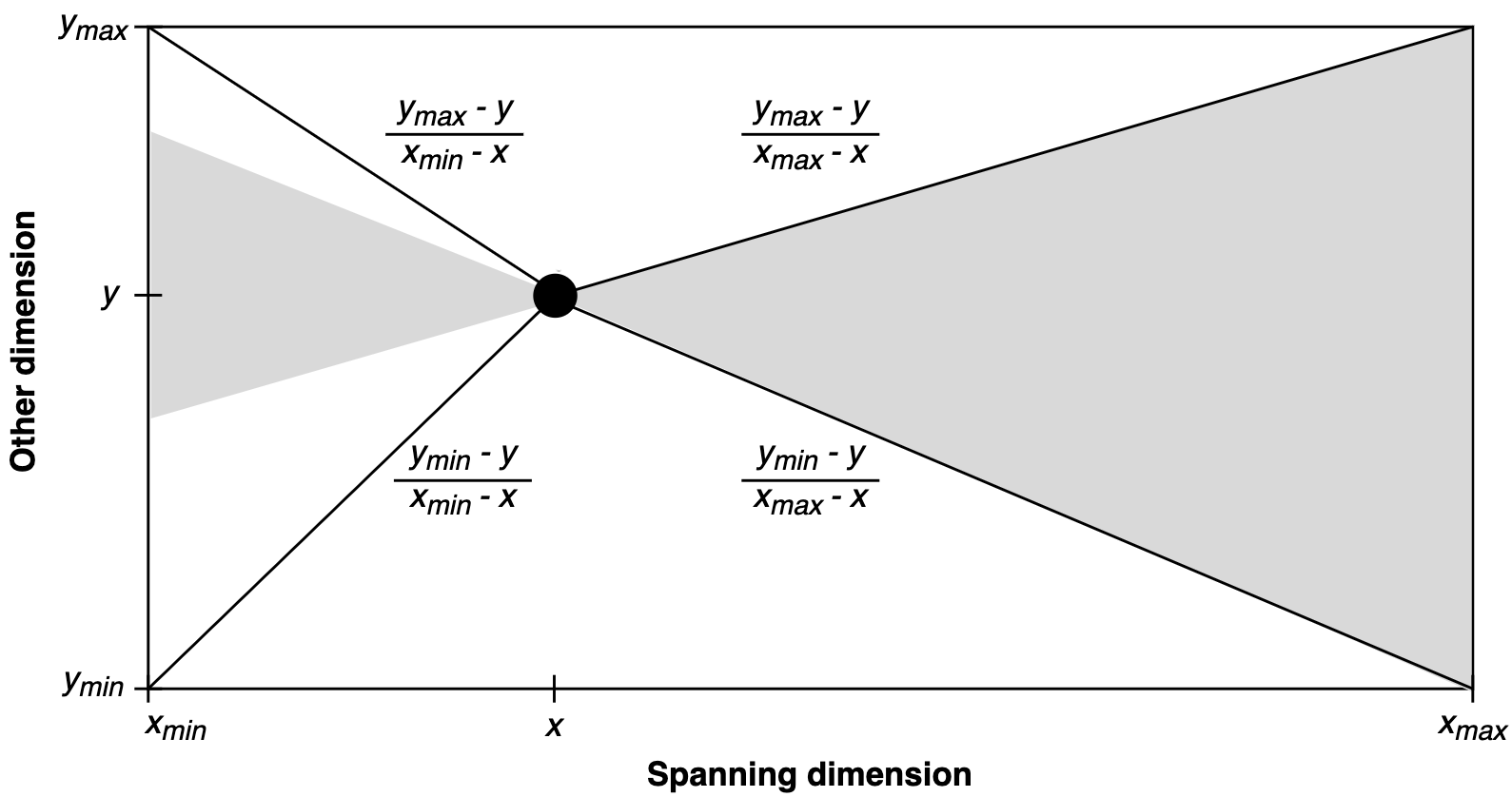}
    \caption{Computing the range of slopes such that the selected dimension is spanned.}
    \label{F:slopes}
\end{figure}

\section{Bayesian Optimization Benchmark}\label{A:BO}

Bayesian Optimization (BO) is a sequential search strategy for optimizing low-dimensional black-box functions \citep{m2012}. Typically, the function being optimized takes a long time to evaluate and has no analytical form, precluding access to gradients. We first provide a high-level overview of BO, and then we summarize the design choices that we use in this paper. In particular, we adapt the setup of a recent dedicated study by \cite{lbds2019} on using Bayesian Optimization to select MOD system service parameters subject to passenger mode choice. Interested readers are referred to \cite{lbds2019} for a more detailed description of BO in this application context. Readers interested in a general BO tutorial are referred to \cite{bcd2010}. In the absence of a closed-form representation of our black-box function $f: \bm{x} \to \mathbbm{R}$, the BO framework first imposes a prior belief upon $f$ via a {\it probabilistic surrogate model}. Given this surrogate model, BO iteratively (i) updates the likelihood of historical observations with new evaluations of $f$ to obtain a more informative posterior, and (ii) queries an {\it acquisition function} that uses the updated posterior to recommend the next value of $\bm{x}$ that should be evaluated. A very common surrogate model for the black-box function is called a Gaussian Process (GP), which is a stochastic process in which any finite set of random variables follows a multivariate Gaussian distribution \citep{m2012}. After building the GP with historical observations, the GP maps a given point in the search space to a univariate Gaussian distribution. We interpret this output distribution as a set of potential values of $f(\bm{x})$, accounting for noise. The acquisition function is a BO design choice intended to help in navigating the trade-off between exploration and exploitation, i.e., exploring more of the search space versus exploiting regions where a globally optimal solution is suspected to exist. As in \cite{lbds2019}, we use a GP upper confidence bound as our acquisition function, which characterizes the BO optimization process as a multi-armed bandit problem.

\section{Additional Computational Insights}\label{App:component_performance}

\subsection{Performance of SOS2-CD Subroutines.}
In this section, we examine the performance of individual components of the algorithm. 
Figure \ref{Fa:solvetime} illustrates the distribution of the second-stage model's solution times (in sec). Each calculation of $W$ requires solving the second-stage model once, and needs on average 5.7 seconds of CPU time, which is fast enough to be useful, but slow enough to warrant judicious selection of candidate first-stage solution points. Each of the 31,553 observations was obtained in under one minute. The observations were aggregated across all solutions of the second-stage model represented in the paper, including runs of SOS2-CD, BF-CD, and BO.

\begin{figure}[htbp!]
     \centering
     \begin{subfigure}[b]{0.45\textwidth}
         \centering
         \includegraphics[width=\textwidth]{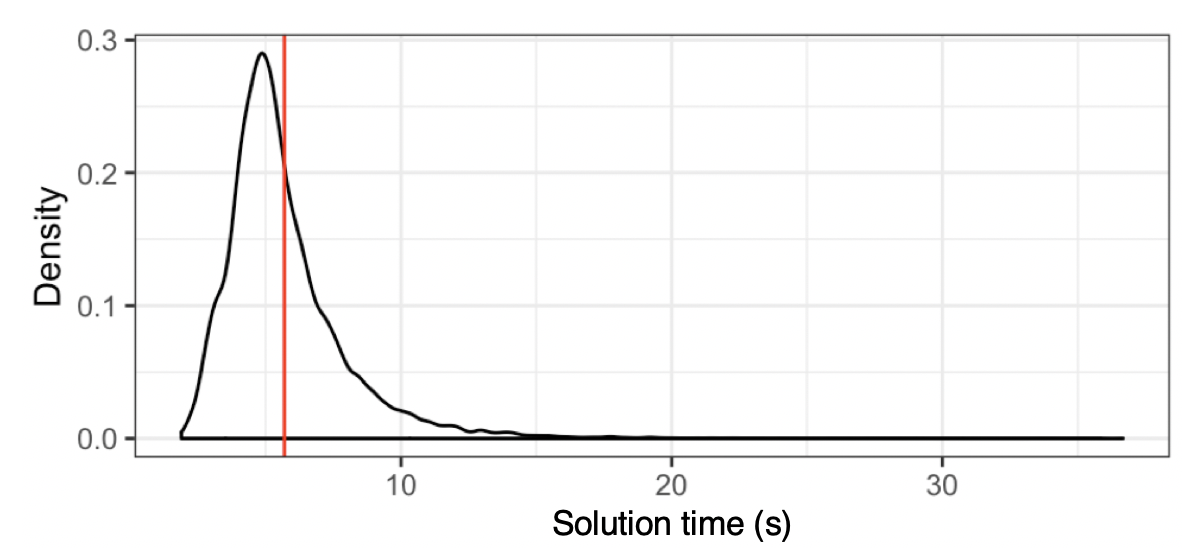}
         \caption{Second-Stage Model Solution Times}
         \label{Fa:solvetime}
     \end{subfigure}
     \hfill
     \begin{subfigure}[b]{0.45\textwidth}
         \centering
         \includegraphics[width=\textwidth]{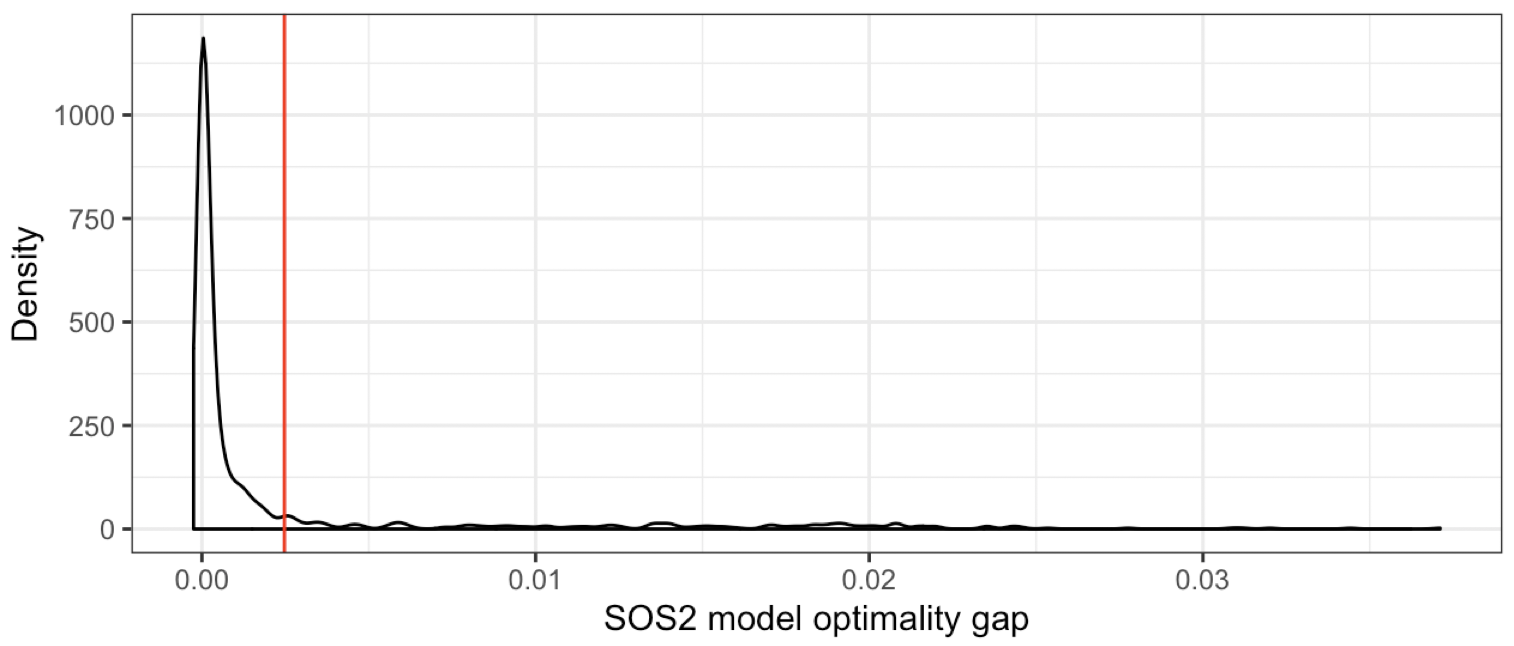}
         \caption{SOS2 Model Optimality Gaps}
         \label{Fa:sos2_val}
     \end{subfigure}
        \caption{Distributions of (a) 31,553 second-stage model solution times (with mean 5.7 sec), and (b) SOS2 model optimality gaps (with mean gap $= 0.25$\%, across 927 solutions. Exact solutions (with 0\% gaps) obtained in 276 of the 927 cases.)}
        \label{fig:three graphs}
\end{figure}


Figure \ref{Fa:sos2_val} depicts the accuracy of the SOS2 model with anchors placed at 10\% intervals for the discount multiplier axis and at \$1 intervals for the rest of the fare parameters. In each trial, a fare parameter combination and a search dimension were selected uniformly at random. The ``true optimal'' fare parameters for each trial (i.e., each combination of current solution and search direction) were computed using a brute-force approach, through exhaustive enumeration of every solution along that line at a 1\% granularity for discount multiplier and a \$0.10 granularity for the other four fare parameters. Then the $W$ obtained through $SOS2^*$ procedure was evaluated and compared with this ``true optimal'' solution, to compute the SOS2 optimality gap. The mean optimality gap was 0.25\%, and it was 0\% in 29.8\% of the instances. This establishes the trustworthiness of the SOS2 model outputs, especially given the drastic CPU time reduction they provide. We repeated this procedure for 10 hours, resulting in 927 solutions.
Figure \ref{F:bfcd} shows the distribution of BF-CD solution times for a single trajectory across 50 trials. BF-CD never terminates before an hour elapses.
\begin{figure}[htbp!]
    \centering
    \includegraphics[width=0.45\textwidth]{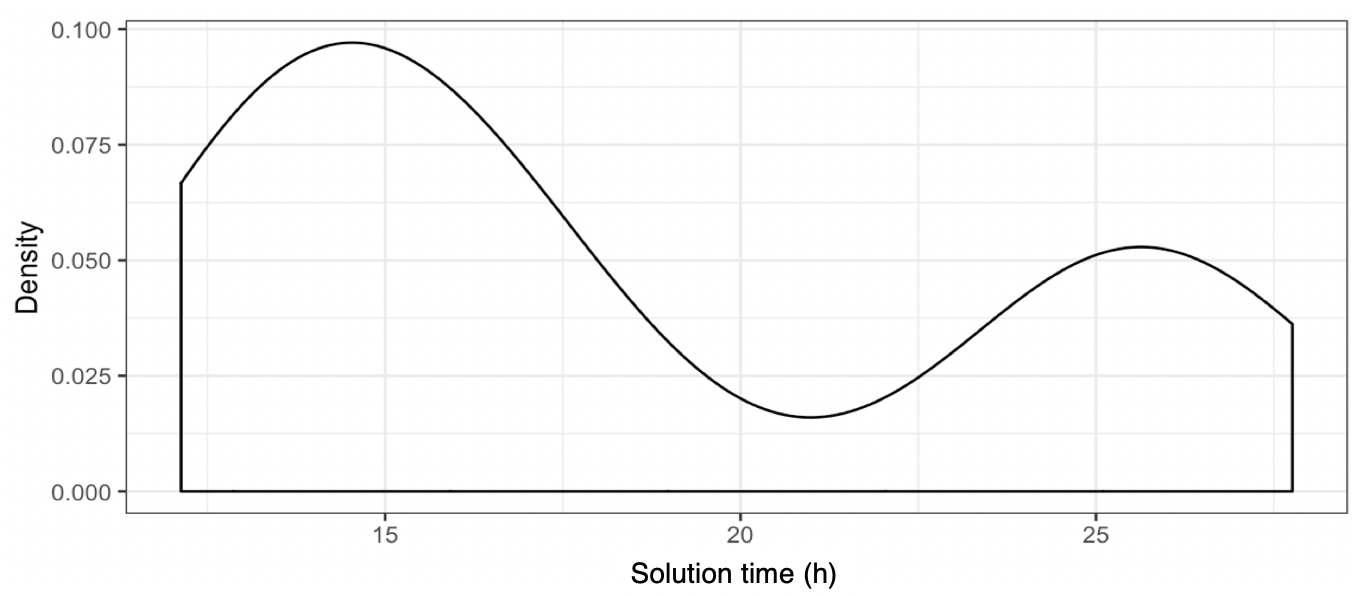}
    \caption{\footnotesize BF-CD single-trajectory solution times across 50 trials.}
    \label{F:bfcd}
\end{figure}

\subsection{Approximate Pareto Frontiers.}\label{App:pareto}

{\color{black} As discussed in Section \ref{S:computational_analysis}, Pareto-efficient curves cannot be obtained using the SOS2-CD approach. Instead, we now approximate a Pareto-efficient curve with incumbent solutions generated throughout SOS2-CD for the solutions in Table \ref{T:benchmark}. We compile the PX, PR, and VM performance metrics associated with each incumbent solution and plot them in Figure \ref{F:approximate_pareto}, together with the final solutions in each regime.  }

\begin{figure}[htbp!]
    \centering
    \includegraphics[width=0.5\textwidth]{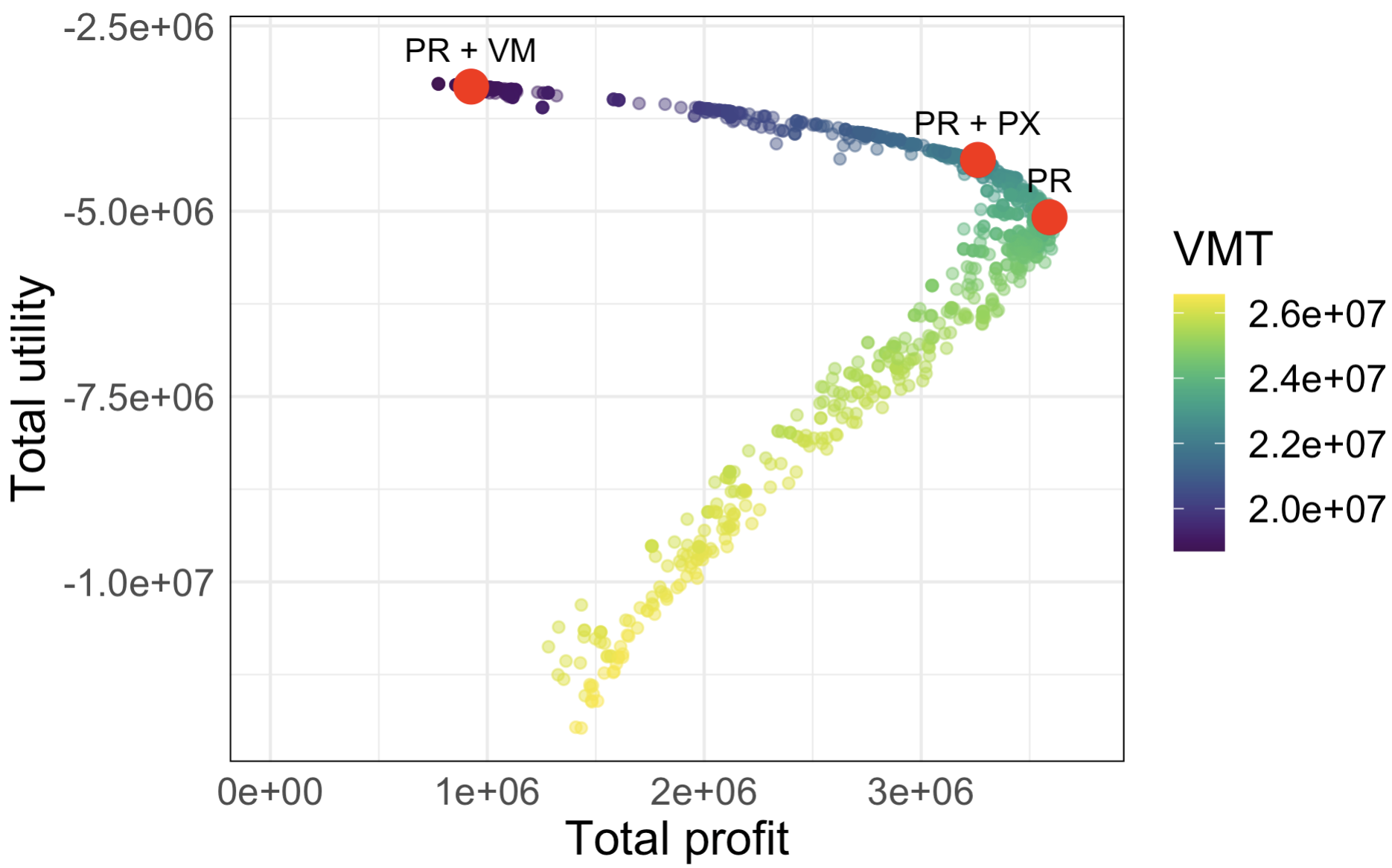}
    \caption{Approximate Pareto frontier: total profit vs. total utility vs. VMT. Constructed from intermediate incumbent solutions collected throughout SOS2-CD runs that generated the PR, PR+PX, and PX+VM benchmarks in Table \ref{T:benchmark}. Final solutions of each benchmark are shown in red.}
    \label{F:approximate_pareto}
\end{figure}

{\color{black} The PR solution achieves among the highest total profits of all incumbent solutions, as to be expected from the main case study in Section \ref{S:computational_analysis}. The PR+PX solution lies on the outer hull, trading-off total utility and profit. Comparing against the color legend, the PR+VM solution nearly achieves the lowest observed VMT in the batch of intermediate solutions. From the lack of very low-profit solutions, we conclude that SOS2-CD avoids exploring the min-VMT solution (corresponding to no profit) due to the profit priority weight. Though the three benchmarks were generated from a linear scalarization of the multiobjective problem rather than with Pareto-efficient optimization, they are still non-dominated by the observed intermediate solutions. }

\section{Greater Boston Area Case Study}\label{App:case_study}



\begin{figure}[h]
    \centering
    \includegraphics[width=0.4\textwidth]{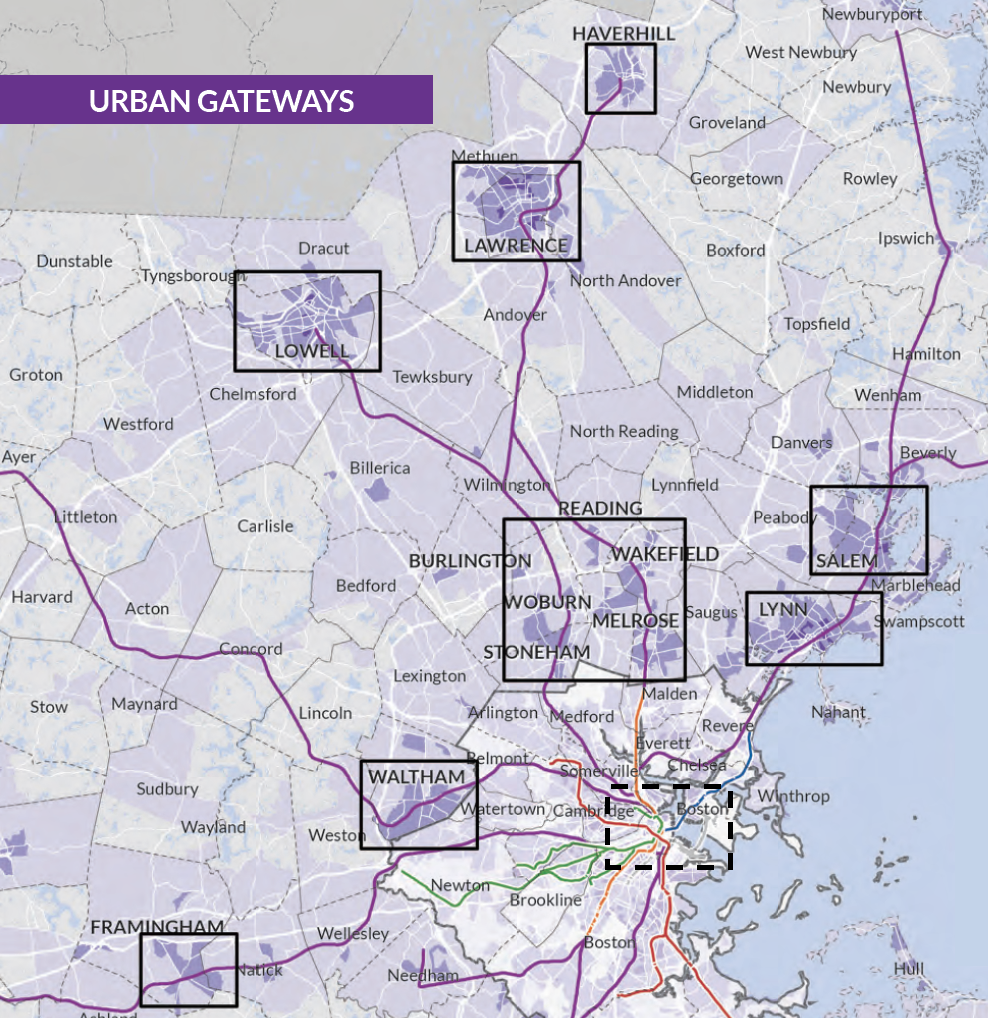}
    \caption{\footnotesize Partial map of urban gateways identified by the MBTA \citep{mbta2019}. The regions are identified with solid-line bounding boxes and capital text. The dashed-line bounding box demarcates the region we identify as the inner city. (Urban gateway not depicted: Brockton, located south of the inner city.)}
    \label{F:pplace}
\end{figure}

We consider a pricing alliance for the morning commute in the Greater Boston Area, restricting our time window to 6-10AM on a typical weekday in Fall 2018, {\color{black} and assuming negligible marginal costs}. We define a \textit{local commute} as either working in the town of residence or in an adjacent town that is also part of the service region. For example, commutes between Salem and Lynn or between Burlington and Melrose are considered local. We use the Longitudinal Employer-Household Dynamics (LODES) datasets from the U.S. Census Bureau to approximate commuting volumes for each origin/destination pair of census tracts \citep{lodes}. For the calculation of the ``Real fares'' in Table \ref{T:benchmark}, we compute the MOD base fare by combining Lyft's published minimum fare and service fee, and we compute their markup by combining the published markups per unit distance and time, assuming an average vehicle speed of 25 mph \citep{lyft}. We ignore fare multipliers utilized to manage the two-sided market, since they are outside the scope of this work; note that this may lead to slight undercounting of real MOD fares and slight overcounting of the real ridership. The real-world MBTA commuter rail base fare and markup are interpolated from its zone-based pricing structure, which assigns higher prices to farther zones \citep{mbta_fare}.

We obtain MBTA's commuter rail network data using MBTA General Transit Feed Specification (GTFS) data \citep{gtfs}, while the MOD operator corresponds to all potential direct travel options and first-mile connections in the service region. We construct route choice sets for each passenger type by first executing Yen's $k$-shortest paths algorithm \citep{yen}. We then include in each passenger type's route choice set their fastest option of each mode: transit-only, MOD-only, hybrid (MOD first mile to transit), and driving (i.e., the outside option). We represent the utility of each route option as a linear combination of in-vehicle, walking and wait time; incurred costs including fare, gasoline, and parking fees as appropriate; and mode discomfort relative to the convenience of driving. Total travel times are the sum of in-vehicle travel time, expected wait time, and walking time. Transit wait times are computed using trip frequencies from GTFS data, assuming uniformly distributed departures over the model’s time horizon. MOD wait times are uniformly distributed in a 5-15 minute interval. The discount activation categories correspond to town pairs. In other words, a discount might be activated for any MOD or hybrid-mode trip leg originating from a town in the service region to the inner city, to an adjacent town that is also part of the service region, or to another destination in the same town. There are 77 discount activation categories in the case study. We allow each operator to set fares up to a maximum of \$10 for base fares, \$5/mile for markups, and 0.5 for discount multipliers.

\color{black}

\section{Additional Practical Insights}\label{app:practical}

\subsection{3-Objective Regimes.}\label{app:objTableExtension}
\vspace{-5mm}
\begin{table}[htbp!]
    \centering
    \scriptsize
        \caption{Aggregate metrics for regimes prioritizing all objectives. Same nomenclature as in Table \ref{T:testsuite}. }
    \label{T:testsuite_extended}
    \begin{tabular}{ccc|ccc|ccc|c}
    \toprule[1pt]
    \multicolumn{3}{c|}{Objective weights} & \multicolumn{3}{c|}{Route price (\$)} & \multicolumn{3}{c|}{Performance metrics} & System  \\    
    $\pi^{PX}$ & $\pi^{PR}$ & $\pi^{VM}$ & Min. & Mean & Max.  & PX  & PR  & VM  & util. \%  \\\hline 
0.2&	1.0&	0.2&	\$7.77&	\$17.02&	\$38.61&	52.91\%&	94.58\%&	124.48\%&	30.70\%\\
0.4&	1.0&	0.4&	\$6.00&	\$13.89&	\$26.57&	67.40\%&	81.11\%&	117.92\%&	33.01\%\\
0.6&	1.0&	0.6&	\$5.50&	\$10.63&	\$18.32&	79.85\%&	59.00\%&	111.10\%&	35.73\%\\
0.8&	1.0&	0.8&	\$2.16&	\$8.08&	\$12.42&	89.72\%&	31.73\%&	104.84\%&	38.11\%\\
1.0&	1.0&	1.0&	\$0.00&	\$5.66&	\$8.44&	96.23\%&	9.45\%&	100.70\%&	41.02\%\\
    \bottomrule[1pt]
    \end{tabular}
\end{table}
\vspace{-8mm}

\subsection{Incorporating Relationship Between Travel Time, Fleet Size, and Passenger Volume.}\label{A:wait_time}

In an MOD network, average travel time increases with the volume of passengers who choose to travel on that mode---a higher number of MOD and hybrid route passengers leads to additional detours and wait times for the MOD and hybrid routes. Further, average travel time can be reduced by increasing the MOD fleet size. Finally, longer travel times reduce passengers' travel utility making the MOD and hybrid alternatives less attractive, which in turn reduces travel demand for these routes. Each of these dynamics is complex and very difficult to model. Fortunately, some prior studies, such as \cite{do2019}, have developed approximate versions of these models, particularly suitable for a design study such as ours. Therefore, in this section, we will demonstrate how we leveraged the work by \cite{do2019} to first right-size the fleet, and then to evaluate the impact of passenger volume changes on the average travel times.

Specifically, \cite{do2019} (DO, for short) have developed an analytical approach for approximating travel time, as a function of passenger volume and fleet size. We first apply the DO approach to the non-cooperative setting prevalent in the real-world case study to calibrate fleet size. We pull non-cooperative fares from Figure \ref{F:noncoop_fares} and record passenger arrival rates to the MOD subnetwork in each contiguous cluster of urban gateways under non-cooperative Nash equilibrium fares. We assume a detour ratio of 1.2 (i.e., average door-to-door time equal to 120\% of direct driving time) to compute fleet sizes using the DO approach. Then using these fleet sizes and the DO approach, we recalculate the alliance objective values. We find that the new objective values have a 0.05\% gap with the original ones, confirming that the model is well-calibrated.

Now, holding this fleet size constant, for the allied setting in the profit maximization regime, we iteratively solve a fixed-point problem to ensure that the resulting optimal fare parameters and discount activations, passenger flow volumes, and travel times are in equilibrium with each other.
We begin the iterations using the optimal fares from Figure \ref{F:testsuite_fares} and the corresponding passenger flows for each urban gateway cluster. We use the previously calibrated fleet sizes to calculate the corresponding detours and travel times. We then retrieve updated passenger flows by re-solving the fare-setting model with the new passenger utilities. We repeated this process until convergence, as depicted in Table \ref{T:wait_detour_table}.

\begin{table}[htbp!]
    \centering
    \scriptsize
        \caption{Convergence of detours and passenger volumes on MOD subnetwork.$^*$ ``Large'' stands for Burlington +  Melrose + Reading + Stoneham + Wakefield + Woburn. }     
    \resizebox{\textwidth}{!}{\begin{tabular}{l|ccc|cc|cc|cc|cc}
    \toprule[1pt] 
    & \multicolumn{3}{c|}{Non-cooperative calibration}&\multicolumn{2}{c|}{Alliance Iter. 0} &		\multicolumn{2}{c|}{Alliance Iter. 1} &	\multicolumn{2}{c|}{Alliance  Iter. 2} &		\multicolumn{2}{c}{Alliance Iter. 3} \\
\hline
Urban gateway cluster & Volume & Fleet size & Detour &	Volume	&Detour	&Volume	&Detour	&Volume	&Detour	&Volume	&Detour \\ \hline
Large$^*$	      &7,817 &	388	&1.2 &7,355	&1.54	&7,369	&1.12	&7,456	&1.13	&7,451	&1.13 \\
Lawrence	&1,722	&55	&1.2 &1,582	&1.73	&1,581	&1.19	&1,654	&1.20	&1,653	&1.19 \\
Salem + Lynn	&5,454	&219	&1.2 &5,156	&1.49	&5,157	&1.15	&5,386	&1.18	&5,366	&1.18 \\
Framingham	&1,557	&76	&1.2 &1,417	&1.53	&1,414	&1.18	&1,482	&1.19	&1,480	&1.19 \\
Haverhill	&1,189	&70	&1.2 &1,041	&1.67	&1,041	&1.18	&1,098	&1.19	&1,097	&1.19\\
Winchester	&871	&39	&1.2 &633	&1.74	&631	&1.19	&667	&1.19	&667	&1.19\\
Lowell	     &2,183	&77	&1.2 &1,966	&1.68	&1,965	&1.18	&2,063	&1.19	&2,061	&1.19\\
Waltham	    &2,513	&84	&1.2 &2,500	&1.62	&2,502	&1.20	&2,031	&1.16	&2,038	&1.16\\
Brockton	&2,994	&112	&1.2 &2,677	&1.58	&2,673	&1.17	&2,807	&1.18	&2,803	&1.18\\
         \bottomrule[1pt]
    \end{tabular}}
    \label{T:wait_detour_table}
\end{table}

Finally, we compared the results across the fixed-point iterations to each other and with those obtained in Section \ref{S:results}. The iterative process converged in 3 iterations. In general, the passenger volumes are highly stable. Equilibrium passenger flows in the MOD subnetwork increased by 1.19\% over the original values, and total profit increased by 1.30\%. The average final detour was 117.6\% (i.e., average passenger door-to-door travel times were 17.6\% longer than driving directly)---a realistic level of service. Overall, we conclude that the effects of this nonlinear relationship between passenger volume, fleet size, and travel time are negligible for pricing alliance design purposes. Thus our key findings from the main body of the paper continue to hold.

\vspace{-8mm}

\subsection{Incorporating Marginal Costs.}\label{A:marginal_costs}
First, we explain how an operator-specific marginal cost can be calculated from data. Then, we compare the feasible region and objective level curves obtained using these marginal costs against those under zero marginal costs. Finally, we perform a sensitivity analysis of marginal costs.  

Table \ref{T:costData} summarizes average miles per gallon (MPG), fuel prices, and assumed vehicle occupancy for different transportation modes. Based on these values, the marginal cost can be computed as 
\begin{equation}
    \text{Marginal cost per mile per passenger} = \frac{\text{Cost of gasoline per gallon}}{\text{(Miles per gallon)} \cdot \text{(Passengers per vehicle)}}.
\end{equation}
\vspace{-8mm}

\begin{table}[htbp!]
    \centering
    \caption{Marginal cost computation for different transit modes. $^*$Fuel cost data from 2022 \citep{fuelData}. $^{**}$Fuel mileage data \citep{mpgData}.}
    \resizebox{\textwidth}{!}{\begin{tabular}{llcccc}
    \toprule[1pt]
       Mode & Vehicle  & Vehicle  & Cost of fuel  & Miles per  & Marginal cost per \\
       & type & occupancy & per gallon &  gallon &  mile per passenger ($c_k$) \\\hline
       Transit & Transit bus & 10 & \$4.99$^*$  & 3.26$^{**}$ & \$0.15 \\ 
        High-capacity ride-share  & Shuttle & 4 & \$3.95$^*$  & 7.10$^{**}$ & \$0.14 \\ 
        Single-occupancy ride-share & Sedan & 1 & \$3.95$^*$ & 25.5$^{**}$ & \$0.15 \\ 
       \bottomrule[1pt]
    \end{tabular}}
    \label{T:costData}
\end{table}

The marginal costs in Table \ref{T:costData} are all very close to each. Setting $c_k = \$0.15$ for each operator, we compare the feasible regions and level curves of the synthetic example \eqref{E:synthetic} with no marginal costs in Figure \ref{F:costObj}.

\begin{figure}[htbp!]
    \centering
    \includegraphics[width=0.45\textwidth]{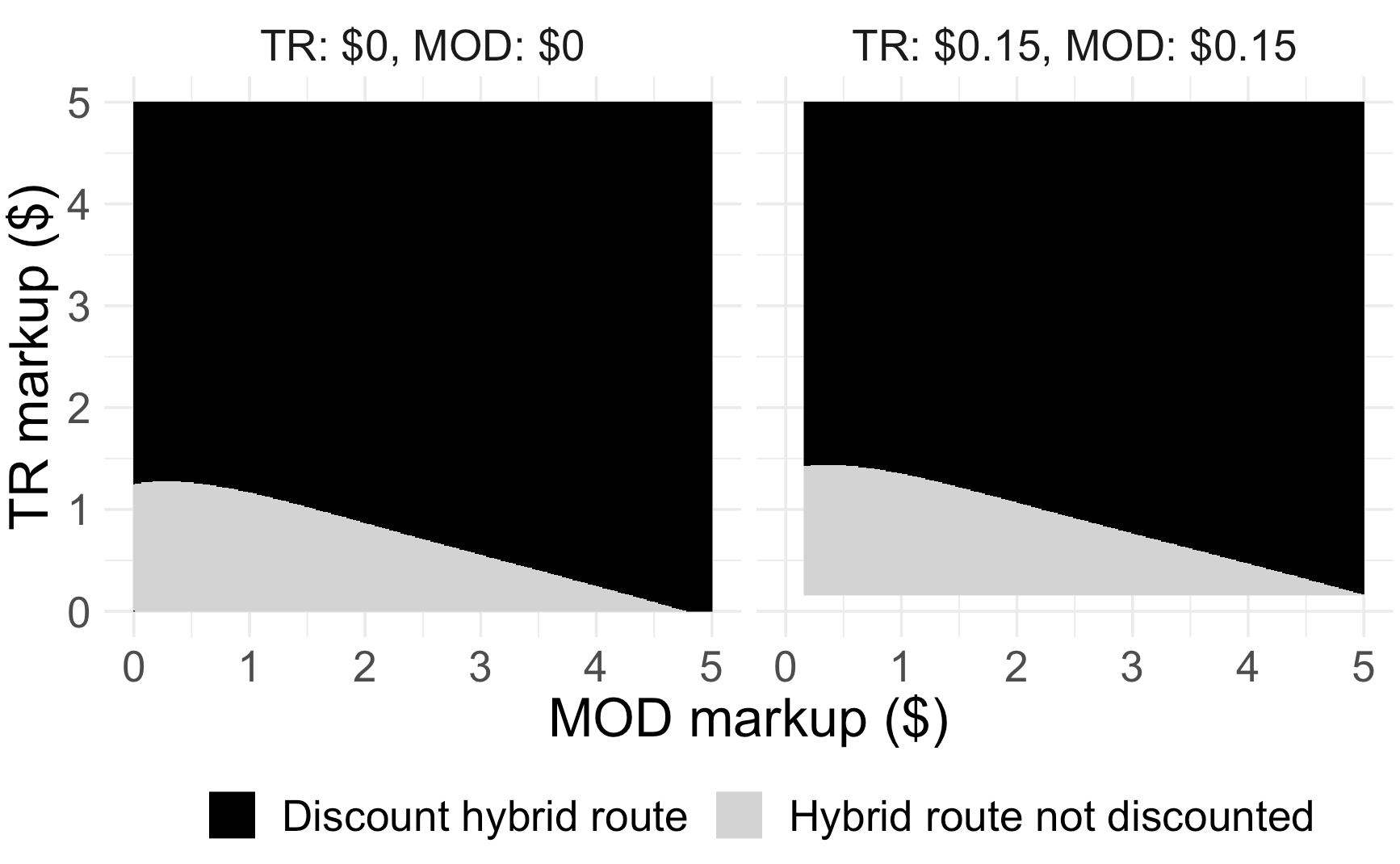}
    \includegraphics[width=0.45\textwidth]{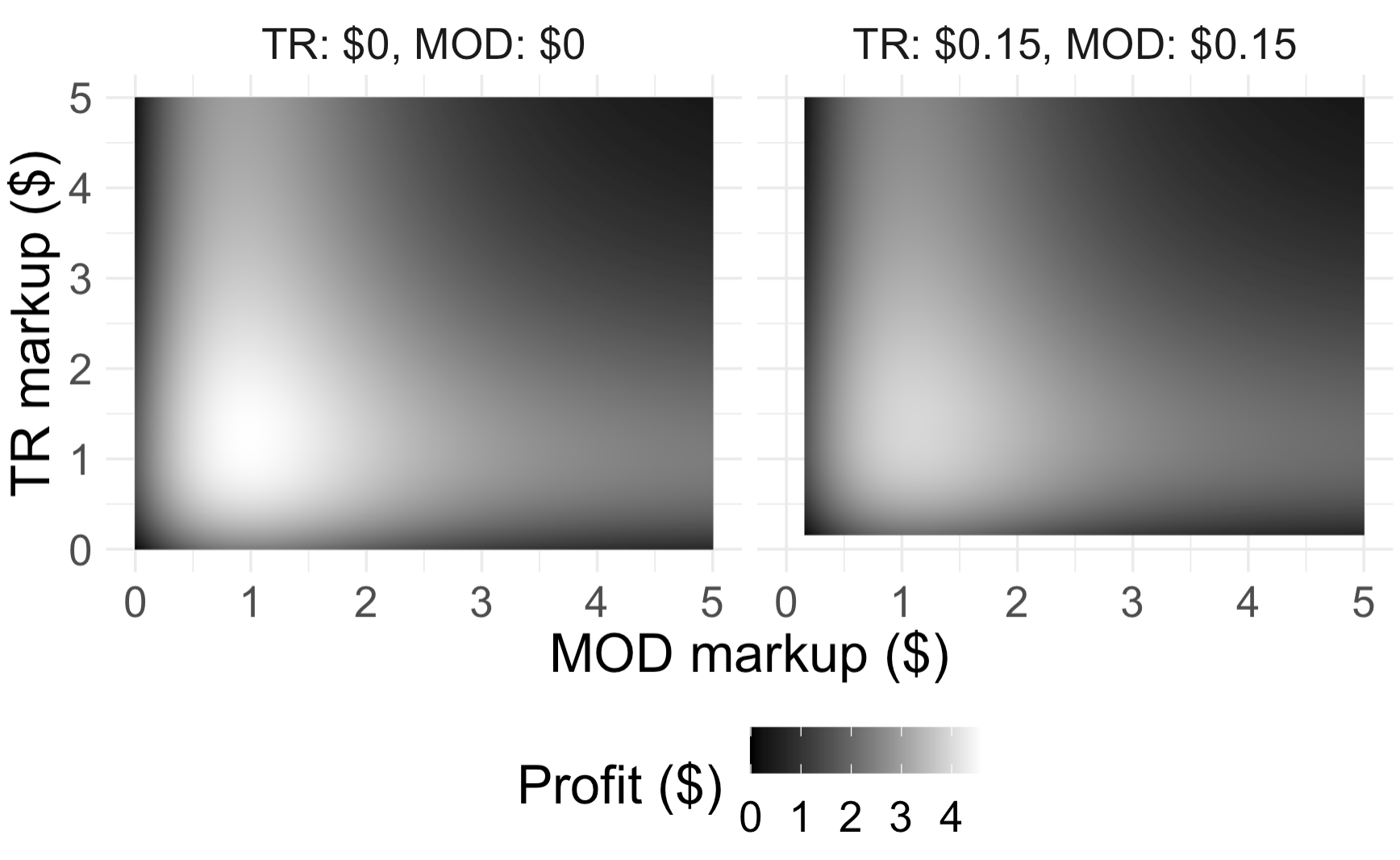}
    \caption{Optimal discounts and profits for synthetic example \eqref{E:synthetic} with $c_k = \$0$ (left) and $c_k = \$0.15$ (right).}
    \label{F:costObj}
\end{figure}

Predictably, the optimal discount decision retains a similar shape as before: the hybrid route is discounted when the transit markup is moderate or when both operators’ markups are high. The state transition boundary has simply pushed inward (i.e., toward top right) with the non-zero marginal cost. The objective level curves also have a similar shape, with optimal profit obtained for moderately low markups. As to be expected, when marginal costs are non-zero, the optimal profit is lower (\$3.95 vs. \$4.59), and the optimal markups are higher ($\mu_{MOD} = \$1.12$ and $\mu_{TR} = \$1.40$, vs. $\mu_{MOD} = \$0.99$ and $
\mu_{TR} = \$1.27$). 

Next, we examine a broader range (\$0.10 to \$0.20 at \$0.02 increment) of marginal costs for each operator. We find that each operator's optimal markup increases linearly with their marginal cost and decreases slightly with the other operator's marginal cost. Running a simple linear regression on each operator's markup, we obtain with $R^2 = 0.99$ in each case: $\mu_{TR} = \$1.25 - \$0.07 \cdot c_{MOD} + \$1.00 \cdot c_{TR}$, and $\mu_{MOD} = \$0.99 + \$0.90 \cdot c_{MOD} - \$0.05 \cdot c_{TR}$. If we denote the optimal markup under $c_{TR} = c_{MOD} = \$0$ by $\mu_{TR}^0=\$1.27$ and $\mu_{MOD}^0=\$0.99,$ we can approximate each optimal markup by $\mu_{TR} \approx \mu_{TR}^0 + c_{TR}$ and $\mu_{MOD} \approx \mu_{MOD}^0 + c_{MOD}.$ 

The incorporation of marginal costs predictably increases prices, and as to be expected, system utilization and profits both decrease in response. However, the general outcomes are comparable to scenarios with zero marginal costs. In particular, the feasible region and objective level curves are analogously shaped. Further, optimal markups can be coarsely approximated by adding marginal costs to optimal markups under zero marginal cost, within a reasonable range of values. The significance of this conclusion is that optimal markups vary linearly with marginal costs---a highly interpretable relationship. By omitting marginal costs in broader experiments, we can observe pricing alliance design outcomes in their most ``noiseless'' form. We conclude that the practical insights from our experiments still hold with the incorporation of marginal costs.

\subsection{Implementing a Multi-Period Model.}\label{app:timePeriods}

We now investigate an extended case study with three time periods. We divide the overall time horizon into three equal parts, with 50\% of each passenger type traveling in the middle period, 25\% in the first, and 25\% in the third. We assume that the driving and MOD options incur a 20\% extra travel time due to peak traffic congestion in the middle period, with no travel time change in the first and last periods. Finally, we assume a profit-oriented alliance. Implementing the model for multiple time periods increases the computational burden. Average solve times increase to 1-2 minutes per second-stage model, as opposed to 5-10 seconds. 

We input optimal fares from the single-period case study and compare second-stage outcomes. The optimal profit for the single-period case study was \$3,658,161, as compared to a profit of \$3,739,441 for the single-period optimal fares in the three-period case study. This 2.2\% profit increase can be attributed to a slight increase in alliance passenger volume---0.60\% on average across towns. Interestingly, the same discount activation options were optimal in each case study. So, the slight increase in total alliance passenger volumes can be attributed to the travel time detours that were added to the middle time period. As a result, fewer people chose the more congested outside option and more passengers elected to travel via the alliance. 

Next, we assessed the impacts on optimal fares. We solved for optimal fares under three regimes, as in Table \ref{T:benchmark}: PR, PR+PX, and PR+VM. We solved each regime with the SOS2-CD-MD-R method and a 24-hour time limit. In all three regimes, we found the three-period optimal fares to be the same as those under the single-period setting. We conclude that increased time granularity enables more fine-grained analysis of passenger decisions, but leads to similar high-level strategic outcomes, indicating a single period to be appropriate for our case study. Thus, throughout the main results section, we utilize a single time period expediting solutions without compromising the accuracy of our findings and the validity of our insights.

\subsection{Investigating the Impacts of Demand Uncertainty.}\label{app:stochastic}

We assess the impacts of perturbed demand on optimal profits for the PADP-FS. Recall from Appendix \ref{App:case_study} that passenger volumes on each OD pair were calibrated using commuter census data. We now set $\tilde{N}_{it} \sim \big(\text{Normal}(N_{it}, \tfrac{N_{it}}{6})\big)^+$ for each $i \in \mathcal{N}, t \in \mathcal{T}$ (still with a single time period, i.e., $|\mathcal{T}|=1$) and solved the corresponding PADP-FS with SOS2-CD-MD-R and a 12-hour time limit, repeating this experiment five times. For each experiment, we retrieved the optimal objective of a profit-oriented alliance and measured the extent to which the fares optimized under perturbed demand information were suboptimal. Across the five repetitions of the experiment, we observed an average 0.024\% gap between the optimal profits obtained under the accurate knowledge of the true demand and those obtained with fares optimized assuming perturbed demand but evaluated under the true demand. We conclude from these experiments that the final fares are highly robust to small variations in demand.

\subsection{Coalition Formation with Three Operators.}\label{S:coalition}

The service regions of competing MOD operators often overlap. We examine the effects of a second MOD operator on pricing alliance design, in terms of optimal prices, passenger decisions, and operator outcomes. In particular, we are interested in investigating whether some alliances are better than others, and how the excluded MOD operator is affected by cooperation. Specifically, we add a second MOD operator to the synthetic example \eqref{E:synthetic}, and assume all three operators to be profit-oriented. The passengers' choice set includes driving, transit-only, MOD1-only, MOD2-only, transit-MOD1 hybrid and transit-MOD2 hybrid. A candidate alliance can discount their collectively operated hybrid route, but not the other hybrid route (if it exists). We assume that the MOD operators are not interested in collaborating with each other (which excludes the grand coalition---three-operator pricing alliances are uncommon in practice). Fixing MOD1's alternative specific constant utility (ASC), we tune MOD2's ASC to be either twice as negative (worse), equivalent (substitutable), or half as negative (better). As a benchmark, we also include the scenario where MOD2 does not exist (setting $ASC_2$ to $-\infty$). Table \ref{T:coalitions} catalogs optimal markups, whether a 25\% discount is applied to the allied hybrid route, the profits of each coalition or operator, the total utilization (i.e., the percentage of passengers who do not drive themselves), and average prices across all route options.

\begin{table}[htbp!]
    \centering
    \resizebox{\textwidth}{!}{\begin{tabular}{lc|cccc|c|c|c|cc}
    \toprule[1pt]
         Coalition & MOD ASCs & Transit  & MOD 1  & MOD 2  & Discount  & MOD 1  & Transit  & MOD 2  & Utilization & Average  \\
         & & markup & markup & markup & & profit & profit & profit & & price\\\hline 
Transit + MOD1 & $ASC_2 \rightarrow -\infty$       & \$1.27 & \$0.99 & --- & 25\% & \multicolumn{2}{c|}{\$4.59} & --- & 18.63\% & \$24.30 \\\cline{7-9}
None          &                         & \$1.11 & \$0.94 & --- & ---  & \$2.58 & \$2.00            & --- & 19.57\% & \$24.20 \\ \hline
Transit + MOD1 & $ASC_2 = ASC_1 \cdot 2$ & \$1.26 & \$0.99 & \$0.87 & 25\% & \multicolumn{2}{c|}{\$4.55} & \$1.21 & 23.63\% & \$24.77 \\\cline{7-9}
Transit + MOD2 &                         & \$1.25 & \$0.93 & \$0.93 & 25\% & \$2.43 & \multicolumn{2}{c|}{\$3.27} & 24.09\% & \$23.91\\\cline{7-9}
None &                         & \$1.11 & \$0.93 & \$0.87 & --- & \$2.42 & \$2.11            & \$1.15 & 24.55\% & \$24.12 \\ \hline
Transit + MOD1 & $ASC_2 = ASC_1 $         & \$1.26 & \$0.98 & \$0.92 & 25\% & \multicolumn{2}{c|}{\$4.53} & \$2.29 & 28.10\% & \$25.02 \\\cline{7-9}
None          &                         & \$1.11 & \$0.92 & \$0.92 & --- & \$2.28 & \$2.23 & \$2.28 & 29.03\% & \$24.36 \\\hline
Transit + MOD1 & $ASC_2 = ASC_1 / 2$     & \$1.25 & \$0.98 & \$0.96 & 25\% & \multicolumn{2}{c|}{\$4.53} & \$3.17 & 31.26\% & \$25.14\\\cline{7-9}
Transit + MOD2 &                         & \$1.32 & \$0.91 & \$1.02 & 25\% & \$2.20 & \multicolumn{2}{c|}{\$5.50} & 30.97\% & \$25.85 \\\cline{7-9}
None          &                         & \$1.12 & \$0.92 & \$0.96 & --- & \$2.18 & \$2.33 & \$3.15 & 32.06\% & \$24.72\\
\bottomrule[1pt]
    \end{tabular}}
    \caption{Coalition solutions for the synthetic example \eqref{E:synthetic} with an additional MOD operator.}
    \label{T:coalitions}
\end{table}

All operators benefit from horizontal cooperation, even if they do not participate.
As to be expected, coalitions earn strictly more than they would by competing, which is in alignment with Lemma \ref{L:profit_allocation}c. Despite lower system utilization, alliances increase total profits by raising prices and by collecting more earnings per passenger (except for sub-optimal coalitions---when MOD2 is ``worse'' than MOD1, average route prices for the transit-MOD2 alliance are lower than in the non-cooperative scenario). When a coalition is formed, the external MOD operator still benefits over the fully non-cooperative scenario---i.e., they earn more when the other two operators are allied than when all three operators compete independently. The alliance sets higher prices together than when those same operators work independently. In response, the non-allied MOD operator attracts passengers who leave allied operators after their prices increase. We would also expect external operators to benefit if the alliance were altruistically oriented---if the alliance lowers their prices to prioritize passengers or VMT, then the profit-oriented external MOD operator benefits from the additional volume on hybrid mode routes due to the increased room in passenger budgets. 

Further, optimal pricing alliance design itself is negligibly influenced by the existence of an external MOD operator, regardless of their level of dominance in the market. To observe this, note that the ASC of MOD1 is fixed throughout all experiments. When the transit operator aligns with MOD1, their outcomes are fairly stable: markups vary by \$0.02 at most, and the range of total alliance profits is \$0.06 ($\sim$1\% difference). Utilization varies a lot based on MOD2 existence or dominance, but optimal prices do not. Importantly, the transit operator prefers to collaborate with the dominant, rather than the non-dominant, MOD operator when they have the opportunity to choose.

\color{black}

\end{APPENDICES}

\end{document}